\documentclass[12pt]{amsart}

\usepackage{amssymb, amsthm, enumerate, comment, color, mathrsfs}
\usepackage{amsmath, mathtools}
\newtheorem{theorem}{Theorem}[section]
\newtheorem{lemma}[theorem]{Lemma}
\newtheorem{proposition}[theorem]{Proposition}
\newtheorem{corollary}[theorem]{Corollary}
\theoremstyle{definition}
\newtheorem{ex}[theorem]{Example}
\theoremstyle{remark}
\newtheorem{remark}[theorem]{Remark}

\usepackage{geometry}

\newcommand{\lz}{\left(}
\newcommand{\pz}{\right)}

\newcommand{\E}{\mathcal{E}}

\renewcommand{\epsilon}{\varepsilon}
\newcommand{\C}{\mathbb{C}}

\newcommand{\sumstar}{\mathop{\sum\nolimits^{\mathrlap{*}}}}
\newcommand{\sumplus}{\mathop{\sum\nolimits^{\mathrlap{+}}}}
\newcommand{\N}{\mathbb{N}}
\newcommand{\R}{\mathbb{R}}

\newcommand{\lab}{\left\lvert}
\newcommand{\rab}{\right\rvert}
\newcommand{\M}{\mathcal{M}}

\newcommand{\re}{\mathrm{Re}}

\newcommand{\im}{\mathrm{Im}}

\newcommand{\bfrac}[2]{\lz\frac{#1}{#2}\pz}

\renewcommand{\mod}[1]{\text{ (mod $#1$)}}

\renewcommand{\phi}{\varphi}

\newcommand{\norm}[1]{\left\|#1\right\|}

\newcommand{\sumd}{\sideset{}{^d}\sum}

\numberwithin{equation}{section}

\title{On optimality of mollifiers}

\author{Martin \v Cech}
\address{Charles University, Faculty of Mathematics and Physics, Department of Mathematical analysis and Department of Algebra, Sokolovská 83, 186 75 Praha 8, Czech Republic}
\email{martin.cech@matfyz.cuni.cz}

\author{Kaisa Matom\"aki}
\address{Department of Mathematics and Statistics, University of Turku, 20014 Turku, Finland}
\email{ksmato@utu.fi}

\begin{document}
	\begin{abstract}
		Mollifiers are used in a variety of contexts, for instance to study the non-vanishing of $L$-functions. In this paper, we study the general question of finding optimal mollifiers and provide criteria to identify them provided the corresponding mollified moments can be computed. As an application, we study the non-vanishing of central values of Dirichlet $L$-functions. In particular we show that the Michel-Vanderkam mollifier is optimal in a wide class of balanced two-piece mollifiers as well as provide a new proof that the Iwaniec-Sarnak mollifier is optimal in a wide class of one-piece mollifiers.
	\end{abstract}
	
	\maketitle
	\section{Introduction}
	
	In the methods used to study the non-vanishing or zeroes of interesting objects, it is important to ``dampen'' the large values of the objects using a suitable mollifier. Examples of such situations include the study of non-vanishing of central values of Dirichlet $L$-function $L(1/2, \chi)$ and the study of critical zeroes of the Riemann zeta function $\zeta(s)$ (i.e. zeroes on the critical line $\re(s) = 1/2$).
	
	In particular in the context of the zeta zeroes, more and more complicated mollifiers have been introduced (see e.g.~\cite[Section 1.3]{PRZZ} for discussion of different choices), leading to small improvements on the proportion of critical zeroes.
	
	Except for the simplest cases, the optimal choice of mollifier is not known. The purpose of this work is to study some aspects of the optimality question. We work in the non-vanishing setting, but similar ideas work also in the context of the critical zeroes, see Remark~\ref{remark:zetazeroes} for more information.
	
	For simplicity, we postpone discussing our general results to Section~\ref{sec:general mollifying}; they include for instance Corollary~\ref{cor:MNMMcomp}, where we describe a simple criterion to find out in a very general setting whether adding a multiple of an extra mollifier to an existing one is useful, provided the corresponding moments can be calculated. For example, we will find that adding a new piece is almost always useful if it is not too correlated with the original mollifier. We hope our general results will be useful in various contexts. In this paper we apply them to prove optimality of existing mollifiers, but in the future one might be able to utilize similar ideas for instance to find better mollifiers.
	
	For now we focus on the non-vanishing of central values of Dirichlet $L$-functions associated to primitive characters modulo $q \in \mathbb{N}$ and the consequences of our general results in this case. In this context, Chowla~\cite{Ch} conjectured that $L(1/2,\chi)\neq 0$ for every Dirichlet character $\chi$. Partial results can be obtained towards this conjecture by showing that a positive proportion of central values of Dirichlet $L$-functions in some family are non-zero.
	
	The usual strategy for proving that a positive proportion of the central values $L(1/2,\chi)$ are non-zero is to apply the Cauchy-Schwarz inequality, together with the knowledge of the first two moments. Indeed, we have
	\begin{equation*}
		\lz \sumstar_{\chi\mod q} L(1/2,\chi)\pz^2\leq \sumstar_{\chi\mod q} \mathbf{1}_{L(1/2,\chi)\neq 0}\cdot\sumstar_{\chi\mod q}|L(1/2,\chi)|^2,
	\end{equation*}
	which yields the lower bound
	\begin{equation*}
		\sumstar_{\chi\mod q} \mathbf{1}_{L(1/2,\chi)\neq 0}\geq \frac{\lz \sumstar\limits_{\chi\mod q} L(1/2,\chi)\pz^2}{\sumstar\limits_{\chi\mod q}|L(1/2,\chi)|^2}.
	\end{equation*}
	This by itself is not enough to yield a positive proportion, as the second moment is too large because it is dominated by a few large central values $L(1/2,\chi)$. 
	
	To deal with this problem, one notices that, for any $M(\chi)$,
	\[
	L(1/2,\chi)M(\chi) \neq 0\implies L(1/2,\chi)\neq 0,
	\]
	so one can choose a suitable $M(\chi)$ and study the non-vanishing of $L(1/2, \chi) M(\chi)$. This is useful as one can hope to choose $M(\chi)\approx L(1/2,\chi)^{-1}$; such an $M(\chi)$ designed to dampen the large values of $L(1/2,\chi)$ is called a mollifier. The Cauchy-Schwarz inequality now gives
	\begin{equation*}
		\begin{aligned}
			\sumstar_{\chi\mod q} \mathbf{1}_{L(1/2,\chi)\neq 0}&\geq \sumstar\limits_{\chi\mod q} \mathbf{1}_{L(1/2,\chi)M(\chi)\neq 0} \geq \frac{\lz\sumstar\limits_{\chi\mod q} L(1/2,\chi)M(\chi)\pz^2}{\sumstar\limits_{\chi\mod q}|L(1/2,\chi)M(\chi)|^2},
		\end{aligned}
	\end{equation*}
	and one may now hope to obtain a positive proportion of non-vanishing provided one can find an efficient mollifier and compute the first two mollified moments. A very natural choice for a mollifier is
	\begin{equation}\label{eqn: one-piece mollifier}
		M(\chi)=\sum_{b \leq q^\theta} \frac{x_b\chi(b)}{\sqrt b}
	\end{equation} for some coefficients $x_b \in \mathbb{R}$, with $\theta$ as large as possible.
	
	For the family of Dirichlet characters modulo $q$, this strategy was successfully pursued by Iwaniec and Sarnak~\cite{IS}. It turns out that the optimal choice of $x_b$ in~\eqref{eqn: one-piece mollifier} is essentially $x_b=\mu(b)(1-\frac{\log b}{\log q^\theta})$ (actually Iwaniec and Sarnak had an asymptotically equivalent choice, see~\cite[formulas (6.17)--(6.18)]{IS}, noting also their differently normalized definition of $M(\chi)$; nonetheless, we will call the mollifier with $x_b=\mu(b)(1-\frac{\log b}{\log q^\theta})$ the Iwaniec-Sarnak mollfier). Iwaniec and Sarnak were able to take any $\theta < 1/2$, which led to the proportion of non-vanishing at least $1/3-\varepsilon$ for any $\varepsilon > 0$.

	Here we will deduce from our general results and moment calculations the optimality of the choice of Iwaniec and Sarnak~\cite{IS} as they do not give a rigorous proof of the optimality in their paper but rather mention it is possible to prove it, in a "reasonable" class of mollifiers (see Remark~\ref{rem:IS} below for more information about their work, and Remark~\ref{rem:previousOptimality} below for more general information about optimality proofs). 
	
 For a mollifier $M$, we write $\beta_q(M)$ for the non-vanishing proportion obtained from it, i.e.
	\begin{equation*}
		\beta_q(M):=\frac{\lz\sumstar\limits_{\chi\mod q} L(1/2,\chi)M(\chi)\pz^2}{\sumstar\limits_{\chi\mod q}|L(1/2,\chi)M(\chi)|^2}
	\end{equation*}
	(when the denominator vanishes, we define $\beta_q(M) = 0$), and will prove the following theorem.
	
	\begin{theorem}\label{thm:ISOptimality}
		Let $\theta < 1/2$, let $q \geq 3$, and let
		\begin{equation*}
			M_{\mathrm{IS}}(\chi):=\sum_{b\leq q^{\theta}}\frac{\mu(b)\chi(b)}{\sqrt b}\left(1-\frac{\log b}{\log q^\theta}\right)
		\end{equation*} 
		be the Iwaniec-Sarnak mollifier. Let 
		\begin{equation*}
			N_{\mathrm{G}}(\chi):=\sum_{b \leq q^\theta} \frac{x_{b} \chi(b) }{\sqrt{b}}
		\end{equation*} 
		be a general mollifier of the same type, where $x_b$ are arbitrary complex coefficients normalized so that $x_1 = 1$. Assume that $x_b \ll q^\varepsilon$ for every $\varepsilon > 0$ and $b \in \mathbb{N}$. Then
		\begin{equation} \label{eq:I-Soptimalitybound}
			\beta_q(N_{\mathrm{G}})\leq \beta_q(M_{\mathrm{IS}})+o(1) = \frac{1}{1+\frac{1}{\theta}} + o(1).
		\end{equation}
	\end{theorem}
	\begin{remark} \label{rem:previousOptimality}
	It seems that proofs of optimality of mollifiers are very sparse in the
literature. The only rigorous optimality proof we are aware of is a different proof of the optimality of the Iwaniec-Sarnak type mollifier due to Soundararajan which can be found from Radziwi{\l\l}'s unpublished manuscript (see~\cite[Proposition B]{Radziwill}, which is stated in the context of the Riemann zeta function but the same argument, combined with the moment calculations from~\cite{IS}, yields the optimality in our context). However, the proof utilizes the exact form of the mollified moments, while our proof should be more adaptable to different settings.

In other cases, one often pre-determines the coefficients of the mollifier, at least apart from optimizing polynomials, instead of looking for the optimal choice. This was done for example by Bui~\cite{Bui} who worked with the mollifier~\eqref{eq:BuiMollifier} below, and by Michel and Vanderkam in \cite{M-V}, who used the mollifier \eqref{eq:MVblanced} below (where a polynomial has been optimized).
		
		Alternatively one picks up a convenient part of the second moment and finds the minimal mollifer for that, see for instance the discussion in~\cite[Below the proof of Proposition 6 in Section 2.4.7]{KM2000}. Similarly, in Soundararajan's celebrated paper \cite{sound}, when going from the exact formula (6.6) to (6.7), only the diagonal part of the quadratic form is retained for which the optimal mollifier is found.
	\end{remark}
	
	Motivated by the approximate functional equation and Soudararajan's work~\cite{Sound95}, Michel and Vanderkam~\cite{M-V} considered the two-piece mollifier
	\begin{equation}\label{eq:MVblanced}
		\sum_{b\leq q^{\theta}}\frac{\mu(b)\chi(b)}{\sqrt b}\left(1-\frac{\log b}{\log q^\theta}\right)+\overline\epsilon_\chi\sum_{b \leq q^{\theta}}\frac{\mu(b)\overline\chi(b)}{\sqrt b}\left(1-\frac{\log b}{\log q^\theta}\right)
	\end{equation}
	that we will call the Michel-Vanderkam mollifier. The appearance of cross terms involving the root number makes computing the second moment more difficult, and they were able to take any $\theta < 1/4$, which lead to the same non-vanishing proportion $1/3-\varepsilon$ as that of Iwaniec and Sarnak. However, the two-piece mollifier showed up to be more efficient in studying the non-vanishing proportion of the central values of higher derivatives of the completed $L-$functions $\Lambda^{(k)}(1/2,\chi)$ (with $1-\frac{\log b}{\log q^\theta}$ in~\eqref{eq:MVblanced} replaced by $P_k\lz1-\frac{\log b}{\log q^\theta}\pz$   for some polynomial $P_k$, which was optimized).
	
	The proportion $1/3-\varepsilon$ was first improved to $0.3411-\varepsilon$ by Bui~\cite{Bui}, who introduced a new two-piece mollifier of the form
	\begin{equation} \label{eq:BuiMollifier}
		\sum_{b\leq q^{\theta_1}}\frac{\mu(b)\chi(b)}{\sqrt b} P_1\left(\frac{\log (q^{\theta_1}/b)}{\log q^{\theta_1}}\right) +\frac{1}{\log q}\sum_{a b \leq q^{\theta_2}}\frac{\Lambda(a)\mu(b)\overline\chi(a)\chi(b)}{\sqrt{a b}}  P_2\left(\frac{\log (q^{\theta_2}/(ab))}{\log q^{\theta_2}}\right),
	\end{equation}
	with polynomials $P_1(x), P_2(x)$ such that $P_1(0) = P_2(0) = 0$. Thanks to the avoidance of the root number, Bui was able to take any $\theta_1,\theta_2<1/2$, which led to the improvement of the non-vanishing proportion. It is not clear what is an optimal mollifier of this shape --- for instance one can obtain a (very small) improvement by replacing $\frac{\Lambda(a)}{\log q}$ by coefficients that are supported also on integers with at most $k$ prime factors for some fixed $k$.
	
	For prime moduli $q$, Khan and Ngo~\cite{Khan-Ngo} obtained the proportion $3/8-\varepsilon$ by improving the length of the Michel-Vanderkam mollifier~\eqref{eq:MVblanced} to any $\theta < 3/10$. In a subsequent paper with Mili\'cevi\'c~\cite{Khan-Mil-Ngo}, they considered the more general ``unbalanced'' two-piece mollifier
	\begin{equation} \label{eq:KMNmol}
		\sum_{b\leq q^{\theta_1}}\frac{\mu(b)\chi(b)}{\sqrt b} \left(1-\frac{\log b}{\log q^{\theta_1}}\right) + \alpha \cdot \overline\epsilon_\chi\sum_{b \leq q^{\theta_2}}\frac{\mu(b)\overline\chi(b)}{\sqrt b} \left(1-\frac{\log b}{\log q^{\theta_2}}\right)
	\end{equation} 	
	for an arbitrary $\alpha \in \mathbb{R}$ (whose optimal value was found to be $\frac{\theta_2}{\theta_1}$; we show how to obtain the same optimal value from our general results in Section \ref{se:unbalanced}), and managed to compute the moments for any lengths $\theta_1<3/8, \theta_2<1/4$, which lead to the currently best proportion $5/13-\varepsilon$. 
	
	If one averages also over $q$, one may obtain a better proportion since the extra averaging allows one to take longer mollifiers. In particular, the extra average over $q$ should allow one to take the Iwaniec-Sarnak mollifier of any length $\theta<1$, or the Michel-Vanderkam two-piece mollifier with any lengths $\theta_1, \theta_2 < 1/2$, both of which lead to proportion $1/2-\varepsilon$ (see~\cite[Section 1]{IS} and Theorem~\ref{thm:main general result DirL} below). When we average over $q\asymp Q$, we replace $q$ by $Q$ in the definitions of these mollifiers, so that the definition is independent of $q$.
	
	Pratt \cite{Pratt} claimed that, in the averaged case, one can improve the proportion to slightly above $1/2$ by adding an extra Bui piece to the Michel-Vanderkam mollifier, but there turns out to be a crucial calculation mistake in the published version of the argument (see the arXiv version of~\cite{Pratt} for up-to-date information). Assuming the Generalized Riemann Hypothesis (GRH), proportion $>0.51$ was obtained by Drappeau, Pratt and Radziwi{\l\l} \cite{DPR} by computing the one-level density rather than using mollified moments (note that their main result that extends the support of the Fourier transform of admissible test functions is unconditional, but the application to non-vanishing requires GRH).
	
	In this paper, we prove that actually in the balanced case the Michel-Vanderkam mollifier is optimal among a broad class of mollifiers, i.e., that adding any extra mollifier of a general ``Bui'' type does not lead to an improvement (and thus the strategy in the published version of~\cite{Pratt} cannot work). On the other hand, in the unbalanced case~\eqref{eq:KMNmol} with $\theta_1 \neq \theta_2$, our results imply that adding an additional Bui piece of length $\max\{q^{\theta_1},q^{\theta_2}\}$ would be helpful, assuming one can compute the mollified moments and the main terms are as expected. See Section~\ref{se:unbalanced} for more discussion on the unbalanced case.
	
	We prove the optimality result in the balanced case (Theorem~\ref{thm:main general result DirL} below) with the averaging over $q$ and for convenience insert a weight to $q$.  A completely similar result should hold also without the average over $q$, assuming that the main terms of the mollified moments are as expected, which one might hope to be provable for $\theta < 1/4$.
	
	For a mollifier $M$ we write $\beta_Q(M)$ for the (weighted) non-vanishing proportion obtained from it, i.e.
	\begin{equation*}
		\beta_Q(M):=\frac{\lz\sum\limits_{q\geq 1}\Phi\bfrac qQ \frac{q}{\phi(q)}\sumstar\limits_{\chi\mod q} L(1/2,\chi)M(\chi)\pz^2}{\sum\limits_{q\geq 1}\Phi\bfrac{q}{Q}\frac{q}{\phi(q)}\sumstar\limits_{\chi\mod q}|L(1/2,\chi)M(\chi)|^2},
	\end{equation*} where $\Phi(x)$ is a non-negative smooth (bounded and infinitely differentiable with bounded derivatives) function with $\Phi(1)\geq 1$ and support in $[1/2,2]$. When the denominator vanishes, we define $\beta_Q(M)=0$.
	
	\begin{theorem}\label{thm:main general result DirL}
		Let $\theta < 1/2$. For every $\varepsilon_1 > 0$, there exists a $c = c(\varepsilon_1) > 0$ such that the following holds. Let $Q \geq 3$ and
		\begin{equation*}
			M_{\mathrm{MV}}(\chi):=\sum_{b \leq Q^{\theta}}\frac{\mu(b)\chi(b)}{\sqrt b} \left(1-\frac{\log b}{\log Q^{\theta}}\right) +\overline\epsilon_\chi\sum_{b\leq Q^{\theta}}\frac{\mu(b)\overline\chi(b)}{\sqrt b} \left(1-\frac{\log b}{\log Q^{\theta}}\right)
		\end{equation*} be the Michel-Vanderkam mollifier.
		Let
		\begin{equation*}
			N_{\mathrm{B}}(\chi):=\sum_{ab \leq Q^\theta} \frac{x_{a, b}\overline{\chi}(a)\chi(b)}{\sqrt{ab}} + \overline\varepsilon_\chi \sum_{ab \leq Q^\theta} \frac{y_{a, b}\chi(a)\overline{\chi}(b)}{\sqrt{ab}} 
		\end{equation*} 
		be a general mollifier of Bui type, where $x_{a, b}, y_{a, b}$ are arbitrary real coefficients such that the following two conditions are satisfied for $z_{a, b} := x_{a, b} + y_{a, b}$.
		\begin{enumerate}[(a)]
			\item We have, for all $a, b \in \mathbb{N}$ and $\varepsilon > 0$,
			\[
			z_{a, b}\ll_\varepsilon Q^\varepsilon.
			\]
			\item We have
			\begin{equation*} 
			\sum_{ab \leq Q^\theta} \frac{|z_{a, b}|^2}{ab} \ll \exp\left(c \frac{(\log Q)^{3/5}}{(\log \log Q)^{1/5}}\right).
			\end{equation*}
			\item We have
			\begin{equation*}
				\sum_{q \geq 1} \Phi\left(\frac{q}{Q}\right) \frac{q \varphi^+(q)}{\varphi(q)}\sum\limits_{\substack{k^2b\leq Q^\theta \\ (kb, q) = 1}}\frac{z_{kb,k}}{kb}\asymp Q^2.
			\end{equation*}
		\end{enumerate} Then
		\begin{equation*}
			\beta_Q(N_{\mathrm{B}})\leq \beta_Q(M_{\mathrm{MV}})+\varepsilon_1 = \frac{1}{1+\frac{1}{2\theta}} + \varepsilon_1 + o(1).
		\end{equation*}
		
		If the quasi-Riemann hypothesis\footnote{By the quasi-Riemann hypothesis we mean the claim that there exists an absolute constant $\delta > 0$ such that the Riemann zeta function has no zeroes $\rho$ with $\re(\rho) \geq 1-\delta$.} holds, we do not need to assume (b).
	\end{theorem}
	
	The normalization condition (c) above is to ensure that the first mollified moment has size $\asymp 1$ (that the left-hand side corresponds to the first mollified moment follows from Proposition \ref{prop:Psi_N_1} below). In the case of an Iwaniec-Sarnak type mollifier \eqref{eqn: one-piece mollifier} the first moment is the first coefficient $x_1$, and the condition (c) in Theorem \ref{thm:main general result DirL} corresponds to the normalization $x_1=1$ in~\cite{IS} and in Theorem~\ref{thm:ISOptimality}. The conditions (a) and (b) are needed for bounding the error terms in the moment calculations.
	
	Note that in Theorem \ref{thm:main general result DirL}, $M_{\mathrm{MV}}(\chi)$ is a special case of $N_{\mathrm{B}}(\chi)$ with 
	\[ 
	x_{a,b}=y_{a,b}=\mathbf{1}_{a=1} \mu(b) \left(1-\frac{\log b}{\log Q^{\theta}}\right),
	\] 
	which satisfies the conditions (a), (b), and (c), and Bui's mollifier~\eqref{eq:BuiMollifier} corresponds to 
	\[
	x_{a,b} =\mathbf{1}_{a=1} \mu(b) P_1\left(\frac{\log (Q^\theta/b)}{\log Q^\theta}\right) + \frac{\Lambda(a)}{\log q} \mu(b) P_2\left(\frac{\log (Q^\theta/ab)}{\log Q^\theta}\right) \text{ and } y_{a, b} = 0,
	\]
	where the conditions (a) and (b) are clearly satisfied, and (c) follows from Bui's results (or Proposition~\ref{prop:Psi_N_1} below).
	
	In Section \ref{sec:general mollifying}, we state general results that enable us to compare the efficiency of given mollifiers $M,N$ in a general context and tell whether a combined mollifier $M+\alpha N$ is more efficient than the individual ones, provided one can compute the associated mollified moments. The most novel part of this is asking the relevant questions and finding the correct answers. Once the results are stated and properly digested, the proofs are elementary and rather short; the proofs of all results of Section~\ref{sec:general mollifying} will be provided in Section~\ref{sec:general mollifying proofs}.  
	
	We then apply the general results in Sections~\ref{sec:reductionIS} and~\ref{sec:reductionMV}, where we reduce the proofs of Theorems~\ref{thm:ISOptimality} and~\ref{thm:main general result DirL} to a computation of certain mollified moments. After establishing some auxiliary results in Section \ref{sec:prelim}, we adapt the methods of Iwaniec and Sarnak \cite{IS}, Michel and Vanderkam \cite{M-V} and Pratt \cite{Pratt} to compute the mollified moments with the general mollifier in Sections \ref{sec:PsiN1}--\ref{sec:PsiM2N1}. 
	
	In Section~\ref{sec:PsiM1N1}, we study the cross-term between $M_{\mathrm{MV}}$ and $N_{\mathrm{B}}$ that does not involve the twist by $\varepsilon_\chi$. We need to be very careful with the error terms --- to be able to bound them successfully, we transfer in the proof of Theorem~\ref{thm:main general result DirL} to the case where the Michel-Vanderkam mollifier is a tiny bit longer than the general Bui type mollifier (see Remark~\ref{rem:increase} below for more discussion). 
	
	In Section \ref{se:unbalanced}, we briefly discuss the unbalanced Michel-Vanderkam mollifier.
	
	\section*{Acknowledgments}
	We are grateful to Hung Bui and Roger Heath-Brown for an e-mail conversation concerning the optimality of the Iwaniec-Sarnak mollifier, thanks to which we revisited the optimality of the Iwaniec-Sarnak and Michel-Vanderkam mollifiers, which led to improved results and exposition. We would also like to thank Henryk Iwaniec and Peter Sarnak for an e-mail conversation concerning~\cite{IS}, Maksym Radziwi{\l\l} for helpful discussions concerning mollifiers, and Debmalya Basak for questions that prompted us to clarify the exposition. We also thank the anonymous referee for pointing out Radziwi{\l\l}'s paper \cite{Radziwill}.
	
	The first author was supported by the Research Council of Finland grant no. 333707 and the Charles University grants PRIMUS/24/SCI/010 and PRIMUS/25/SCI/017. The second author was supported by the Research Council of Finland grants no. 333707, 346307, and 364214.
	
	\section{General results} \label{sec:general mollifying}
	Consider a family of objects $L \colon G \to \mathbb{C}$ with some (indexing) set $G$, which depends on some parameter $Q$ that goes to infinity.
	\begin{ex} \label{ex:DirL} In the set-up of Theorem~\ref{thm:ISOptimality} we can take $G = \{\chi \mod q \colon \chi \text{ primitive}\}$ and consider $L(\chi) = L(1/2, \chi)$ for $\chi \in G$.
	\end{ex}
	
	We shall study the non-vanishing of $L(\pi)$ using the usual Cauchy-Schwarz strategy described in the previous section. In practice, one often uses weighted averages, so for $w,M, N \colon G \to \mathbb{C}$ where $w$ is a non-negative weight function (that is not identically 0) and $M,N$ are mollifiers, we write
	\begin{equation*}
		\begin{aligned}
			\Psi_M &:= \mathbb{E}_{\pi \in G}^w L(\pi) M(\pi), \\
			\Psi_{M, N} &:=\mathbb{E}_{\pi \in G}^w L(\pi) M(\pi) \overline{L(\pi) N(\pi)},
		\end{aligned}
	\end{equation*} where for $f:G\rightarrow \C$, we define
	$$
	\mathbb{E}_{\pi\in G}^w f(\pi):=\frac{\sum_{\pi\in G}w(\pi) f(\pi)}{\sum_{\pi\in G}w(\pi)}.
	$$
	By the Cauchy-Schwarz inequality, we have
	\begin{equation*}
		\Psi_{N,N} \geq |\Psi_N|^2, \quad \Psi_{M, M} \geq |\Psi_M|^2, \quad\hbox{ and }\quad 	|\Psi_{M,N}|^2\leq\Psi_{M,M}\Psi_{N,N}.
	\end{equation*}	
	
	The Cauchy-Schwarz inequality also gives
	\begin{align} \label{eq:C-Sargument}
		|\Psi_M|^2 &= \left|\mathbb{E}^w_{\pi \in G} L(\pi) M(\pi)\right|^2 \leq \mathbb{E}^w_{\substack{\pi \in G}} \mathbf{1}_{L(\pi) \neq 0} \cdot \mathbb{E}^w_{\pi \in G} |L(\pi) M(\pi)|^2 \leq \mathbb{E}^w_{\substack{\pi \in G}} \mathbf{1}_{L(\pi) \neq 0} \cdot \Psi_{M,M},
	\end{align}
	and similarly 
	\begin{align*}
		|\Psi_N|^2 &\leq \mathbb{E}^w_{\substack{\pi \in G}} \mathbf{1}_{L(\pi) \neq 0} \cdot \Psi_{N,N}.
	\end{align*}
	In the case $\Psi_{M,M}, \Psi_{N, N} \neq 0$ this implies that the (weighted) non-vanishing proportion is at least
	\begin{equation*}
		\max\left\{\frac{|\Psi_M|^2}{\Psi_{M, M}}, \frac{|\Psi_N|^2}{\Psi_{N, N}}\right\}. 
	\end{equation*} We say that a mollifier $M$ is \emph{efficient} if it leads to a positive proportion of non-vanishing, which is equivalent to $|\Psi_M|^2 \asymp \Psi_{M,M}.$ We generalize the definition of $\beta_Q(M)$ from the previous section to
	\begin{align*}
		\beta(M) := \frac{|\Psi_M|^2}{\Psi_{M, M}},
	\end{align*}
	and we also define $\beta(M)=0$ if $\Psi_{M,M}=0$.
	\begin{remark}\label{remark:zetazeroes} 
		In addition to studying non-vanishing, mollifiers are also used for instance to study the critical zeroes of the Riemann zeta function. In this case the best results are obtained by Levinson's method (for an exposition of Levinson's method, see for example~\cite{Bombieri} or~\cite[Appendix A]{CIS}). The optimization in the context of Levinson's method is slightly different from ours, but one still wants to use a mollifier which leads to a minimal mollified second moment under certain normalization. For instance in Levinson's original work~\cite{Levinson} one roughly wants to minimize
		\[
		\frac{1}{R} \int_{-T}^T \left|\left(\zeta(\sigma+it) + \frac{1}{\log T} \zeta'(\sigma+it)\right) \sum_{m \leq T^\theta} \frac{a_m}{m^{\sigma+it}}\right|^2 dt
		\]
		with $\sigma = 1/2-R/\log T$, and the normalization $a_1 = 1$ (and for instance $a_m \ll m$). More generally, the normalization $\Psi_M = 1$ in our set-up would correspond to $M(s)$ being very close to $1$ for $\re(s)$ large. Variants of our general results can be given in this context.
	\end{remark}
	
	\begin{remark} \label{rem:normalization}
		Let $M, N \colon G \to \mathbb{C}$. Notice that, for any non-zero $u, v \in \mathbb{C}$, we have
		\begin{equation*}
			\Psi_{u M} = u \Psi_M, \quad \Psi_{u M, u M} = |u|^2 \Psi_{M, M}, \quad  \Psi_{u M, v N} = u \overline{v} \Psi_{M, N}, \quad \beta(u M) = \beta(M).
		\end{equation*}
		Due to this, mollifiers $M$ and $u M$ are essentially same. For this reason, whenever we have two mollifiers $M, N$ and make assumptions or statements concerning them, we formulate them in such a way that replacing $M$ by $u M$ and $N$ by $v N$ leads to essentially the same conclusion. This often allows us to assume that $\Psi_M = 1$ and $\Psi_N \in \{0, 1\}$, and the reader may wish to read the proofs in Section~\ref{sec:general mollifying proofs} assuming this.
	\end{remark}
	
	The results in this section concern comparing the efficiency of mollifiers $M, N,$ and $M + \alpha N$. Our first result provides a sufficient condition which implies that $N$ is not a more efficient mollifier than $M$.
	\begin{theorem}\label{thm:criterion1}
		Let $\delta \in [0, 1/10)$ and $M,N \colon G \to \mathbb{C}.$ Assume that $\Psi_{M,M} \neq 0$ and
		\begin{equation}\label{eq:crucial inequality}
			|\Psi_M\Psi_{M,N}| \geq (1-\delta) |\Psi_N|\Psi_{M,M}.
		\end{equation}
		Then
		\begin{equation*}
			\beta(N)\leq (1+4\delta) \beta(M).
		\end{equation*}
	\end{theorem}
	
	Hence, if we have a class of mollifiers $\mathcal{C}$ that contains a mollifier $M$ and show that 
	\begin{equation}
		\label{eq:crucialagain}
		|\Psi_M \Psi_{M, N}| = |\Psi_N| \Psi_{M, M} + o(|\Psi_N| \Psi_{M, M}),
	\end{equation}
	for all $N\in\mathcal{C},$ we can conclude that $M$ is asymptotically the most efficient among all mollifiers in $\mathcal{C}$. In particular Theorem~\ref{thm:main general result DirL} follows by showing that~\eqref{eq:crucialagain} holds in the corresponding set-up.
	
	In many applications one combines two or more mollifiers. Hence given an efficient mollifier $M$ and a different mollifier $N$, it is very natural to ask if combining them leads to a more efficient mollifier. Similarly to Khan, Mili\'cevi\'c and Ngo~\cite{Khan-Mil-Ngo}, we consider the combined mollifier $M+\alpha N$ and optimize $\alpha$. 
	
	Now the non-vanishing proportion obtained from the combined mollifier $M+\alpha N$ is
	\begin{equation*}
		\beta(M+\alpha N) = \frac{|\Psi_M+\alpha\Psi_N|^2}{\Psi_{M,M}+2\re(\overline{\alpha} \Psi_{M,N})+|\alpha|^2\Psi_{N,N}}
	\end{equation*}
	(provided the denominator is non-zero; otherwise by definition $\beta(M+\alpha N)=0$). We first want to find the optimal choice of $\alpha$. 
	
	By symmetry, we consider the case where $M$ is at least as efficient as $N$, so $\beta(M) \geq \beta(N)$. If $\Psi_M = 0$, then by $\beta(M) \geq \beta(N)$ we also have that $\Psi_N = 0$ and thus $\beta(M+\alpha N) = 0$ for every $\alpha \in \mathbb{C}$. Similarly if $\Psi_{N, N} = 0$, then $w(\pi) L(\pi) N(\pi) = 0$ for every $\pi \in G$ and thus $\beta(M + \alpha N) = \beta(M)$ for every $\alpha \in \mathbb{C}$. Hence we can restrict to the case $\Psi_{N, N} \Psi_M \neq 0$.
	
	\begin{theorem}\label{thm:criterion2}
		Assume that $M,N \colon G \to \mathbb{C}$ satisfy $\beta(M)\geq\beta(N)$ and $\Psi_M \Psi_{N, N} \neq 0$. 
		\begin{enumerate}[(i)]
			\item Assume that
			\begin{equation*}
				\Psi_N \Psi_{M,N} = \Psi_M\Psi_{N,N}.
			\end{equation*} Then for all $\alpha\in\C$ with $\alpha\neq-\frac{\Psi_M}{\Psi_N}$, we have $\beta(M+\alpha N)=\beta(M)$.
			\item Assume that
			\begin{equation}\label{eq:corrlowbound}
				\Psi_N \Psi_{M,N} \neq \Psi_M\Psi_{N,N}.
			\end{equation}
			Then, for $\alpha\in\C$, the maximum of the non-vanishing proportion $\beta(M+\alpha N)$ obtained from $M+\alpha N$ occurs for
			\begin{equation} \label{eq:alpha_1def}
				\alpha=\alpha_1:=\frac{\overline{\Psi_M}\Psi_{M,N}-\overline{\Psi_{N}}\Psi_{M,M}}{\overline{\Psi_N}\overline{\Psi_{M,N}}-\overline{\Psi_M}\Psi_{N,N}}
			\end{equation}
			and is
			\begin{align} \label{eq:betaformulas}
				\begin{aligned}
					\beta(M+\alpha_1N) &= \frac{|\Psi_N|^2 \Psi_{M, M} + |\Psi_M|^2 \Psi_{N, N} - 2\re(\overline{\Psi_M}\Psi_N\Psi_{M, N})}{\Psi_{M, M} \Psi_{N, N} - |\Psi_{M, N}|^2} \\
					&= \beta(M) + \frac{|\overline{\Psi_M}\Psi_{M,N}-\overline {\Psi_N}\Psi_{M,M}|^2}{\Psi_{M,M}(\Psi_{M,M}\Psi_{N,N}-|\Psi_{M,N}|^2)}.
				\end{aligned}
			\end{align}
			The denominators in these expressions for $\beta(M+\alpha_1 N)$ are positive.
		\end{enumerate}
	\end{theorem}
	
	From this theorem we will quickly deduce the following corollary giving a simple general criterion on whether there exists $\alpha \in \mathbb{C}$ such that $M+\alpha N$ is a more efficient mollifier than $M$.	
	
	\begin{corollary} \label{cor:MNMMcomp}
		Assume that $M,N \colon G \to \mathbb{C}$ satisfy $\beta(M)\geq\beta(N)$ and $\Psi_M \Psi_{N, N} \neq 0$. 
		\begin{enumerate}[(i)]
			\item Assume that
			\begin{equation*}
				\Psi_M \overline{\Psi_{M,N}} \neq \Psi_N \Psi_{M,M}.
			\end{equation*}
			Then~\eqref{eq:corrlowbound} holds, and hence so does the conclusion of Theorem~\ref{thm:criterion2}(ii). In particular, for $\alpha_1$ as in~\eqref{eq:alpha_1def}, we have 
			\[
			\beta(M+\alpha_1 N) > \beta(M).
			\]
			\item Assume that
			\begin{equation*}
				\Psi_M \overline{\Psi_{M,N}} = \Psi_N \Psi_{M,M}.
			\end{equation*}
			Then $\beta(M + \alpha N) \leq \beta(M)$ for every $\alpha \in \mathbb{C}$.
		\end{enumerate}
	\end{corollary}
	
	In Theorem \ref{thm:criterion1}, we saw a criterion that is useful to show that a mollifier is the most efficient from some class. Corollary~\ref{cor:MNMMcomp}(i) implies the following converse, namely that in an additive class of mollifiers, the most efficient mollifier satisfies (a strong form of) the criterion from Theorem \ref{thm:criterion1}.
	
	\begin{corollary}\label{cor:criterion converse}
		Let $\mathcal{C}$ be a class of mollifiers $M:G\rightarrow \C$ such that if $M,N\in \mathcal C$, then $M+\alpha N\in\mathcal C$ for every $\alpha\in\C$.
		
		Then if $M\in\mathcal C$ is such that $\beta(M)\geq\beta(N)$ for each $N\in\mathcal C$, it must also satisfy $\Psi_M\overline{\Psi_{M,N}}=\Psi_N\Psi_{M,M}$ for each $N\in\mathcal C.$
	\end{corollary}
	
	Our next aim is to use Theorem~\ref{thm:criterion2}(ii) to figure out when $M + \alpha_1 N$ is truly more efficient than the individual mollifiers $M$ and $N$ and when not, i.e. we aim to prove a quantitative version of Corollary~\ref{cor:MNMMcomp}. We formulate our results precisely in terms of a parameter $\delta$, but we do not aim for best possible dependency on it.
	
	The first corollary works best in the case when $M$ and $N$ are both pretty efficient mollifiers.
	
	\begin{corollary}\label{cor:EfficientN}
		Let $\delta \in (0, 1/2)$. Assume that $M,N \colon G \to \mathbb{C}$ satisfy $\Psi_M \Psi_{N, N} \neq 0$ and $\beta(M) \geq \beta(N) > 0$. Then the following hold.
		\begin{enumerate}[(i)]
			\item Assume that
			\begin{equation}\label{eq:corrlowboundCoreffi}
				|\overline{\Psi_M} \Psi_{M,N}-\overline{\Psi_N}\Psi_{M,M}|\geq \delta^2 \cdot |\Psi_N|\Psi_{M, M} .
			\end{equation}
			Then the conclusion of Theorem~\ref{thm:criterion2}(ii) holds. Furthermore
			\[
			\beta(M+\alpha_1N) - \beta(M) \geq \delta^4 \beta(N).
			\]
			Hence, if $\beta(N) \gg 1$ and~\eqref{eq:corrlowboundCoreffi} holds with $\delta \gg 1$, then the combined mollifier $M+\alpha_1N$ is truly more efficient than the individual mollifiers $M$ and $N$.
			\item Assume that
			\begin{equation}\label{eq:EfficientCorAssumpiiup}
				|\overline{\Psi_M} \Psi_{M,N}-\overline{\Psi_N} \Psi_{M,M}| < \delta^2 \cdot |\Psi_N| \Psi_{M, M}
			\end{equation}
			and
			\begin{equation}\label{eq:EfficientCorAssumpiilow}
				|\Psi_N \Psi_{M,N} - \Psi_M \Psi_{N,N}|\geq \delta \cdot \Psi_{N, N} |\Psi_M|.
			\end{equation}
			Then the conclusion of Theorem~\ref{thm:criterion2}(ii) holds. Furthermore
			\[
			\beta(M+\alpha N) - \beta(M) \leq \delta^2 \frac{\beta(N)}{\beta(M)} \leq \delta^2
			\quad \text{for every $\alpha\in\C.$}
			\]
			Hence, if~\eqref{eq:EfficientCorAssumpiiup} and~\eqref{eq:EfficientCorAssumpiilow} hold with $\delta = o(1)$, then the combined mollifier is not significantly more efficient than the individual mollifier $M$.
		\end{enumerate}
	\end{corollary}
	Notice that we are in case (i) in particular if $M$ and $N$ are so uncorrelated that $|\Psi_M \Psi_{M, N}| = o(|\Psi_N| \Psi_{M, M})$ (although such total uncorrelation is probably rare in practice).

	\begin{remark}
		If, $\beta(M),\beta(N)> 0$, we can always normalize so that $\Psi_M = \Psi_N = 1$. If $M$ and $N$ are both normalized and efficient, Corollary~\ref{cor:EfficientN}(i) shows that combining them is useful as soon as $\Psi_{M, N}$ is not very close to $\min\{\Psi_{M, M}, \Psi_{N, N}\}$. There is an alternative and more intuitive way to see this: by symmetry, we can assume that
		\[
		\frac{1}{\Psi_{M, M}} = \beta(M) \geq \beta(N) = \frac{1}{\Psi_{N, N}}.
		\]
		Now the Cauchy-Schwarz argument~\eqref{eq:C-Sargument} yielding the non-vanishing proportion $1/\Psi_{M, M}$ has equality if and only if, for every $\pi \in G$ with $w(\pi) \neq 0$, we have
		\[
		L(\pi) M(\pi) =
		\begin{cases}
			\frac{1}{\beta(M)}, & \text{if $L(\pi) M(\pi) \neq 0,$} \\
			0, & \text{if $L(\pi) M(\pi) = 0$.}
		\end{cases}
		\]
		But if this holds, then $\Psi_{M, N} = \frac{1}{\beta(M)} \Psi_N = \frac{1}{\beta(M)} = \Psi_{M, M}$. Hence if $\Psi_{M, N}$ is not very close to $\Psi_{M, M}$, the Cauchy-Schwarz argument~\eqref{eq:C-Sargument} is not sharp, and one can improve on the non-vanishing proportion (for instance by replacing the Cauchy-Schwarz inequality by an exact formula as in Lemma~\ref{le:exactCS} below). Our results give a more quantitative version of this heuristic.
	\end{remark}
	Corollary~\ref{cor:EfficientN} is useful unless $\beta(N)$ is very small, i.e. $N$ is inefficient. In particular it works well if $\beta(N) > \frac{\delta}{4}\beta(M)$, say. Next we investigate the case when $N$ is significantly less efficient than $M$, more precisely the case $\beta(N) \leq \frac{\delta}{4}\beta(M)$. Such an inefficient second mollifier piece has been useful in the context of Levinson's method, see~\cite{BCY} (the normalization in that case is $\Psi_{N, N} \asymp 1$ and $\Psi_N = o(1)$, as \cite[(2.4)]{BCY} corresponds to having $\Psi_N = o(1)$ in our setting).
	\begin{corollary}\label{cor:NonEfficientN}
		Let $\delta \in (0, 1/2)$. Assume that $M,N \colon G \to \mathbb{C}$ satisfy $\Psi_{M} \Psi_{N, N} \neq 0$ and $\beta(N) \leq \frac{\delta}{4} \beta(M)$. Then the conclusion of Theorem~\ref{thm:criterion2}(ii) holds. Furthermore:
		\begin{enumerate}[(i)]
			\item If
			\begin{equation*}
				|\Psi_{M,N}|^2 \geq \delta \cdot \Psi_{M, M} \Psi_{N, N},
			\end{equation*}
			then
			\[
			\beta(M+\alpha_1N) - \beta(M) \geq \frac{\delta}{4} \beta(M).
			\]
			Hence, for $\delta \gg 1$, the combined mollifier $M+\alpha_1N$ is truly more efficient than the individual mollifiers $M$ and $N$.
			\item If
			\begin{equation*}
				|\Psi_{M,N}|^2 \leq \delta \cdot \Psi_{M, M} \Psi_{N, N},
			\end{equation*}
			then
			\[
			\beta(M+\alpha N) - \beta(M) \leq 5\delta \cdot \beta(M) \quad \text{for every $\alpha \in \mathbb{C}.$}
			\]
			Hence, for $\delta = o(1)$, the combined mollifier $M+\alpha_1N$ is not significantly more efficient than the individual mollifiers $M$ and $N$.
		\end{enumerate}
	\end{corollary}
	
	\begin{remark}
		If $\beta(M)\geq\beta(N)$, Corollaries~\ref{cor:EfficientN} and~\ref{cor:NonEfficientN} essentially cover all possible cases except when $\beta(N) \asymp 1$ and both 
		\begin{equation*}
			\lab\overline{\Psi_M}\Psi_{M,N}-\overline{\Psi_N}\Psi_{M,M}\rab=o(|\Psi_N| \Psi_{M, M}) \ \ \  \text{and} \ \ \  \lab\Psi_N\Psi_{M,N}-\Psi_M\Psi_{N,N}\rab = o(|\Psi_M| \Psi_{N, N}).
		\end{equation*}
		Furthermore, if $\Psi_N \Psi_{M, N} - \Psi_M \Psi_{N, N} = 0$, we can apply Theorem~\ref{thm:criterion2}(i), and if $\overline{\Psi_M} \Psi_{M, N} - \overline{\Psi_N} \Psi_{M, M} = 0$, we can apply Corollary~\ref{cor:MNMMcomp}(ii). Hence the only remaining case is 
		$\beta(M) \geq \beta(N) \gg 1$ and both
		\begin{align*}
			0 &\neq \overline{\Psi_M}\Psi_{M,N}-\overline{\Psi_N}\Psi_{M,M}=o(|\Psi_N| \Psi_{M, M}) \\ \text{and} \quad  0 &\neq \Psi_N\Psi_{M,N}-\Psi_M\Psi_{N,N} = o(|\Psi_M| \Psi_{N, N}).
		\end{align*}
		
		It is not possible to say anything in general in this situation without further assumptions, as it can lead to some degenerate setting.
		Indeed, if we also normalize to $\Psi_M=\Psi_N=1$, the above gives $\Psi_{M,N}=\Psi_{M,M}+o(1)=\Psi_{N,N}+o(1)$. In other words, $M$ and $N$ are essentially the same mollifiers in the sense that they are not distinguishable by considering only the first and second moments, and they are also very correlated: we have almost equality in the Cauchy-Schwarz inequality $\Psi_{M,N}^2 \leq \Psi_{M,M}\Psi_{N,N}$. 
		
		In this case when $M$ and $N$ are essentially equal, note that $\alpha_1=\frac{o(1)}{o(1)}$, and much cannot be said without further assumptions, since the lower order terms can come into play. To illustrate this point, consider two mollifiers $M_1$, $M_2$, with $\Psi_{M_j}, \Psi_{M_j,M_j}\asymp 1$, and set
		\begin{equation*}
			\begin{aligned}
				M&=M_1+\frac{1}{Q^{100}}M_2,\\
				N&=M_1.
			\end{aligned}
		\end{equation*} In this case, the individual mollifiers $M,N$ are essentially equal to $M_1$, and so is the combined mollifier $M+\alpha N$ for $\alpha\neq -1$, but $M-N$ is essentially $M_2$. Therefore $\beta(M+\alpha N)$ is essentially constant, except around $\alpha\approx -1$, where it takes either a maximum or a minimum, depending on the relative efficiency of $M_1$ and $M_2.$
	\end{remark}

	\section{Proofs of the general results} \label{sec:general mollifying proofs}
	In this section we prove Theorems~\ref{thm:criterion1} and~\ref{thm:criterion2} as well as Corollaries~\ref{cor:MNMMcomp},~\ref{cor:criterion converse},~\ref{cor:EfficientN}, and~\ref{cor:NonEfficientN}.	Let us first prove a very simple lemma.
	
	\begin{lemma}\label{le:SimpleLemma}
		Assume that $M,N \colon G \to \mathbb{C}$ satisfy $\Psi_{M} \Psi_{N, N} \neq 0$, and $\beta(N) \leq \beta(M)$. If 
		\begin{equation}
			\label{eq:SimpleLemmaAssu}
			\Psi_N \Psi_{M, N} = \Psi_M \Psi_{N, N},
		\end{equation}
		then 
		\[
		\Psi_M \overline{\Psi_{M, N}} = \Psi_N \Psi_{M, M}.
		\]
	\end{lemma}
	\begin{proof}
		By the assumptions we must have $\Psi_N \neq 0$, and $\beta(M) \geq \beta(N)$ implies $|\Psi_N|^2\Psi_{M, M} \leq |\Psi_M|^2\Psi_{N, N}$. Using also the Cauchy-Schwarz inequality and~\eqref{eq:SimpleLemmaAssu}, we obtain
		\[
		|\Psi_N|^2|\Psi_{M, N}|^2 \leq |\Psi_N|^2\Psi_{M, M} \Psi_{N, N} \leq |\Psi_M|^2\Psi_{N, N}^2 =|\Psi_N|^2 |\Psi_{M, N}|^2.
		\]
		Therefore we have equalities in the inequalities above, and the second equality with~\eqref{eq:SimpleLemmaAssu} (and the fact that $\Psi_{N, N} \in \mathbb{R}$) gives
		$$
		\Psi_N\Psi_{M,M}=\frac{\Psi_M\overline{\Psi_M}\Psi_{N,N}}{\overline{\Psi_N}}=\Psi_M\overline{\Psi_{M,N}}.
		$$
	\end{proof}
	
	The proofs of Theorem~\ref{thm:criterion1} and~\ref{thm:criterion2}(i) are very quick, and it is also quick to deduce Corollaries~\ref{cor:MNMMcomp} and~\ref{cor:criterion converse} from  Theorem~\ref{thm:criterion2}.
	
	\begin{proof}[Proof of Theorem~\ref{thm:criterion1}]
		We can assume that $\Psi_{N,N} \neq 0$ since otherwise $\beta(N) = 0$ and the claim is trivial. Now by the assumption~\eqref{eq:crucial inequality} and the Cauchy-Schwarz inequality
		\[
		(1-\delta)^2 |\Psi_N|^2 \Psi_{M, M}^2 \leq |\Psi_M|^2|\Psi_{M, N}|^2 \leq |\Psi_M|^2\Psi_{M, M} \Psi_{N, N}.
		\]
		Thus
		\[
		\beta(N) = \frac{|\Psi_N|^2}{\Psi_{N, N}} \leq \frac{|\Psi_M|^2}{(1-\delta)^2 \Psi_{M, M}} = \frac{\beta(M)}{(1-\delta)^2} \leq (1+4\delta) \beta(M).
		\]
	\end{proof}
	
	\begin{proof}[Proof of Theorem \ref{thm:criterion2}(i)]
		Since $\Psi_M \Psi_{N, N} \neq 0$, the condition
		\begin{equation*}
			\Psi_N\Psi_{M,N}=\Psi_M\Psi_{N,N}
		\end{equation*} can hold only if $\Psi_N \Psi_{M,N}\neq 0$.  Lemma~\ref{le:SimpleLemma} implies that $\Psi_{M, N} = \frac{\overline{\Psi_N}}{\overline{\Psi_M}}\Psi_{M,M}$ and $\Psi_{N,N}=\frac{\Psi_N}{\Psi_{M}}\Psi_{M,N}=\frac{|\Psi_N|^2}{|\Psi_M|^2}\Psi_{M,M}$.
		
		Hence, provided the denominator below is non-zero (which is obvious from the second line since $\alpha \neq -\frac{\Psi_M}{\Psi_N}$),
		$$
		\begin{aligned}
			\beta(M+\alpha N)&=\frac{|\Psi_M+\alpha\Psi_N|^2}{\Psi_{M,M}+2\re(\overline{\alpha} \Psi_{M,N})+|\alpha|^2\Psi_{N,N}}\\&= \frac{|\Psi_M|^2|1+\alpha\Psi_N/\Psi_M|^2}{\Psi_{M, M}(1+2\re(\alpha\Psi_N/\Psi_M)+|\alpha|^2|\Psi_N|^2/|\Psi_M|^2)} = \frac{|\Psi_M|^2}{\Psi_{M, M}} = \beta(M).
		\end{aligned}
		$$ 
	\end{proof}
	
	\begin{proof}[Proof of Corollary~\ref{cor:MNMMcomp} assuming Theorem~\ref{thm:criterion2}]
		The claim (i) follows immediately from Lemma~\ref{le:SimpleLemma} (in the contrapositive) and Theorem~\ref{thm:criterion2}(ii). 
		
		If $\Psi_N \Psi_{M,N} = \Psi_M\Psi_{N,N}$, the claim (ii) follows from Theorem~\ref{thm:criterion2}(i). Otherwise it follows from Theorem~\ref{thm:criterion2}(ii).
	\end{proof}
	
	\begin{proof}[Proof of Corollary~\ref{cor:criterion converse} assuming Corollary~\ref{cor:MNMMcomp}]
		Assume that $\mathcal C, M$ are as in the statement of Corollary~\ref{cor:criterion converse}. Now suppose there exists some $N \in \mathcal{C}$ such that
		\begin{equation}
			\label{eq:corcounter}
			\Psi_M\overline{\Psi_{M,N}}\neq \Psi_N\Psi_{M,M}.
		\end{equation}	
		
		If $\Psi_{N, N} = 0$, then $\Psi_{M, N} = \Psi_N = 0$ and~\eqref{eq:corcounter} cannot hold. Similarly if we had $\Psi_M = 0$, then $\beta(M) \geq \beta(N)$ would also imply $\Psi_N = 0$ and thus~\eqref{eq:corcounter} cannot hold. Hence we can assume that $\Psi_{N, N} \Psi_M \neq 0$, so Corollary~\ref{cor:MNMMcomp}(i) is applicable. It implies that $\beta(M+\alpha_1 N) > \beta(M)$ for $\alpha_1$ as in~\eqref{eq:alpha_1def}. But $M + \alpha_1 N \in \mathcal{C}$ by our assumptions, so this is a contradiction with the assumption that $\beta(M) \geq \beta(N)$ for every $N \in \mathcal{C}$. 
	\end{proof}
	
	Before turning to the proofs of Theorem~\ref{thm:criterion2}(ii) and Corollaries~\ref{cor:EfficientN} and~\ref{cor:NonEfficientN}, we provide two auxiliary lemmas. 
	
	The first lemma is an ``exact'' version of the Cauchy-Schwarz inequality which we will use in the proof of the second lemma. It follows from a well-known proof of the inequality but we include its proof for completeness.
	\begin{lemma}\label{le:exactCS}
		Let $V$ be an inner-product space. Then for any $u,v\in V$, we have
		$$
		\norm u^2\norm v^2-|\langle u,v\rangle|^2 = \frac{1}{\norm{v}^2}\norm{\norm v^2 u-\langle u,v\rangle v}^2.
		$$
	\end{lemma}
	\begin{proof}
		Let $w := u - \frac{\langle u, v \rangle}{\norm v^2} v$, so that $w$ and $v$ are orthogonal. Then  $a_1 := w$ and $a_2 := \frac{\langle u, v \rangle}{\norm v^2} v$ are also orthogonal and by the Pythagorean theorem we have 
		\begin{align*}
			\norm{a_1 + a_2}^2 &= \norm {a_1}^2 + \norm {a_2}^2,
		\end{align*}
		so that
		\begin{align*}
			\norm{u}^2  =\norm{u - \frac{\langle u, v \rangle}{\norm v^2} v}^2 + \norm{\frac{\langle u, v \rangle}{\norm v^2} v}^2 = \frac{1}{\norm{v}^4}\norm{u \norm{v}^2 - \langle u, v \rangle v}^2 + \frac{1}{\norm{v}^2} |\langle u, v \rangle|^2
		\end{align*}
		which implies the claim.
	\end{proof}
	
	The second lemma allows us to ensure that some denominators are non-zero (or satisfy a lower bound).
	\begin{lemma}\label{le:nonzerocorrs}
		Assume that $M,N \colon G \to \mathbb{C}$ satisfy $\Psi_M \Psi_{N, N} \neq 0$. Assume that
		\begin{equation}\label{eq:corrlowbound2qual}
			\Psi_N\Psi_{M,N} \neq \Psi_M\Psi_{N,N}.
		\end{equation}
		Then the following hold.
		\begin{enumerate}[(i)]
			\item We have
			\begin{equation*}
				\Psi_{M, M} \Psi_{N, N} - |\Psi_{M, N}|^2 > 0.
			\end{equation*}
			\item We have, for every $\alpha \in \mathbb{C}$,
			\begin{equation*}
				\Psi_{M + \alpha N, M + \alpha N} \neq 0.
			\end{equation*}
			\item If for some $\delta > 0$, \eqref{eq:corrlowbound2qual} holds in the quantitative form
			\begin{equation}\label{eq:corrlowbound2quant}
				|\Psi_N\Psi_{M,N}-\Psi_M\Psi_{N,N}|\geq\delta |\Psi_M| \Psi_{N, N},
			\end{equation}
			then
			\begin{equation*}
				\Psi_{M, M} \Psi_{N, N} - |\Psi_{M, N}|^2 \geq \delta^2 |\Psi_M|^2 \Psi_{N, N}.
			\end{equation*}
		\end{enumerate}
	\end{lemma}
	\begin{proof}
		Consider the inner-product space, where, for two functions $A, B \colon G \to \mathbb{C}$,
		\begin{align*}
			\langle A, B \rangle = \mathbb{E}^w_{\pi \in G} A(\pi) \overline{B(\pi)}.
		\end{align*}
		Let us consider part (iii) first. By Lemma \ref{le:exactCS},
		\begin{align*}
			\Psi_{M, M} \Psi_{N, N} - |\Psi_{M, N}|^2 = \frac{1}{\Psi_{N, N}} \Vert \Psi_{N, N} L(\cdot)M(\cdot) - \Psi_{M, N} L(\cdot)N(\cdot) \Vert^2.
		\end{align*}
		By~\eqref{eq:corrlowbound2quant} and the Cauchy-Schwarz inequality,
		\begin{align*}
			\delta^2 |\Psi_M|^2 \Psi_{N, N}^2 &\leq |\Psi_N \Psi_{M, N} - \Psi_M \Psi_{N, N}|^2 = |\langle \Psi_{M,N} L(\cdot)N(\cdot) - \Psi_{N, N}L(\cdot)M(\cdot), 1\rangle|^2 \\
			&\leq \Vert \Psi_{N, N} L(\cdot)M(\cdot) - \Psi_{M, N} L(\cdot)N(\cdot) \Vert^2,
		\end{align*}
		and part (iii) follows. Part (i) follows similarly using~\eqref{eq:corrlowbound2qual} in place of~\eqref{eq:corrlowbound2quant}.
		
		Let us turn to part (ii) and assume for contradiction that $\Psi_{M + \alpha N,M+\alpha N} = 0$ for some $\alpha \in \mathbb{C}$. This implies that $w(\pi)L(\pi)(M(\pi) + \alpha N(\pi)) = 0$ for every $\pi \in G$. If this is the case, then $\Psi_M = - \alpha \Psi_N$ and $\Psi_{M, N} = - \alpha \Psi_{N, N}$. Since $\Psi_M \neq 0$, one cannot have $\alpha = 0$. Hence we obtain 
		\[
		\Psi_N \Psi_{M, N} = \frac{-1}{\alpha} \Psi_M \cdot (- \alpha \Psi_{N, N}) = \Psi_M \Psi_{N, N}
		\] 
		which contradicts the assumption~\eqref{eq:corrlowbound2qual}.
	\end{proof}
	
	Now we are ready to turn to the proofs of Theorem~\ref{thm:criterion2}(ii) and Corollaries~\ref{cor:EfficientN} and~\ref{cor:NonEfficientN}.	
	\begin{proof}[Proof of Theorem~\ref{thm:criterion2}(ii)]
		Recall
		\begin{equation*}
			\beta(M+\alpha N) = \frac{\Psi_{M + \alpha N}}{\Psi_{M+\alpha N, M+\alpha N}} = \frac{|\Psi_M+\alpha\Psi_N|^2}{\Psi_{M,M}+2\re(\overline{\alpha} \Psi_{M,N})+|\alpha|^2\Psi_{N,N}}.
		\end{equation*}
		The denominator is always non-zero by Lemma~\ref{le:nonzerocorrs}(ii).
		
		By Theorem~\ref{thm:criterion1} with $\delta = 0$ we see that $\beta(M+\alpha_0 N)$ is maximal for a given $\alpha_0 \in \mathbb{C}$ if we can show that
		\begin{equation}\label{eq:GoalInPfofThm2.5}
			|\Psi_{M+\alpha_0 N} \Psi_{M + \alpha_0 N, M+\alpha_0 N + \Delta N}| = |\Psi_{M+\alpha_0 N + \Delta N}| \Psi_{M+\alpha_0 N, M+\alpha_0 N}
		\end{equation}
		for all $\Delta \in \mathbb{C}$. Our assumptions imply that $\Psi_M\neq 0.$ We will first consider the normalized case $\Psi_M=1$ and $\Psi_N \in \{0, 1\}$ and then deduce the general case from it.
		
		\textbf{Case 1: $\Psi_M = 1$ and $\Psi_N \in \{0, 1\}$.}
		
		\textbf{Case 1.1: $\Psi_M = 1$ and $\Psi_N = 0$.} 
		In this case~\eqref{eq:GoalInPfofThm2.5} reduces to
		\[
		|\Psi_{M+ \alpha_0 N, M+\alpha_0N+\Delta N}| = \Psi_{M+\alpha_0 N, M+\alpha_0 N} 
		\]
		for every $\Delta \in \mathbb{C}$. This holds if and only if $\Psi_{M+\alpha_0 N, N} = 0$ which is in turn equivalent to
		\begin{equation}\label{eq:optimal alpha N=0}
			\Psi_{M, N} + \alpha_0 \Psi_{N, N} = 0 \iff \alpha_0 = -\frac{\Psi_{M, N}}{\Psi_{N, N}} = \frac{\overline{\Psi_M} \Psi_{M, N} - \overline{\Psi_N} \Psi_{M, M}}{\overline{\Psi_N} \overline{\Psi_{M, N}} - \overline{\Psi_M} \Psi_{N, N}} = \alpha_1.
		\end{equation}
		A direct computation gives
		$$
		\beta(M+\alpha_1 N)=\frac{\Psi_{N,N}}{\Psi_{M,M}\Psi_{N,N}-|\Psi_{M,N}|^2}=\frac{1}{\Psi_{M,M}}+\frac{|\Psi_{M,N}|^2}{\Psi_{M,M}\lz\Psi_{M,M}\Psi_{N,N}-|\Psi_{M,N}|^2\pz},
		$$ which gives \eqref{eq:betaformulas} in this case. Furthermore the denominators are positive by Lemma~\ref{le:nonzerocorrs}(i).
		
		\textbf{Case 1.2: $\Psi_M = \Psi_N = 1$.} 
		Now~\eqref{eq:GoalInPfofThm2.5} reduces to
		\begin{align} \label{eq:neededeq}
			\begin{aligned}
				&  |(1+\alpha_0)| \cdot |\Psi_{M+\alpha_0 N, M+\alpha_0 N} + \overline{\Delta}\Psi_{M+\alpha_0 N, N}| \\
				&= |(1+\alpha_0) \Psi_{M+\alpha_0 N, M+\alpha_0N} + \Delta \Psi_{M+\alpha_0N, M+\alpha_0 N}|
			\end{aligned}
		\end{align}
		If this held for every $\Delta \in \mathbb{C}$ for $\alpha_0 = -1$, we would have $\Psi_{M-N, M-N} = 0$ which is not possible by Lemma \ref{le:nonzerocorrs}(ii).
		
		For $\alpha_0 \neq -1$, since $\Psi_{M+\alpha_0 N,M+\alpha_0 N}\in\R_{>0}$,~\eqref{eq:neededeq} is equivalent to
		\begin{align*}
			|\Psi_{M+\alpha_0 N, M+\alpha_0 N} + \Delta \overline{\Psi_{M+\alpha_0 N, N}}| &= \left|\Psi_{M+\alpha_0 N, M+\alpha_0N} + \frac{\Delta}{1+\alpha_0} \Psi_{M+\alpha_0N, M+\alpha_0 N}\right|.
		\end{align*}
		Now this holds for every $\Delta \in \mathbb{C}$ if and only if we have
		\[
		(1+\alpha_0)\overline{\Psi_{M+\alpha_0 N, N}} = \Psi_{M+\alpha_0 N, M+ \alpha_0 N}.
		\]
		This opens to
		\begin{align*}
			&\overline{\Psi_{M, N}} + \alpha_0 \overline{\Psi_{M, N}} + \overline{\alpha_0} \Psi_{N, N} + |\alpha_0|^2 \Psi_{N, N} = \Psi_{M, M} + \alpha_0 \overline{\Psi_{M, N}} + \overline{\alpha_0} \Psi_{M, N} + |\alpha_0|^2 \Psi_{N, N}
		\end{align*}
		which is equivalent to
		\[
		\overline{\alpha_0}(\Psi_{M, N} - \Psi_{N, N}) = \overline{\Psi_{M, N}} - \Psi_{M, M}.
		\]
		Thus
		\begin{equation}\label{eq:optimal alpha normalized}
			\alpha_0 = \frac{\Psi_{M, N} - \Psi_{M, M}}{\overline{\Psi_{M, N}}-\Psi_{N, N}} = \alpha_1
		\end{equation}
		yields the maximum. Since $\Psi_{N, N} \in \mathbb{R}$, the denominator is non-zero by~\eqref{eq:corrlowbound}. 
		
		Substituting the value of $\alpha_1$ into
		\begin{equation*}
			\beta(M+\alpha N) = \frac{|1+\alpha|^2}{\Psi_{M,M}+2\re(\overline{\alpha} \Psi_{M,N})+|\alpha|^2\Psi_{N,N}}.
		\end{equation*}
		and doing some arithmetic yields
		\begin{equation}\label{eq:beta formula normalized}
			\begin{aligned}
				\beta(M + \alpha_1 N) = \frac{\Psi_{M, M} + \Psi_{N, N} - 2 \re(\Psi_{M, N})}{\Psi_{M, M} \Psi_{N, N} - |\Psi_{M, N}|^2}
			\end{aligned}
		\end{equation}
		and	
		\begin{equation} \label{eq:beta formula normalized 2}
			\beta(M + \alpha_1 N)-\beta(M)=\frac{|\Psi_{M,N}-\Psi_{M,M}|^2}{\Psi_{M,M}(\Psi_{M,M}\Psi_{N,N}-|\Psi_{M,N}|^2)},
		\end{equation} which proves \eqref{eq:betaformulas} in this case. Furthermore the denominators are positive by Lemma~\ref{le:nonzerocorrs}(i).
		
		\textbf{Case 2: $\Psi_M \in \mathbb{C} \setminus \{0\}$ and $\Psi_N \in \mathbb{C}$.}
		We let $\widetilde M:=M/\Psi_M$ and 
		$$
		S(M,N):=\{\alpha\in\C: \beta(M+\alpha N) \text{ is maximal}\}.
		$$
		Notice that for $u,v\in\C\setminus\{0\}$, we have $S(uM,N)=u\cdot S(M,N)$, and $S(M,vN)=1/v\cdot S(M,N)$.
		
		\textbf{Case 2.1: $\Psi_M \in \mathbb{C} \setminus \{0\}$ and $\Psi_N = 0$.}
		Now~\eqref{eq:optimal alpha N=0} shows that 
		$$
		-\frac{\Psi_{M,N}}{\Psi_M\Psi_{N,N}}=-\frac{\Psi_{\widetilde M,N}}{\Psi_{N,N}}\in S(\widetilde M,N)=\frac{1}{\Psi_M}S\lz M,N\pz,
		$$
		so $$
		-\frac{\Psi_{M,N}}{\Psi_{N,N}}=\frac{\overline{\Psi_M}\Psi_{M,N}-\overline{\Psi_N}\Psi_{M,M}}{\overline{\Psi_N\Psi_{M,N}}-\overline{\Psi_M}\Psi_{N,N}}\in S(M,N), 
		$$ and also \eqref{eq:betaformulas} follows directly.
		
		\textbf{Case 2.2: $\Psi_M, \Psi_N \in \mathbb{C} \setminus \{0\}.$} Writing $\widetilde N=N/\Psi_N$, \eqref{eq:optimal alpha normalized} implies that
		$$
		\begin{aligned}
			\frac{\Psi_{\widetilde M,\widetilde N}-\Psi_{\widetilde M,\widetilde M}}{\overline{\Psi_{\widetilde M,\widetilde N}}-\Psi_{\widetilde N,\widetilde N}}\in S(\widetilde M,\widetilde N) & \iff		\frac{\frac{\Psi_{M,N}}{\Psi_M\overline{\Psi_N}}-\frac{\Psi_{M,M}}{|\Psi_M|^2}}{\frac{\overline{\Psi_{M,N}}}{\overline{\Psi_M}\Psi_N}-\frac{\Psi_{N,N}}{|\Psi_N|^2}}\in \frac{\Psi_N}{\Psi_M} S(M,N)\\
			& \iff	\frac{\frac{\Psi_{M,N}}{\overline{\Psi_N}}-\frac{\Psi_{M,M}}{\overline{\Psi_M}}}{\frac{\overline{\Psi_{M,N}}}{\overline{\Psi_M}}-\frac{\Psi_{N,N}}{\overline{\Psi_N}}}=\alpha_1 \in S(M,N).
		\end{aligned}
		$$
		To prove \eqref{eq:betaformulas} in the general case, we either make a direct computation, or notice that for $\alpha_1\in S(M,N)$ and $\widetilde\alpha_1 \in S(\widetilde M,\widetilde N)$, we must have $\beta(M+\alpha_1 N)=\beta(\widetilde M+\widetilde\alpha_1\widetilde N).$ Therefore using \eqref{eq:beta formula normalized} and~\eqref{eq:beta formula normalized 2}, we obtain
		$$
		\begin{aligned}
			\beta(M+\alpha_1 N)&=\frac{\Psi_{\widetilde M,\widetilde M}+\Psi_{\widetilde N,\widetilde N}-2\re (\Psi_{\widetilde M,\widetilde N})}{\Psi_{\widetilde M,\widetilde M}\Psi_{\widetilde N,\widetilde N}-|\Psi_{\widetilde M,\widetilde N}|^2}\\
			&=\frac{|\Psi_N|^2\Psi_{M,M}+|\Psi_M|^2\Psi_{N,N}-2\re\lz \overline{\Psi_M}\Psi_N\Psi_{M,N}\pz}{\Psi_{M,M}\Psi_{N,N}-|\Psi_{M,N}|^2}\\
			&=\frac{|\Psi_M|^2}{\Psi_{M,M}}+\frac{\lab\overline{\Psi_{M}}\Psi_{M,N}-\overline{\Psi_N}\Psi_{M,M}\rab^2}{\Psi_{M,M}\lz\Psi_{M,M}\Psi_{N,N}-|\Psi_{M,N}|^2\pz}.
		\end{aligned}
		$$
		
	\end{proof}
	
	\begin{proof}[Proof of Corollary~\ref{cor:EfficientN}]
		Since $\beta(N) > 0$, we have $\Psi_N \neq 0$. Let us first show that Theorem~\ref{thm:criterion2}(ii) is applicable. In part (ii) this is immediate from~\eqref{eq:EfficientCorAssumpiilow} whereas in part (i) this follows from~\eqref{eq:corrlowboundCoreffi} and Corollary~\ref{cor:MNMMcomp}(i).
		
		Hence we can apply Theorem~\ref{thm:criterion2}(ii), so the maximum value of $\beta(M+\alpha N)$ occurs for $\alpha = \alpha_1$ and furthermore
		\begin{align*}
			\beta(M+\alpha_1N) - \beta(M) = \frac{|\overline{\Psi_M}\Psi_{M,N}-  \overline{\Psi_N}\Psi_{M,M}|^2}{\Psi_{M,M}(\Psi_{M,M}\Psi_{N,N}-|\Psi_{M,N}|^2)}.
		\end{align*}
		Next we establish the remaining two claims by bounding the right hand side of the above expression. In case (i), the denominator is $\leq \Psi_{M, M}^2 \Psi_{N, N}$, whereas the numerator is $\geq \delta^4 |\Psi_N|^2\Psi_{M, M}^2$ by~\eqref{eq:corrlowboundCoreffi},  so (i) follows.
		
		In case (ii), combining~\eqref{eq:EfficientCorAssumpiilow} and Lemma~\ref{le:nonzerocorrs}(iii) implies that the denominator is $\geq \Psi_{M, M} \cdot \delta^2 |\Psi_M|^2\Psi_{N, N}$ whereas~\eqref{eq:EfficientCorAssumpiiup} implies that the numerator is $\ll \delta^4 |\Psi_N|^2\Psi_{M, M}^2$, so also (ii) follows.
	\end{proof}
	
	\begin{proof}[Proof of Corollary~\ref{cor:NonEfficientN}]
		By $\beta(N) \leq \frac{\delta}{4} \beta(M) < \beta(M)$ we know that $|\Psi_N|^2\Psi_{M, M} < |\Psi_M|^2\Psi_{N, N}$. Using the Cauchy-Schwarz inequality and this, we see that $|\Psi_N\Psi_{M,N}|^2\leq |\Psi_N|^2\Psi_{M,M}\Psi_{N,N}<|\Psi_M|^2 \Psi_{N,N}^2$, so $\Psi_N\Psi_{M,N}\neq \Psi_M\Psi_{N,N}$. Therefore Theorem~\ref{thm:criterion2}(ii) is applicable and the maximum value of $\beta(M+\alpha N)$ occurs for $\alpha = \alpha_1$ and furthermore
		\begin{align}\label{eq:betadifnoneffcorproof}
			\beta(M+\alpha_1N) - \beta(M) = \frac{|\overline{\Psi_M}\Psi_{M,N}-\overline{\Psi_N}\Psi_{M,M}|^2}{\Psi_{M,M}(\Psi_{M,M}\Psi_{N,N}-|\Psi_{M,N}|^2)}.
		\end{align}
		If now $|\Psi_{M, N}|^2 \leq \delta \Psi_{M, M} \Psi_{N, N}$, the right-hand side above is
		\[
		\leq \frac{2(|\Psi_M|^2|\Psi_{M,N}|^2 + |\Psi_N|^2\Psi_{M, M}^2)}{\frac{1}{2}\Psi_{M, M}^2 \Psi_{N, N}} \leq \frac{4\delta|\Psi_M|^2}{\Psi_{M, M}} + \frac{4|\Psi_N|^2}{\Psi_{N, N}} = 4\delta\beta(M) + 4\beta(N) \leq 5\delta \beta(M). 
		\]
		which proves (ii).
		
		Since $\Psi_{M, M} > 0$, we have $|\overline{\Psi_M}\Psi_{M, N} - \overline{\Psi_N}\Psi_{M, M}| \geq |\Psi_M\Psi_{M, N}| - |\Psi_N|\Psi_{M, M}$. If now $|\Psi_{M, N}|^2 \geq \delta \Psi_{M, M} \Psi_{N, N}$, then        
		\begin{equation}
			\label{eq:ineffcorhelp}
			|\overline{\Psi_M}\Psi_{M, N} - \overline{\Psi_N}\Psi_{M, M}| \geq \frac{|\Psi_M\Psi_{M, N}|}{2} + \frac{|\Psi_M|}{2} \sqrt{\delta \Psi_{M, M} \Psi_{N, N}} - |\Psi_N|\Psi_{M, M}.
		\end{equation}
		Since $\beta(N) \leq \frac{\delta}{4}\beta(M)$, we have $|\Psi_M|^2\Psi_{N, N} \geq \frac{4}{\delta} |\Psi_N|^2\Psi_{M, M}$. Thus here
		\[
		\frac{|\Psi_M|}{2} \sqrt{\delta \Psi_{M, M} \Psi_{N, N}} - |\Psi_N|\Psi_{M, M} \geq 0.
		\]
		Inserting this into~\eqref{eq:ineffcorhelp}, we obtain that the right hand side of~\eqref{eq:betadifnoneffcorproof} is
		\[
		\geq \frac{|\Psi_M|^2|\Psi_{M,N}|^2}{4 \Psi_{M, M}^2 \Psi_{N, N}} > \frac{\delta}{4} \beta(M),
		\]
		and (i) follows.
	\end{proof}

	\section{Reducing the proof of Theorem~\ref{thm:ISOptimality} to computing mollified moments}\label{sec:reductionIS}
	
	In this section, we use Theorem \ref{thm:criterion1} to reduce the proof of Theorem \ref{thm:ISOptimality} to a computation of some mollified moments. For simplicity, we shall concentrate on even primitive characters and denote summations over them by $\sum^+$ and their number by $\phi^+(q)$. Odd characters can be handled similarly. 
	
	Following the notation in Theorem \ref{thm:ISOptimality} and Section \ref{sec:general mollifying},  for arbitrary  mollifiers $M, N$, we write
	\begin{align} \label{eq:Psiqdefs}
		\begin{aligned}
			\Psi_M(q) &= \sumplus_{\chi \mod{q}} L(1/2, \chi) M(\chi), \\
			\Psi_{M, N}(q) &= \sumplus_{\chi \mod{q}} L(1/2, \chi) M(\chi) \overline{L(1/2, \chi) N(\chi)}.
		\end{aligned}
	\end{align}
	Let $M_{\mathrm{IS}}$ and $N_{\mathrm{G}}$ be as in Theorem~\ref{thm:ISOptimality}. To be able to handle some error terms efficiently (see Remark~\ref{rem:increase} below), we work also with a slightly longer Iwaniec-Sarnak mollifier of length $Q^{\theta+\varepsilon_0}$ with a small $\varepsilon_0$. To facilitate this, we prove our moment estimates for mollifiers of general lengths $y_1$ and $y_2$. Let us state the relevant moment estimates concerning
	\begin{equation} \label{eq:ISMNdefs}
		M(\chi) := \sum_{b \leq y_1} \frac{\mu(b) \chi(b) (1-\frac{\log b}{\log y_1})}{\sqrt{b}} \quad \text{and}\quad N(\chi) := \sum_{b \leq y_2} \frac{x_b \chi(b)}{\sqrt{b}}
	\end{equation}
	and then deduce Theorem~\ref{thm:ISOptimality} from these and Theorem~\ref{thm:criterion1}.
	\begin{proposition}\label{prop:Psi_N_1IS}
		Let $\theta \in (0, 1)$ be fixed and let $q \geq 3.$ Then the following hold.
		\begin{enumerate}[(i)]
			\item Let $y_2 \leq q^{\theta}$ and let $N(\chi)$ be as in~\eqref{eq:ISMNdefs} with arbitrary coefficients $x_{b} \in \mathbb{C}$ such that $x_{b} \ll_\epsilon q^\epsilon$ for every $\varepsilon > 0$. We have, for every $\varepsilon > 0$,
			\begin{equation*}
				\Psi_{N}(q) =\phi^+(q) x_1 + O(q^{\theta/2+1/2+\varepsilon}).
			\end{equation*}
			\item Let $y_1 \leq q^{\theta}$ and let $M(\chi)$ be as in~\eqref{eq:ISMNdefs}. We have, for every $\varepsilon > 0$,
			\begin{equation*}
				\Psi_M(q)= \varphi^+(q) +O(q^{\theta/2+1/2+\varepsilon}).
			\end{equation*}
		\end{enumerate}
	\end{proposition}
	
	\begin{proposition}\label{prop:PsiMN_1IS}
		Let $\theta \in (0, 1/2)$ be fixed. Let $q \geq 3$. Then the following hold.
		\begin{enumerate}[(i)]
			\item Let $\varepsilon_0 > 0$ be fixed, let $y_1 \in [q^{2\varepsilon_0}, q^{\theta}]$, and $y_2 \leq y_1 q^{-\varepsilon_0}$. Let $M(\chi)$ and $N(\chi)$ be as in~\eqref{eq:ISMNdefs} with arbitrary coefficients $x_{b} \in \mathbb{C}$ such that $x_{b} \ll_\epsilon q^\epsilon$ for every $\varepsilon > 0$. Then
			\begin{equation*}
				\Psi_{M,N}(q)= \overline{x_1} \lz1+\frac{\log q}{\log y_1}\pz \varphi^+(q)+o\left(\sqrt{\Psi_{N, N}(q) \varphi^+(q)}+\varphi^+(q)\right).
			\end{equation*}
			\item Let $y_1 \leq q^{\theta}$ and let $M(\chi)$ be as in~\eqref{eq:ISMNdefs}. Then
			\begin{equation*}
				\Psi_{M,M}(q)= \left(1+\frac{\log q}{\log y_1}\right)\varphi^+(q) +o(\varphi^+(q)).
			\end{equation*}
		\end{enumerate}
	\end{proposition}
	Propositions~\ref{prop:Psi_N_1IS} and~\ref{prop:PsiMN_1IS}(i) follow immediately from the more general Propositions~\ref{prop:Psi_N_1} and~\ref{prop:PsiM1N1gcd1} below (taking $z_{a, b} = \mathbf{1}_{a = 1} \mathbf{1}_{(b, q) = 1} x_b$)  that we will prove in Sections~\ref{sec:PsiN1} and~\ref{sec:PsiM1N1}. Proposition~\ref{prop:PsiMN_1IS}(ii) is well-known and we prove it in Section~\ref{sec:PsiM1N1} referring to the literature.
	
	\begin{proof}[Proof of Theorem~\ref{thm:ISOptimality} assuming Propositions~\ref{prop:Psi_N_1IS} and~\ref{prop:PsiMN_1IS}]
		First notice that Propositions~\ref{prop:Psi_N_1IS}(ii) and~\ref{prop:PsiMN_1IS}(ii) imply that $\beta_q(M_{\mathrm{IS}}) = \frac{1}{1+\frac{1}{\theta}} + o(1)$ as claimed in~\eqref{eq:I-Soptimalitybound}. Let now
		\[
		M_0(\chi) := \sum_{b \leq q^{\theta+\varepsilon_0}} \frac{\mu(b) \chi(b)}{\sqrt{b}} \left(1-\frac{\log b}{\log q^{\theta+\varepsilon_0}}\right)
		\]
		with a small positive constant $\varepsilon_0 < 1/2-\theta$. Then Propositions~\ref{prop:Psi_N_1IS}(ii) and~\ref{prop:PsiMN_1IS}(ii) also give that $\beta_q(M_0) = \frac{1}{1+\frac{1}{\theta+\varepsilon_0}} + o(1)$. Thus it suffices to show that, for every sufficiently small $\varepsilon_0 > 0$, we have $\beta_q(N_{\mathrm{G}}) \leq \beta_q(M_0) + o(1)$, whenever $N_{\mathrm{G}}$ is as in Theorem~\ref{thm:ISOptimality}.
		
		By Theorem \ref{thm:criterion1}, this follows once we have shown that
		\begin{equation*}
			|\Psi_{M_0}(q) \Psi_{M_0, N_{\mathrm{G}}}(q)| = |\Psi_{N_{\mathrm{G}}}(q) \Psi_{M_0, M_0}(q)| + o(|\Psi_{N_{\mathrm{G}}}(q)| \Psi_{M_0, M_0}(q)),
		\end{equation*}
		Propositions~\ref{prop:Psi_N_1IS} and~\ref{prop:PsiMN_1IS} together with the assumptions of Theorem~\ref{thm:ISOptimality} imply that $|\Psi_{N_{\mathrm{G}}}(q)| \Psi_{M_0, M_0}(q) \gg \varphi^+(q)^2$. Thus it suffices to establish that
		\begin{equation*}
			|\Psi_{M_0}(q) \Psi_{M_0, N_{\mathrm{G}}}(q)| = |\Psi_{N_{\mathrm{G}}}(q) \Psi_{M_0, M_0}(q)| + o(\varphi^+(q)^2).
		\end{equation*}
		But this follows by substituting the moment computations from Propositions~\ref{prop:Psi_N_1IS} and~\ref{prop:PsiMN_1IS}, and noting that we may assume $\Psi_{N_G,N_G}\ll \phi^+(q)$, as otherwise the claim  $\beta_q(N_{\mathrm{G}}) \leq \beta_q(M_0)+o(1)$ is trivial.
	\end{proof}
	
	\begin{remark}\label{rem:IS} 
		Although Iwaniec and Sarnak~\cite{IS} (see e.g. the discussion above (6.10) and also the end of Section 6) indicate that their choice of mollifier is optimal, the argument provided in~\cite[Section 6]{IS} (which is motivated by Selberg's works such as~\cite{Selberg1, Selberg2, Selberg3}) does not seem to prove optimality --- there the inequality~\cite[(6.10)]{IS} is used to diagonalize the main term of the second moment, and then the diagonalized main term is minimized assuming the first moment is $1$. Finally, it is noticed that with this optimized choice, there is equality in the inequality~\cite[(6.10)]{IS}. However, it seems that this argument does not necessarily rule out the existence of a different mollifier for which the inequality in~\cite[(6.10)]{IS} is strict and thus has a possibly smaller second moment. However, the discussion above~\cite[(6.10)]{IS} indicates that the authors of~\cite{IS} have also a rigorous proof of the optimality, at least in a "reasonable" class of mollifiers\footnote{Private communication indicates that this reasonable class of mollifiers is smaller than the class we consider in Theorem~\ref{thm:ISOptimality}.}.
	\end{remark}

	\section{Reducing the proof of Theorem~\ref{thm:main general result DirL} to computing mollified moments}\label{sec:reductionMV}
	
	In this section, we use Theorem \ref{thm:criterion1} to reduce the proof of Theorem \ref{thm:main general result DirL} to a computation of some mollified moments. For simplicity, we shall again concentrate on even primitive characters.
	
	Following the notation in Theorem \ref{thm:main general result DirL} and Section \ref{sec:general mollifying}, we write, for mollifiers $M, N$,
	\begin{align*}
		\Psi_M &= \frac{1}{Q^2}\sum_{q\geq 1}\Phi\bfrac{q}{Q}\frac{q}{\phi(q)} \Psi_M(q), \\ 
		\Psi_{M, N} &= \frac{1}{Q^2}\sum_{q\geq 1}\Phi\bfrac{q}{Q}\frac{q}{\phi(q)} \Psi_{M, N}(q),
	\end{align*}
	where $\Psi_{M}(q)$ and $\Psi_{M, N}(q)$ are as in~\eqref{eq:Psiqdefs} and $\Phi(x)$ is as in Theorem \ref{thm:main general result DirL}. For simplicity and convenience, the normalization is somewhat different from Section \ref{sec:general mollifying}.
	
	Similarly to the previous section, the error estimates become easier when we work with a slightly longer Michel-Vanderkam mollifier, and instead of $M_{\mathrm{MV}}$ we use
	\begin{equation} \label{eq:M0MVdef}
		M_0(\chi) := \sum_{b \leq Q^{\theta+\varepsilon_0}} \frac{\mu(b) \chi(b) (1-\frac{\log b}{\log Q^{\theta+\varepsilon_0}})}{\sqrt{b}} + \overline{\varepsilon_\chi} \sum_{b \leq Q^{\theta+\varepsilon_0}} \frac{\mu(b) \overline{\chi}(b) (1-\frac{\log b}{\log Q^{\theta+\varepsilon_0}})}{\sqrt{b}}.
	\end{equation}
	The starting point of the proof of Theorem~\ref{thm:main general result DirL} will be showing that it suffices to show that
	\begin{equation} \label{eq:genLclaim}
		|\Psi_{M_0} \Psi_{M_0, N_{\mathrm{B}}}| = |\Psi_{N_{\mathrm{B}}} \Psi_{M_0, M_0}| + o(|\Psi_{N_{\mathrm{B}}}| \Psi_{M_0, M_0}),
	\end{equation}
	whenever $N_{\mathrm{B}}$ is as in Theorem~\ref{thm:main general result DirL}. Before turning to the proof of Theorem~\ref{thm:main general result DirL}, we simplify this claim and state propositions concerning relevant mollified moments.
	
	As in the statement of Theorem \ref{thm:main general result DirL}, we let $z_{a,b}=x_{a,b}+y_{a,b}$, and define
	\begin{equation}\label{eq:N0def}
		N_0(\chi) := \sum_{ab \leq Q^\theta} \frac{z_{a, b}\overline{\chi}(a)\chi(b)}{\sqrt{ab}}.
	\end{equation}
	Note that by the functional equation $\overline{\varepsilon_\chi} L(1/2, \chi) = L(1/2, \overline{\chi})$. Thus
	\begin{align*}
		\Psi_{N_{\mathrm{B}}}(q) &= \sumplus_{\chi \mod{q}} L(1/2, \chi) \sum_{ab \leq Q^\theta} \frac{x_{a, b}\overline{\chi}(a)\chi(b)}{\sqrt{ab}} + \sumplus_{\chi \mod{q}} L(1/2, \chi) \overline{\varepsilon_\chi} \sum_{ab \leq Q^\theta} \frac{y_{a, b}\chi(a)\overline{\chi}(b)}{\sqrt{ab}} \\
		&= \sumplus_{\chi \mod{q}} L(1/2, \chi) \sum_{ab \leq Q^\theta} \frac{x_{a, b}\overline{\chi}(a)\chi(b)}{\sqrt{ab}} + \sumplus_{\chi \mod{q}} L(1/2, \overline{\chi}) \sum_{ab \leq Q^\theta} \frac{y_{a, b}\chi(a)\overline{\chi}(b)}{\sqrt{ab}}.
	\end{align*}
	In the second sum above, we can now change variable $\chi \to \overline{\chi}$ and see that
	\begin{equation} \label{eq:PsiNB=PsiN0}
		\Psi_{N_{\mathrm{B}}}(q) = \Psi_{N_0}(q).
	\end{equation}
	Furthermore
	\begin{align*}
		\Psi_{M_0, N_{\mathrm{B}}}(q) &= \sumplus_{\chi \mod{q}} |L(1/2, \chi)|^2 M_0(\chi) \sum_{ab \leq Q^\theta} \frac{\overline{x_{a, b}}\chi(a)\overline\chi(b)}{\sqrt{ab}} \\
		& + \sumplus_{\chi \mod{q}} |L(1/2, \chi)|^2 M_0(\chi)  \varepsilon_\chi \sum_{ab \leq Q^\theta} \frac{\overline{y_{a, b}} \, \overline\chi(a)\chi(b)}{\sqrt{ab}}.
	\end{align*}
	Noting that $\varepsilon_\chi M_0(\chi) = M_0(\overline{\chi})$ and making the change of variable $\chi \to \overline{\chi}$, we see that
	\begin{equation*}
		\Psi_{M_0,N_{\mathrm{B}}}(q)=\Psi_{M_0,N_0}(q).
	\end{equation*}
	Hence~\eqref{eq:genLclaim} is equivalent to
	\begin{equation}\label{eq:goal}
		|\Psi_{M_0} \Psi_{M_0, N_0}| = |\Psi_{N_0} \Psi_{M_0, M_0}| + o(|\Psi_{N_0}| \Psi_{M_0, M_0}).
	\end{equation}
	
	To show this, we shall prove the following propositions concerning the mollified moments with mollifiers
	\begin{equation} \label{eq:MdefMV}
		M(\chi) := \sum_{b \leq y_1} \mu(b) \frac{\chi(b)}{\sqrt{b}} \left(1-\frac{\log b}{\log y_1}\right) + \overline{\varepsilon_\chi} \sum_{b \leq y_1} \mu(b) \frac{\overline{\chi}(b)}{\sqrt{b}} \left(1-\frac{\log b}{\log y_1}\right)
	\end{equation}
	and
	\begin{equation} \label{eq:NdefMV}
		N(\chi) = \sum_{ab \leq y_2} \frac{z_{a, b}\overline{\chi}(a)\chi(b)}{\sqrt{ab}}.
	\end{equation}
	In the propositions below, we only make the necessary assumptions about $z_{a, b}$, so some of them work also for complex $z_{a, b}$. 
	\begin{proposition}\label{prop:Psi_N_1}
		Let $\theta \in (0, 1)$ be fixed and let $q \geq 3.$ Then the following hold.
		\begin{enumerate}[(i)]
			\item  Let $y_2 \leq q^{\theta}$ and let $N(\chi)$ be as in~\eqref{eq:NdefMV} with arbitrary coefficients $z_{a, b} \in \mathbb{C}$ such that $z_{a, b} \ll_\epsilon q^\epsilon$ for every $\varepsilon > 0$. We have, for every $\varepsilon > 0$,
			\begin{equation*}
				\Psi_{N}(q) =\phi^+(q) \sum_{\substack{kb^2\leq y_2 \\ (kb, q) = 1}}\frac{z_{kb,b}}{kb} + O(q^{\theta/2+1/2+\varepsilon}).
			\end{equation*}
			\item Let $y_1 \leq q^{\theta}$ and let $M(\chi)$ be as in~\eqref{eq:MdefMV}. We have, for every $\varepsilon > 0$,
			\begin{equation*}
				\Psi_M(q)=2 \phi^+(q)+O(q^{\theta/2+1/2+\varepsilon}).
			\end{equation*}
		\end{enumerate}
	\end{proposition}
	
	We will quickly establish this Proposition in Section~\ref{sec:PsiN1}: Part (i) was essentially proved by Bui, and part (ii) follows easily from part (i) (or the work of Michel and Vanderkam~\cite{M-V}).
	
	\begin{proposition}\label{prop:PsiMN_1}
		Let $\theta \in (0, 1/2)$ be fixed. Let $Q \geq 3$ and $C_0 := \frac{\Gamma'(1/4)}{\Gamma(1/4)}-\log \pi$. Then the following hold.
		\begin{enumerate}[(i)]
			\item Let $\varepsilon_0 > 0$ be fixed, let $y_1 \in [q^{2\varepsilon_0}, q^{\theta}]$ and $y_2 \leq y_1 q^{-\varepsilon_0}$. Let $M(\chi)$ and $N(\chi)$ be as in~\eqref{eq:MdefMV} and~\eqref{eq:NdefMV} with arbitrary coefficients $z_{a, b} \in \mathbb{R}$ satisfying the conditions (a) and (b) from Theorem~\ref{thm:main general result DirL} for a sufficiently small constant $c > 0$. 
			
			Then
			\begin{equation*}
				\Psi_{M,N}=\frac{1}{Q^2}\sum_{q\geq 1}\Phi\bfrac{q}{Q}\frac{q \phi^+(q)}{\phi(q)} \sum_{\substack{kb^2\leq y_2 \\ (kb, q) = 1}}\frac{z_{kb,b}}{kb}\lz2+\frac{\log q}{\log y_1} + \frac{C_0}{\log y_1}\pz+o(1).
			\end{equation*}
			Assuming the quasi-Riemann hypothesis, we do not need to assume the condition (b) from Theorem~\ref{thm:main general result DirL}.
			\item We have
			\begin{equation*}
				\Psi_{M,M}=\frac{1}{Q^2} \sum_{q \geq 1} \Phi\left(\frac{q}{Q}\right) \frac{q\phi^+(q)}{\varphi(q)} \left(4+2\frac{\log q}{\log y_1}\right) +o(1).
			\end{equation*}
		\end{enumerate}
	\end{proposition}
	
	Before discussing the proof of Proposition~\ref{prop:PsiMN_1}, let us deduce Theorem~\ref{thm:main general result DirL} (similarly as we deduced Theorem~\ref{thm:ISOptimality}).
	\begin{proof}[Proof of Theorem~\ref{thm:main general result DirL} assuming Propositions~\ref{prop:Psi_N_1} and~\ref{prop:PsiMN_1}]
		Propositions~\ref{prop:Psi_N_1}(ii) and~\ref{prop:PsiMN_1}(ii) immediately imply that $\beta_Q(M_{\mathrm{MV}}) = \frac{1}{1+\frac{1}{2\theta}} + o(1)$. Furthermore, with $M_0$ as in~\eqref{eq:M0MVdef}, we obtain similarly that $\beta_Q(M_0) = \frac{1}{1+\frac{1}{2(\theta + \varepsilon_0)}} + o(1)$. Thus it suffices to show that, for every sufficiently small $\varepsilon_0 > 0$ and every $N_\mathrm{B}$ as in Theorem~\ref{thm:main general result DirL} we have $\beta_Q(N_{\mathrm{B}}) \leq \beta_Q(M_0) + o(1)$. By Theorem~\ref{thm:criterion1} this follows once we have shown~\eqref{eq:genLclaim}, and we already showed above that this follows from~\eqref{eq:goal}. By Propositions~\ref{prop:Psi_N_1} and~\ref{prop:PsiMN_1} and the assumptions of Theorem~\ref{thm:main general result DirL} we see that $|\Psi_{N_0}| \Psi_{M_0, M_0} \gg 1$. Thus it suffices to establish that
		\begin{equation*}
			|\Psi_{M_0} \Psi_{M_0, N_0}| = |\Psi_{N_0} \Psi_{M_0, M_0}| + o(1).
		\end{equation*}
		But this follows by substituting the moment computations from Propositions~\ref{prop:Psi_N_1} and~\ref{prop:PsiMN_1} and the assumptions of Theorem \ref{thm:main general result DirL}.
	\end{proof}
	
	Now we discuss the proof of Proposition~\ref{prop:PsiMN_1} (Part (ii) would follow also from the works of Michel and Vanderkam~\cite{M-V} and Pratt~\cite{Pratt} --- see \cite[Section 3]{Pratt}). We write
	\begin{equation} \label{eq:M12defMV}
		M_1(\chi) := \sum_{b \leq y_1} \mu(b) \frac{\chi(b)}{\sqrt{b}} \left(1-\frac{\log b}{\log y_1}\right) \quad \text{and} \quad  M_2(\chi) := \overline{\varepsilon_\chi} \sum_{b \leq y_1} \mu(b) \frac{\overline{\chi}(b)}{\sqrt{b}} \left(1-\frac{\log b}{\log y_1}\right).
	\end{equation}
	Then $M(\chi) = M_1(\chi) + M_2(\chi)$ and consequently $\Psi_{M,N}(q)=\Psi_{M_1,N}(q)+\Psi_{M_2,N}(q)$ and $\Psi_{M, M}(q) = 2\Psi_{M_1, M_1}(q)+ 2\re \Psi_{M_1, M_2}(q)$. The following propositions concerning $\Psi_{M_1, N}(q)$ and $\Psi_{M_2, N}$ are proved in Sections~\ref{sec:PsiM1N1} and~\ref{sec:PsiM2N1} adapting the methods of Iwaniec and Sarnak~\cite{IS} and Pratt~\cite{Pratt}.
	
	\begin{proposition}\label{prop:PsiM1,N1}
		Let $\theta \in (0, 1/2)$ and $\varepsilon_0 > 0$ be fixed, let $C_0 = \frac{\Gamma'(1/4)}{\Gamma(1/4)}-\log \pi$, and let $q \geq 3$.  Let $y_1 \in [q^{2\varepsilon_0}, q^{\theta}]$ and $y_2 \leq y_1 q^{-\varepsilon_0}$. Let $M_1(\chi)$ and $N(\chi)$ be as in~\eqref{eq:M12defMV} and~\eqref{eq:NdefMV} with arbitrary coefficients $z_{a, b} \in \mathbb{C}$ satisfying the conditions (a) and (b) from Theorem~\ref{thm:main general result DirL} for a sufficiently small constant $c > 0$. 
		
		Then
		\begin{equation*}
			\Psi_{M_1,N}(q)=\phi^+(q) \sum_{\substack{kb^2\leq y_2 \\ (kb, q) = 1}}\frac{\overline{z_{kb,b}}}{kb}\lz1+\frac{\log q}{\log y_1}-\frac{\log k}{\log y_1}+\frac{C_0}{\log y_1}\pz+ o\lz \phi^+(q)\pz.
		\end{equation*}
		Assuming the quasi-Riemann hypothesis, we do not need to assume the condition (b) from Theorem~\ref{thm:main general result DirL}.
	\end{proposition}
	
	\begin{proposition}\label{prop:PsiM_2N_1}
		Let $\theta \in (0, 1/2)$ be fixed, let $Q \geq 3$, and let $y_2 \leq y_1 \leq Q^{\theta}$. Let $M_2(\chi)$ and $N(\chi)$ be as in~\eqref{eq:M12defMV} and~\eqref{eq:NdefMV} with arbitrary coefficients $z_{a, b} \in \mathbb{R}$ such that $|z_{a, b}|\ll_\epsilon q^\epsilon$ for every $\varepsilon > 0$. Then
		\begin{equation*}
			\Psi_{M_2,N}=\frac{1}{Q^2}\sum_{q\geq 1}\Phi\bfrac{q}{Q}\frac{q\phi^+(q)}{\phi(q)} \sum_{\substack{kb^2\leq y_2 \\ (kb, q) = 1}}\frac{z_{kb,b}}{kb}\lz1+\frac{\log k}{\log y_1}\pz+o(1).
		\end{equation*}
	\end{proposition}
	Proposition~\ref{prop:PsiMN_1}(i) follows immediately from Propositions~\ref{prop:PsiM1,N1} and~\ref{prop:PsiM_2N_1} whereas Proposition~\ref{prop:PsiMN_1}(ii) follows from Proposition~\ref{prop:PsiMN_1IS}(ii) and Proposition~\ref{prop:PsiM_2N_1} with $z_{a, b} = 1_{a=1}\mu(b)P[b]$.

	\section{Auxiliary results}\label{sec:prelim}
	
	\subsection{Orthogonality of even characters and root numbers}
	
	For $(mn,q)=1$, we have the orthogonality relations (see for instance \cite[Section 3]{IS})
	\begin{equation}\label{eq:orthogonality}
		\sumplus_{\chi \mod q}\chi(m)\overline\chi(n)=\frac12\sum_{\substack{vw=q \\ w \mid m \pm n}}\mu(v)\phi(w),
	\end{equation} where the sum on the right should be understood as a sum of two sums, one with $w \mid m+n$ and one with $w \mid m-n$ (later we shall use similar convention for sums over $m\equiv \pm n\mod w$). We also denote by $\phi^+(q)$ the number of primitive even characters modulo $q$, and observe from above with $m = n = 1$ that
	\begin{equation} \label{eq:evencharscount}
		\phi^+(q)=\tfrac12(\mu*\phi)(q) + O(1).
	\end{equation}
	
	We shall further need the following orthogonality relation with the root numbers, which again holds whenever $(mn,q)=1$ (see for example \cite[Section 3]{IS}):
	\begin{align}
		\label{eq:orthogwitheps}
		\sumplus_{\chi \mod{q}} \varepsilon_\chi \chi(m) \overline{\chi}(n) = \frac{1}{q^{1/2}} \sum_{\substack{vw = q \\ (v, w) = 1}} \mu^2(v) \varphi(w) \cos\left(\frac{2\pi n \overline{mv}}{w}\right).
	\end{align}
	Orthogonality of all characters, together with an expression of $\epsilon_\chi$ as a normalized Gauss sum gives for any $(mn,w)=1$ that
	\begin{equation}\label{eq:orthogonality of all characters with root number}
		\sum_{\chi\mod w}\overline\epsilon_\chi\chi(m)\overline\chi(n)=\frac{\phi(w)}{\sqrt w}e\bfrac{m\overline n}{w}.
	\end{equation}
	
	\subsection{Sums with the M\"obius function}
	In order to evaluate $\Psi_{M_1, N_0}(q)$ we need a more precise version of a special case of a standard lemma concerning M\"obius sums (for the original version, see Conrey's paper \cite[Lemma 10]{Conrey}).
	
	In its proof we will use a variant of Perron's formula, which is both truncated and contains a smooth weight. Mimicking the proof of \cite[Lemma 7.1]{Kou}, we find that for $y,c>0$, $T>2$,
	$$
	\int_{c-iT}^{c+iT}\frac{y^{-s}}{s^2}ds=\begin{cases}
		-\log y+O\left( \frac{y^{-c}}{T\max\{1,T|\log y|\}}\right),&\text{if $0<y<1$;}\\
		O\left(\frac{y^{-c}}{T\max\{1,T|\log y|\}}\right),&\text{if $y>1$.}
	\end{cases}
	$$ Using this with $y = nj/x$ in the proof of \cite[Theorem 7.2]{Kou}, we obtain the following lemma.
	\begin{lemma}\label{le:Perron}
		Let $\varepsilon > 0$. Let $f(n)$ be an arithmetic function with $|f(n)|\ll 1$ and $F(s) = \sum_{n \in \mathbb{N}} \frac{f(n)}{n^s}$ be its Dirichlet series. Let $x,T\geq 2$, let $j \in [1, x)$, and let $c=1/\log (x/j)$. Then
		$$
		\sum_{n\leq x/j}\frac{f(n)}{n}\lz1-\frac{\log (nj)}{\log x}\pz=\frac{1}{2\pi i\log x}\int_{c-iT}^{c+iT} F(s+1)\frac{(x/j)^s}{s^2}ds+O\lz\frac{1}{T^2 }+\frac1{T \cdot \frac{x}{j} \log x}\pz.
		$$
	\end{lemma}
	We shall encounter the function $\eta \colon \mathbb{N} \to \mathbb{R}$ defined via
	\begin{align*}
		\eta(d) := -\frac{1}{\varphi(d)} \sum_{ac = d} \mu(a) c \log a.
	\end{align*}
	Note that
	\begin{align*}
		\eta(d) &= -\frac{d}{\varphi(d)} \sum_{a \mid d} \frac{\mu(a) \log a}{a} = -\frac{d}{\varphi(d)} \sum_{a \mid d} \frac{\mu(a) \sum_{p \mid a} \log p}{a} = -\frac{d}{\varphi(d)} \sum_{p \mid d} \log p \sum_{\substack{a \mid d \\ p \mid a}} \frac{\mu(a)}{a} \\
		&= -\frac{d}{\varphi(d)} \sum_{p \mid d} \log p \sum_{\substack{a \mid \frac{d}{p}}} \frac{\mu(pa)}{pa} = \frac{d}{\varphi(d)} \sum_{p \mid d} \frac{\log p}{p} \sum_{\substack{\substack{a \mid d \\ (a, p) = 1}}} \frac{\mu(a)}{a} = \sum_{p \mid d} \frac{\log p}{p-1}.
	\end{align*}
	Hence we see that $\eta(d)$ is additive and 
	
	\begin{align}\label{align:etaformulas}
		\eta(d) = -\frac{1}{\varphi(d)} \sum_{ac = d} \mu(a) c \log a =  \sum_{p \mid d} \frac{\log p}{p-1} .
	\end{align}

	\begin{lemma}\label{lemma: Conrey}
		Let $\varepsilon > 0$ and let $c > 0$ be sufficiently small in terms of $\varepsilon$. Let $y \geq 2$ and let $j, q \in \mathbb{N}$ with $j\leq y^{1-\varepsilon}$ and $q = y^{O(1)}$. Then the following hold.
		\begin{enumerate}[(i)]
			\item We have
			\begin{equation*}
				\sum_{\substack{n\leq y/j,\\(n,jq)=1}}\frac{\mu(n)}{n}\lz1-\frac{\log (jn)}{\log y}\pz=\frac{jq}{\phi(jq)\log y}+O\left(\exp\left(-c\frac{(\log y)^{3/5}}{(\log \log y)^{1/5}}\right)\right).
			\end{equation*}
			\item We have
			
			\begin{equation}\label{eqn: Conrey lemma 2}
				\begin{aligned}
					\sum_{\substack{n\leq y/j,\\(n,jq)=1}}\frac{-\mu(n)\log n}{n}\lz1-\frac{\log (jn)}{\log y}\pz &=\frac{jq}{\phi(jq)\log y}\lz\log\frac yj-2\gamma-2\eta(jq)\pz\\
					&+O\left(\exp\left(-c\frac{(\log y)^{3/5}}{(\log \log y)^{1/5}}\right)\right).
				\end{aligned}   
			\end{equation}
			\item Assuming the quasi-Riemann hypothesis, the error terms in parts (i) and (ii) can be replaced by $y^{-\delta}$ for some $\delta > 0$ (depending on $\varepsilon$ and the zero-free region region in the quasi-Riemann hypothesis).
		\end{enumerate}                
	\end{lemma}
	
	In the proof we shall need, for $m \in \mathbb{N}$ and $s \in \mathbb{C}$, the function
	\begin{align*}
		F(m, s) := \prod_{p|m}\lz1-\frac{1}{p^s}\pz
	\end{align*}
	which we can bound using the following rough estimates.
	
	\begin{lemma}\label{le:ConreyError}
		For every $c \in [0, 1/10]$ and every integer $m\geq 5$, we have
		$$
		\frac{1}{F(m,1-c)} \ll \exp((\log m)^{1/5}) \quad \text{and} \quad \frac{F'(m,1-c)}{F(m,1-c)}\ll \log m,
		$$
	\end{lemma}
	\begin{proof}
		We can clearly assume that $m$ is sufficiently large. Notice first that since $\prod_{p \leq 2 \log m} p > m$ by the prime number theorem, for any decreasing $f \colon \mathbb{N} \to \mathbb{R}_{\geq 0}$, we have
		\[
		\sum_{p \mid m} f(p) \leq \sum_{p \leq 2\log m} f(p).
		\]
		Consequently
		\begin{align} \label{eq:logFest}
			\begin{aligned}
				& \log \left|F\left(m,1-c\right)\right|^{-1} = \sum_{p \mid m} \log \left|1-\frac{1}{p^{1-c}}\right|^{-1} \\
				&\leq \sum_{p|m}\frac{1}{p^{1-c}} +O(1) \leq \sum_{p \leq 2\log m} \frac{1}{p^{1-c}} + O(1) \ll (\log m)^{1/10}.
			\end{aligned}
		\end{align}
		and the first claim follows.
		
		To prove the second claim, we note that by differentiating $\log F(m,\sigma)$, we obtain
		\begin{equation}\label{eq:Flogder}
			\frac{F'(m,s)}{F(m,s)}=\frac{d}{ds}\sum_{p|m}\log\lz1-\frac{1}{p^s}\pz=\sum_{p|m}\frac{\log p}{p^s-1}.
		\end{equation}
		Hence, similarly to~\eqref{eq:logFest}, 
		\[
		\frac{F'(m,1-c)}{F(m,1-c)}\ll \sum_{p \leq 2\log m} \frac{\log p}{p^{1-c}} \ll \log m.
		\]
	\end{proof}
	
	\begin{proof}[Proof of Lemma~\ref{lemma: Conrey}]
		We can clearly assume that $y$ is sufficiently large.
		
		For $j \in \mathbb{N}$ and $s \in \mathbb{C}$, define $\zeta_{(j)}(s) := \zeta(s) F(j, s)$. Let $d > 0$ be a small constant. In particular we choose $d$ to be so small that, for any $T \geq 2$, $\zeta(s+1)$ has no zeroes with $|\im(s)| \leq T$ and $\re(s) > \frac{-2d}{(\log T)^{2/3} (\log \log T)^{1/3}}$ (for such a zero-free region, see e.g.~\cite[Theorem 6.1]{IvicBook}). In the proofs of (i) and (ii) we take
		\begin{equation} \label{eq:Tsigdefs}
			T := \exp\left(\frac{(\log y)^{3/5}}{(\log \log y)^{1/5}}\right) \text{ and }  \sigma:=\frac{d}{(\log T)^{2/3} (\log \log T)^{1/3}} \geq \frac{d}{(\log y)^{2/5} (\log \log y)^{1/5}}.
		\end{equation}
		
		\textbf{Proof of (i):} By Lemma \ref{le:Perron}, we have
		\begin{align*}
			\begin{aligned}
				\sum_{\substack{n\leq y/j,\\(n,jq)=1}}\frac{\mu(n)}{n}\lz1-\frac{\log (jn)}{\log y}\pz
				&=\frac{1}{2\pi i\log y}\int_{\frac{1}{\log(y/j)}-iT}^{\frac1{\log(y/j)}+iT}\frac{1}{\zeta_{(jq)}(s+1)}\frac{(y/j)^s}{s^2}ds \\
				& \qquad \qquad +O\lz \frac{1}{T^2} + \frac{1}{T \cdot y/j}\pz.
			\end{aligned}
		\end{align*}
		We now shift the integral to the contour with vertices at $-\sigma-iT, -\sigma-2i, -1/2-2i, -1/2+2i, -\sigma+2i,-\sigma+iT$. By the choice of $\sigma$, we cross only a simple pole at $s=0$. To calculate the residue, we note the following approximations for $s$ around $0$:
		\begin{equation*}
			\begin{aligned}
				\zeta_{(jq)}(s+1)&=\frac{F(jq,1)}{s}+O(1),\\
				\frac{1}{\zeta_{(jq)}(s+1)}&= \frac{s}{F(jq, 1) + O(|s|)} = \frac{s}{F(jq,1)}+O\left(\frac{|s|^2}{F(jq, 1)^2}\right).
			\end{aligned}
		\end{equation*}
		Hence the residue at $s=0$ equals $\frac{1}{F(jq,1)}=\frac{jq}{\phi(jq)}$ which leads to the claimed main term.
		
		On the new path of integration, we have $\frac{1}{|\zeta(s+1)|}\ll \log T$ (see for instance \cite[Lemma 12.3]{IvicBook}), so the integral between $-\sigma + 2i$ and $-\sigma + iT$ contributes
		$$
		\ll \frac{(y/j)^{-\sigma}}{\log y}\int_{2}^{T}\frac{1}{F(jq,1-\sigma)}\frac{\log T}{1+t^2}dt\ll\frac{y^{-\varepsilon \sigma} \log T}{F(jq,1-\sigma)\log y}.
		$$
		Other shifted integrals make similar or smaller contribution.  The claim follows for $c = \varepsilon d/2$ after applying Lemma~\ref{le:ConreyError} and~\eqref{eq:Tsigdefs}.
		
		\textbf{Proof of (ii):}  The proof is similar as part (i). By Lemma~\ref{le:Perron}, the left-hand side of \eqref{eqn: Conrey lemma 2} equals
		\begin{equation}\label{eq:ConreybPerron}
			\frac{1}{2\pi i\log y}\int_{\frac{1}{\log(y/j)}-iT}^{\frac{1}{\log(y/j)}+iT}\lz\frac{1}{\zeta_{(jq)}(s+1)}\pz'\cdot\frac{(y/j)^s}{s^2}ds+O\lz \frac{1}{T^2} + \frac{1}{T \cdot y/j}\pz.
		\end{equation}
		We now shift the integral to the same contour as in part (i). To compute the residue, we use the Laurent expansion
		\begin{equation*}
			\begin{aligned}
				\zeta_{(jq)}(s)&=\zeta(s)F(jq,s)\\
				&=\lz \frac{1}{s-1}+\gamma+E(s-1)\pz\lz F(jq,1)+(s-1)F'(jq,1)+E((s-1)^2)\pz\\
				&=\frac{F(jq,1)}{s-1}+\gamma F(jq,1)+F'(jq,1)+E(s-1),
			\end{aligned}
		\end{equation*}
		where $E((s-1)^m)$ means that the remaining terms in the Laurent expansion have order at least $m$. This gives
		\begin{equation}\label{eqn: Laurent series 1}
			\begin{aligned}
				\frac{1}{\zeta_{(jq)}(s)}&=\frac{1}{\frac{F(jq,1)}{s-1}\left(1+\gamma (s-1)+\frac{F'(jq,1)}{F(jq, 1)}(s-1) +E((s-1)^2)\right)} \\
				&= \frac{s-1}{F(jq,1)}\left(1-\gamma (s-1)-\frac{F'(jq,1)}{F(jq, 1)}(s-1) +E\left((s-1)^2\right)\right) \\
				&= \frac{s-1}{F(jq,1)}+\frac{-\gamma F(jq,1)-F'(jq,1)}{F(jq,1)^2}(s-1)^2+E\left((s-1)^3\right).
			\end{aligned}
		\end{equation} 
		Differentiating \eqref{eqn: Laurent series 1} gives
		\begin{equation*}
			\lz\frac{1}{\zeta_{(jq)}(s+1)}\pz'=\frac{1}{F(jq,1)}-2s\frac{\gamma F(jq,1)+F'(jq,1)}{F(jq,1)^2}+ E(s^2),
		\end{equation*} so the final Laurent expansion of the integrand in~\eqref{eq:ConreybPerron} is
		\begin{equation*}
			\begin{aligned}
				&\frac{1}{s^2}\lz1+s\log(y/j)+E(s^2)\pz\lz\frac{1}{F(jq,1)}-2s\frac{\gamma F(jq,1)+F'(jq,1)}{F(jq,1)^2}+E\lz s^2\pz\pz\\
				&=\frac{1}{s^2F(jq,1)}+\frac{1}{s F(jq,1)}\cdot\lz\log (y/j)-2\gamma-2F'/F(jq,1)\pz+E(s^0).
			\end{aligned}
		\end{equation*}
		To obtain the claimed main term from the residue, we notice that, by~\eqref{eq:Flogder} and~\eqref{align:etaformulas},
		\begin{equation*}
			\frac{F'(jq,1)}{F(jq,1)}=\sum_{p|jq}\frac{\log p}{p-1}=\eta(jq).
		\end{equation*}
		To bound the error term similarly as in part (i), we note that
		$$
		\bfrac{1}{\zeta_{(jq)}(s+1)}'=\frac{-1}{\zeta(s+1)}\frac{F'(jq,s+1)}{F(jq,s+1)^2}-\frac{\zeta'(s+1)}{\zeta(s+1)^2}\frac{1}{F(jq,s+1)},
		$$ and use that on the new line of integration $1/\zeta(s+1),\frac{\zeta'(s+1)}{\zeta(s+1)}\ll (\log T)^3$ (see for instance~\cite[Lemmas 12.3 and 12.4]{IvicBook}.
		
		\textbf{Proof of (iii):} The proof of (iii) follows completely similarly, except now there exists $d > 0$ such that $\zeta(s+1)$ has no zeroes with $\re(s) > -2d$, and we can choose $\sigma = \varepsilon_0 d$ and $T=y^\sigma$ for a sufficiently small $\varepsilon_0 > 0$. When $\re(s) = -\sigma$ and $\im(s) \leq T$, we can bound $\zeta'(s+1)/\zeta(s+1)$ by $\log T$ using e.g.~\cite[Theorem 9.6(A)]{Titchmarsh}. To bound $1/\zeta(s+1)$ we can use e.g.~\cite[Theorem 9.6(B)]{Titchmarsh} to  get the bound $|\im(s)|^{O(1)}$ for $\re(s) \geq -d$. On the other hand, we have, for any $\varepsilon_1 > 0$ and $\re(s) \geq \varepsilon_1$, the bound $1/\zeta(s+1) = O(1)$. We can then use a convexity argument (e.g. Hadamard three circle theorem) to obtain the bound $1/\zeta(s+1) \ll |\im(s)|^{1/10}$ for $\re(s) \geq -d \varepsilon_0$ once $\varepsilon_0$ is sufficiently small. Using these bounds we obtain (iii) with $\delta = d \varepsilon_0 \varepsilon /2$.
		
	\end{proof}
	\subsection{A mean value result}
	We will use the following lemma, which is a consequence of the orthogonality of characters. 
	\begin{lemma}\label{le:general orthogonality}
		Let $A,B\geq 1$, and $\alpha_a, \beta_b$ be divisor-bounded sequences of complex numbers, that is $\alpha_n, \beta_n \ll \tau_k(n)$ for some $k\geq 1$, where $\tau_k(n)$ denotes the $k$-fold divisor function. Then we have
		\begin{equation}\label{eq:mean value lemma}
			\begin{aligned}
				\sum_{\chi\mod q}\lab\sum_{\substack{a \leq A\\ b\leq B}} \alpha_a\overline{\chi}(a)\beta_b \chi(b) \rab^2\ll (\log AB)^{O_{k}(1)}\lz q+AB\pz AB.
			\end{aligned}
		\end{equation}
	\end{lemma}
	\begin{proof}
		Let $k$ be such that $|\alpha_n|,|\beta_n|\ll \tau_k(n)$. Expanding the square and using the orthogonality of characters, the left-hand side of \eqref{eq:mean value lemma} is
		$$
		\phi(q)\sum_{\substack{a_1,a_2\leq A,\\ b_1,b_2\leq B,\\ b_1a_2\equiv a_1b_2\mod q}}\alpha_{b_1} \beta_{a_2} \overline{\alpha_{a_2} \beta_{b_1}}\ll \phi(q)\lz\sum_{\substack {u,v\leq AB\\ u\equiv v\mod q}}\tau_{2k}(u)\tau_{2k}(v)\pz.
		$$
		
		Using the inequality $|ab|\ll |a|^2+|b|^2 $ and symmetry, the above can be bounded by
		$$
		\ll \phi(q)\sum_{u\leq AB}\tau_{2k}^2(u)\sum_{\substack{v\leq AB \\ v\equiv u\mod q}}1\ll \phi(q)\sum_{u\leq AB}\tau_{2k}^2(u)\lz 1+\frac{AB}{q}\pz,
		$$ and the result follows upon estimating the divisor sum.
	\end{proof}
	
	\subsection{A bound for sums of incomplete Kloosterman sums} To estimate the error term of the cross term $\Psi_{M_2,N_0}$, we shall need the following result, which is a direct consequence of the work of Deshouillers and Iwaniec \cite[Theorem 12]{DI} (with $S=1$).
	\begin{lemma}
		\label{le:Klo}
		Let $C, D, N, R \geq 1/2$ and let $b_{n, r} \in \mathbb{C}$. Let $g \colon \mathbb{R}^2 \to \mathbb{R}$ be a smooth compactly supported function such that
		\begin{equation*}
			\left|\frac{\partial^{\nu_1 + \nu_2}}{\partial x_1^{\nu_1} \partial x_2^{\nu_2}} g(x_1, x_2)\right| \ll_{\nu_1, \nu_2} 1 \quad \quad \text{for every $\nu_1, \nu_2 \geq 0.$}
		\end{equation*}
		Then, for any $\varepsilon > 0$,
		\begin{align*}
			&\sum_{\substack{R < r \leq 2R}} \sum_{0 < n \leq N} b_{n, r}\sum_{\substack{c, d \\ (c, rd) = 1}} g\left(\frac{c}{C}, \frac{d}{D}\right) e\left(n \frac{\overline{rd}}{c}\right) \\
			&\ll (CDNR)^\varepsilon \left(C(R+N)(C+DR)+C^2D\sqrt{(R+N)R} + D^2NR\right)^{1/2} \\
			&\quad\cdot \left(\sum_{\substack{0 < n \leq N \\ R < r \leq 2R}} |b_{n, r}|^2\right)^{1/2}.
		\end{align*}
	\end{lemma}
	Note that there are some recent refinements of such bounds, see in particular Pascadi's work~\cite[Corollary 5.7]{Pascadi}. However, the classical bound of Deshouillers and Iwaniec is sufficient for our purposes.
	
	\section{Proof of Proposition \ref{prop:Psi_N_1}}
	\label{sec:PsiN1}
	We have
	\begin{align*}
		\Psi_{N}(q) = \sum_{ab \leq y_2} \frac{z_{a, b}}{\sqrt{ab}}  \sumplus_{\chi \mod{q}} L(1/2, \chi) \overline{\chi}(a) \chi(b).
	\end{align*}
	Using \cite[Lemma 3.3]{Bui} (noting that $V(x)=1+O(x^5)$ for $x\in(0,1),$ where $V(x)$ is from \cite[equation (9)]{Bui}, by shifting the integral to the left), we have, for any $\varepsilon > 0$,
	\begin{equation*}
		\Psi_{N}(q) = \phi^+(q)\sum_{\substack{kb = a \\ ab \leq y_2 \\(ab,q)=1}} \frac{z_{a, b}}{\sqrt{kab}} + O \left(Q ^{\theta/2+1/2+\varepsilon}\right) =\phi^+(q) \sum_{\substack{k b^2 \leq y_2 \\ (kb,q)=1}} \frac{z_{kb, b}}{kb} +  O\left(Q^{\theta/2+1/2+\varepsilon}\right).
	\end{equation*}
	Thus part (i) follows. For part (ii), we apply this with $z_{a, b} = 2 \cdot \mathbf{1}_{a=1}\mu(b)P[b]$ and recall~\eqref{eq:PsiNB=PsiN0}.
	
	\section{Proof of Propositions~\ref{prop:PsiMN_1IS} and \ref{prop:PsiM1,N1}}\label{sec:PsiM1N1}	
	In the rest of the paper, $\varepsilon > 0$ will be small but fixed. In particular, we assume that $\theta < 1/2-100\varepsilon$.
	We first give the proof of Proposition \ref{prop:PsiMN_1IS}(i) and Proposition \ref{prop:PsiM1,N1}, and then deduce Proposition \ref{prop:PsiMN_1IS}(ii) from the literature at the end of this section.
	
	 Let us first reduce Proposition~\ref{prop:PsiM1,N1} to the following variant, where the coefficients $z_{a,b}$ are supported on $a,b\in\N$ with $(a,b)=(ab,q)=1$. This simple observation will significantly simplify the following calculations.
	
	\begin{proposition} \label{prop:PsiM1N1gcd1}
		Let $\varepsilon_0 > 0$ and $\theta < 1/2$. Let $q \in \mathbb{N},$ let $y_1 \in [q^{2\varepsilon_0}, q^{\theta}]$ and let $y_2 \leq y_1 q^{-\varepsilon_0}$. Let
		\[
		M_1(\chi) := \sum_{b \leq y_1} \mu(b) \frac{\chi(b)}{\sqrt{b}} \left(1-\frac{\log b}{\log y_1}\right) \quad \text{and} \quad N(\chi) := \sum_{ab \leq y_2} \frac{z_{a, b} \overline{\chi}(a) \chi(b)}{\sqrt{ab}}
		\]
		with some complex coefficients $z_{a, b}$ such that, for every $\varepsilon > 0$,
			\begin{equation*}
				z_{a, b}\ll_\epsilon q^\epsilon \mathbf{1}_{(a, b) = 1} \mathbf{1}_{(ab, q) = 1}.
			\end{equation*}
		 Let $C_0 = \frac{\Gamma'(1/4)}{\Gamma(1/4)}-\log \pi$. Then the following hold.
		\begin{enumerate}[(i)]
			\item We have
			\begin{equation*}
				\begin{aligned}
					\Psi_{M_1,N}(q)&=\phi^+(q) \sum_{\substack{a\leq y_2 \\ (a, q) = 1}}\frac{\overline{z_{a,1}}}{a}\lz1+\frac{\log q}{\log y_1}-\frac{\log a}{\log y_1} + \frac{C_0}{\log y_1}\pz \\
					& \qquad +o\left(\sqrt{\Psi_{N, N} \varphi^+(q)}+\varphi^+(q)\right).
				\end{aligned}
			\end{equation*}
			\item  Assume that
			\begin{equation} \label{eq:zabexpbound}
				\sum_{ab \leq y_2} \frac{|z_{a, b}|^2}{ab} \ll \exp\left(c\frac{(\log q)^{3/5}}{(\log \log q)^{1/5}}\right)
			\end{equation}
			for a sufficiently small constant $c > 0$. Then
			\begin{equation*}
				\begin{aligned}
					\Psi_{M_1,N}(q)&=\phi^+(q) \sum_{\substack{a\leq y_2 \\ (a, q) = 1}}\frac{\overline{z_{a,1}}}{a}\lz1+\frac{\log q}{\log y_1}-\frac{\log a}{\log y_1} + \frac{C_0}{\log y_1}\pz +o\left(\varphi^+(q)\right).
				\end{aligned}
			\end{equation*}
			\item Assuming the quasi-Riemann hypothesis, part (ii) holds without the assumption~\eqref{eq:zabexpbound}.
		\end{enumerate}
	\end{proposition}
	Proposition~\ref{prop:PsiMN_1IS}(i) is immediate from Proposition~\ref{prop:PsiM1N1gcd1}(i) with $z_{a, b} = \mathbf{1}_{a=1} \mathbf{1}_{(b, q)=1} x_b$.
	\begin{proof}[Proof of Proposition~\ref{prop:PsiM1,N1} assuming Proposition~\ref{prop:PsiM1N1gcd1}]
		Let $z_{a, b}$ be as in Proposition~\ref{prop:PsiM1,N1}. The claim follows by applying Proposition~\ref{prop:PsiM1N1gcd1}(ii) with
		\[
		\mathbf{1}_{(a, b) = 1} \mathbf{1}_{(ab, q) = 1} \sum_{\substack{k^2 \leq y_2/ab \\ (k, q) = 1}} \frac{z_{ka, kb}}{k},
		\]
		in place of $z_{a, b}$.
	\end{proof}	
	
	Before turning to the proof of Proposition~\ref{prop:PsiM1N1gcd1} we introduce some notation and provide some auxiliary results. Inspired by Iwaniec and Sarnak~\cite[Section 6]{IS}, we define, for a sequence $x_{a, b}$ supported on $ab \leq y$ for some $y \geq 1$, the sequences $X_{u, v}$ and $X'_{u, v}$ via
	\[
	X_{u, v} = \sum_{ab \leq y/(uv)} \frac{x_{au, bv}}{abuv} \quad \text{and} \quad X'_{u, v} = \sum_{ab \leq y/(uv)} \frac{x_{au, bv}}{abuv} \log(ab). 
	\]
	Using that $\mathbf{1}_{n=1}(n) = (1 \ast \mu)(n)$ we see that
	\[
	\frac{x_{u, v}}{uv} = \sum_{abcd \leq y/(uv)} \mu(a)\mu(b) \frac{x_{acu, bdv}}{abcduv} = \sum_{ab \leq y/(uv)} \mu(a) \mu(b) X_{au, bv},
	\]
	and using that $\log n = (1 \ast \Lambda)(n)$, we see that
	\begin{equation} \label{eq:X'dec}
		X'_{u, v} =  \sum_{cdb \leq y/(uv)} \frac{x_{cd u, bv}}{bcduv} \Lambda(c) + \sum_{acd \leq y/(uv)} \frac{x_{au, cd v}}{acduv} \Lambda(c) = \sum_{c \leq y/(uv)} \Lambda(c)(X_{cu, v} + X_{u, cv}).
	\end{equation}
	Note that if $x_{a,b}$ is supported on $ab \leq y$ and $(a, b) = (ab,q)=1$, then so are $X_{a,b}$ and $X'_{a,b}$.
	
	Next we shall prove a lemma which takes care of the first steps of evaluating $\Psi_{N_1, N_2}$ for general Bui-type mollifiers $N_1(\chi)$ and $N_2(\chi)$. 
	
	\begin{lemma} \label{le:PsiN1N2FirstSteps}
		Let $q \in \mathbb{N}$, let $y_2 \leq y_1 \leq q^\theta$ with $\theta < 1/2$, and let
		\begin{equation*}
			N_1(\chi) := \sum_{ab \leq y_1} \frac{x_{a, b}\overline{\chi}(a)\chi(b)}{\sqrt{ab}} \quad \text{and} \quad N_2(\chi) := \sum_{ab \leq y_2} \frac{y_{a, b}\overline{\chi}(a)\chi(b)}{\sqrt{ab}}.
		\end{equation*}
		for some complex coefficients $x_{a, b}, y_{a, b}$ such that, for every $\varepsilon > 0$, 
		\begin{equation}\label{eq:assumption coprimality}
			x_{a, b} \ll q^\varepsilon \mathbf{1}_{(a, b) = 1} \mathbf{1}_{(ab, q)=1} \mathbf{1}_{ab \leq y_1} \, \text{ and } \, y_{a, b} \ll q^\varepsilon \mathbf{1}_{(a, b) = 1} \mathbf{1}_{(ab, q)=1} \mathbf{1}_{ab \leq y_2}.
		\end{equation}
		Then, for every $\varepsilon > 0$,
		\begin{align*}
			\Psi_{N_1, N_2}(q) &= \frac{\varphi(q)\varphi^+(q)}{q}\sum_{uv \leq y_2} \varphi(u) \varphi(v) (\log L^2(q) + 2\eta(uv)) X_{u, v} \overline{Y_{u, v}} \\
			& \qquad - \frac{\varphi(q)\varphi^+(q)}{q}\sum_{uv \leq y_2} \varphi(u) \varphi(v) \left(X'_{u, v} \overline{Y_{u, v}}+ X_{u, v} \overline{Y'_{u, v}}\right)+ O(q^{3/4+\theta/2+\varepsilon}),
		\end{align*}
		where $\eta(d)$ is as in \eqref{align:etaformulas} and $L(q)$ is defined by
		\begin{equation}\label{eq:definingL}
			\log L(q) = \frac{1}{2} \log \frac{q}{\pi} + \frac{\Gamma'(1/4)}{2\Gamma(1/4)} + \gamma + \eta(q).
		\end{equation}
	\end{lemma}
	
	\begin{proof}
		We generalize the calculations of Iwaniec and Sarnak~\cite[Sections 5, 6]{IS}. We have 
		\begin{align*}
			\Psi_{N_1, N_2}(q) &= \sum_{a_1 b_1 \leq y_1} \frac{x_{a_1, b_1}}{\sqrt{a_1 b_1}} \sum_{a_2 b_2 \leq y_2} \frac{\overline{y_{a_2, b_2}}}{\sqrt{a_2 b_2}} \sumplus_{\chi \mod{q}} |L(1/2, \chi)|^2 \chi(a_2 b_1) \overline{\chi}(a_1 b_2).
		\end{align*}
		We want to apply \cite[Lemma 3.2]{IS}, for which we need to ensure that $(a_2b_1,a_1 b_2)=1$. Thanks to \eqref{eq:assumption coprimality}, we can already assume that $(a_1, b_1) = (a_2, b_2) = 1$. 
		
		We write $(a_1, a_2) = r$ and $(b_1, b_2) = s$ and replace $a_1, a_2, b_1, b_2$ respectively by $ra_1,ra_2,bs_1,bs_2$, obtaining 
		\begin{align*}
			\Psi_{N_1, N_2}(q) &= \sum_{rs \leq y_2} \frac{1}{rs} \sum_{\substack{a_1 b_1 \leq y_1/(rs)}} \frac{x_{a_1r, b_1 s}}{\sqrt{a_1 b_1}} \sum_{\substack{a_2 b_2 \leq y_2/(rs) \\ (a_2, a_1) = 1 \\ (b_2, b_1) = 1}} \frac{\overline{y_{a_2r, b_2 s}}}{\sqrt{a_2 b_2}} \\
			& \quad \cdot \sumplus_{\chi \mod{q}} |L(1/2, \chi)|^2 \chi(a_2 b_1) \overline{\chi}(a_1 b_2).
		\end{align*}
		
		Now we have secured that $(a_2 b_1, a_1 b_2) = 1$ and can apply \cite[Lemma 3.2 and the formula for $\mathcal{B}_0(m_1, m_2)$ above (4.6)]{IS} to the sum over $\chi\mod q$, obtaining, for any $\varepsilon > 0$,
		\[
		\Psi_{N_1, N_2}(q) = \Psi_{N_1, N_2}^{\mathcal{M}}(q) + O\left(\Psi_{N_1, N_2}^{\mathcal{E}}(q) + (qy_1y_2)^{1/2+\epsilon}\right),
		\]
		where
		$$
		\begin{aligned}
			\Psi_{N_1, N_2}^{\mathcal{M}}(q) &= \frac{\varphi(q)\varphi^+(q)}{q} \sum_{rs \leq y_2} \frac{1}{rs} \sum_{\substack{a_1 b_1 \leq y_1/(rs)}} \frac{x_{a_1r, b_1 s}}{a_1 b_1} \sum_{\substack{a_2 b_2 \leq y_2/(rs) \\ (a_2, a_1) = 1 \\ (b_2, b_1) = 1}} \frac{\overline{y_{a_2r, b_2 s}}}{a_2 b_2} \log \frac{L^2(q)}{a_1 a_2 b_1 b_2}.
		\end{aligned}
		$$
		and
		$$
		\begin{aligned}
			\Psi_{N_1, N_2}^{\mathcal{E}}(q) &= \sum_{w \mid q} w \sum_{\substack{a_1 b_1 \leq y_1, a_2 b_2 \leq y_2 \\ \ell_1 \ell_2 \leq q^{1+\varepsilon} \\ a_2 b_1 \ell_1 \equiv a_1 b_2 \ell_2 \pmod{w} \\ a_2 b_1 \ell_1 \neq a_1 b_2 \ell_2}} \frac{1}{\sqrt{a_1 a_2 b_1 b_2 \ell_1 \ell_2}}\sum_{\substack{rs\leq y_2/(a_2 b_2) \\ rs \leq y_1/(a_1 b_1)}}\frac{|x_{a_1r,b_1s}y_{a_2r,b_2s}|}{rs}.
		\end{aligned}
		$$
		Let us first deal with the error term. By \eqref{eq:assumption coprimality}, the innermost sum is $\ll q^\varepsilon$ for any $\varepsilon > 0$. Next we write $k_1 = a_2 b_1 \ell_1$ and $k_2 = a_1 b_2 \ell_2$. Then $k_1 k_2 \leq q^{1 + \varepsilon} y_1 y_2$.  By symmetry we can assume that $k_1 > k_2$. Thus, for any $\varepsilon > 0$,
		\begin{align*}
			\Psi_{N_1, N_2}^{\mathcal{E}}(q) &\ll q^\varepsilon \cdot \sum_{w \mid q} w \sum_{k_2 \leq q^{1/2 +\varepsilon/2} \sqrt{y_1 y_2}} \frac{1}{\sqrt{k_2}} \sum_{\substack{k_2 < k_1 \leq  q^{1+\varepsilon}y_1 y_2/k_2 \\ k_1 \equiv \pm k_2 \pmod{w}}} \frac{1}{\sqrt{k_1}} \\
			&\leq q^\varepsilon \cdot \sum_{w \mid q} w \sum_{k_2 \leq q^{1/2+ \varepsilon/2} \sqrt{y_1 y_2}} \frac{1}{\sqrt{k_2}} \cdot \left(\frac{1}{\sqrt{w}} + \frac{q^{1/2+\varepsilon/2} \sqrt{y_1 y_2}}{w\sqrt{k_2}} \right) \\
			&\ll q^{3/4+2\varepsilon} (y_1 y_2)^{1/4} +  q^{1/2+2\varepsilon} \sqrt{y_1 y_2},
		\end{align*}
		which is acceptable.
		
		Now consider the main term $\Psi_{N_1, N_2}^{\mathcal{M}}(q)$.	We remove the conditions $(a_2, a_1) = 1$ and $(b_2, b_1) = 1$ using M\"obius inversion, introducing $\mu(c)$ and $\mu(d)$. We obtain
		$$
		\begin{aligned}
			\Psi_{N_1, N_2}^{\mathcal{M}}(q) &= \frac{\varphi(q)\varphi^+(q)}{q} \sum_{cdrs \leq y_2} \frac{\mu(c) \mu(d)}{c^2d^2rs}  \\
			& \cdot \sum_{\substack{a_1 b_1 \leq y_1/(cdrs)}} \frac{x_{a_1cr, b_1 ds}}{a_1 b_1} \sum_{\substack{a_2 b_2 \leq y_2/(cdrs)}} \frac{\overline{y_{a_2cr, b_2 ds}}}{a_2 b_2} \log \frac{L^2(q)}{a_1 a_2 b_1 b_2 c^2 d^2}.
		\end{aligned}
		$$ 
		Now we combine the variables $cr = u$ and $ds = v$, obtaining
		$$
		\begin{aligned}
			&\Psi_{N_1, N_2}^{\mathcal{M}}(q) = \frac{\varphi(q)\varphi^+(q)}{q} \sum_{uv \leq y_2} \frac{1}{u^2v^2} \sum_{\substack{a_1 b_1 \leq y_1/(uv)}} \frac{x_{a_1u, b_1 v}}{a_1 b_1} \sum_{\substack{a_2 b_2 \leq y_2/(uv)}} \frac{\overline{y_{a_2u, b_2 v}}}{a_2 b_2} \\
			& \cdot \left(\log \frac{L^2(q)}{a_1 a_2 b_1 b_2}\sum_{\substack{u = cr \\ v = ds}} \mu(c)r\mu(d)s - 2 \sum_{u = cr} \mu(c) r \log c \sum_{v =ds} \mu(d)s - 2 \sum_{u=cr} \mu(c)r \sum_{v = ds} \mu(s) d \log s\right).
		\end{aligned}
		$$
		Recalling~\eqref{align:etaformulas} and noting that $(u,v)=1$ by \eqref{eq:assumption coprimality} we see that,
		$$
		\begin{aligned}
			\Psi_{N_1, N_2}^{\mathcal{M}}(q) &= \frac{\varphi(q)\varphi^+(q)}{q} \sum_{uv \leq y_2} \varphi(u)\varphi(v) \sum_{\substack{a_1 b_1 \leq y_1/(uv)}} \frac{x_{a_1u, b_1 v}}{a_1 b_1 uv} \sum_{\substack{a_2 b_2 \leq y_2/(uv)}} \frac{\overline{y_{a_2u, b_2 v}}}{a_2 b_2 uv} \\
			& \cdot \left(\log \frac{L^2(q)}{a_1 a_2 b_1 b_2} + 2\eta(uv)\right).
		\end{aligned}
		$$
		and the claim follows from the definitions of $X_{u, v}, Y_{u, v}, X'_{u, v},$ and $Y'_{u, v}$.
	\end{proof}
	
	Before turning to the proof of Proposition~\ref{prop:PsiM1N1gcd1}, we use Lemma~\ref{le:PsiN1N2FirstSteps} to relate $\Psi_{N, N}$ to the sum of $|Z_{u, v}|^2$. 
	
	\begin{lemma} \label{le:PsiN0N0}
		Let $\theta < 1/2$ be fixed. Let $q \in \mathbb{N}$ and $y_2 \leq q^\theta$, and let
		\begin{equation*}
			N(\chi) := \sum_{ab \leq y_2} \frac{z_{a, b}\overline{\chi}(a)\chi(b)}{\sqrt{ab}} 
		\end{equation*}
		for some complex coefficients $z_{a, b}$ such that, for every $\varepsilon > 0$, 
		\[
		z_{a, b} \ll q^\varepsilon \mathbf{1}_{(a, b) = 1} \mathbf{1}_{(ab, q)=1} \mathbf{1}_{ab \leq y_2}.
		\] 
		Then, for any $\varepsilon > 0$,
		\[
		\Psi_{N, N}(q) = \left(1+O^\ast(2\theta) + o(1)\right) \frac{\varphi(q)\varphi^+(q)}{q} \log q \sum_{uv \leq y_2} \varphi(u) \varphi(v) |Z_{u, v}|^2+ O(q^{\theta+1/2+\varepsilon}),
		\]
		where $O^\ast(H)$ indicates a quantity which has absolute value at most $H$.
	\end{lemma}
	
	\begin{proof}
		By Lemma~\ref{le:PsiN1N2FirstSteps} and the fact that $\eta(n) = O(\log \log n)$, it suffices to prove that 
		\begin{equation} \label{eq:N0N0claim}
			S:= 2\left|\sum_{uv \leq y_2} \varphi(u) \varphi(v) Z_{u, v} \overline{Z'_{u, v}}\right| \leq (2\theta + o(1)) \log q \sum_{uv \leq y_2} \varphi(u) \varphi(v) |Z_{u, v}|^2.
		\end{equation}
		By~\eqref{eq:X'dec},
		\[
		S = 2\left|\sum_{uv \leq y_2} \sum_{c \leq y_2/(uv)} \Lambda(c)\varphi(u) \varphi(v) Z_{u, v} (\overline{Z_{cu, v}} + \overline{Z_{u, cv}}) \right|.
		\]
		Inspired by~\cite[(6.10)]{IS}, we use the inequalities
		\begin{align*}
			2|Z_{u, v}| |Z_{cu, v}| &\leq \frac{1}{\varphi(c)} |Z_{u, v}|^2 + \varphi(c) |Z_{cu, v}|^2, \\
			2|Z_{u, v}| |Z_{u, cv}| & \leq \frac{1}{\varphi(c)} |Z_{u, v}|^2 + \varphi(c) |Z_{u, cv}|^2.
		\end{align*}
		Hence, using also that $\phi(c)\phi(u)\leq \phi(cu)$,
		\begin{align*}
			S &\leq \sum_{uv \leq y_2} \varphi(u) \varphi(v) |Z_{u, v}|^2 \left(2 \sum_{c \leq y_2/(uv)} \frac{\Lambda(c)}{\varphi(c)} +  \sum_{u = rm_2} \Lambda(m_2) +  \sum_{v =rm_2} \Lambda(m_2)\right) \\
			&\leq\sum_{uv \leq y_2} \varphi(u) \varphi(v) |Z_{u, v}|^2 \cdot \left(\log \frac{y_2^2}{uv} + O(1)\right),
		\end{align*}
		and~\eqref{eq:N0N0claim} follows.
	\end{proof}
	
	\begin{proof}[Proof of Proposition~\ref{prop:PsiM1N1gcd1}]
		To avoid writing complex conjugates at each occurrence of the coefficients $z_{a, b}$, we evaluate $\Psi_{N,M_1}(q)$ and note that $\Psi_{M_1,N}(q)=\overline{\Psi_{N,M_1}(q)}$. 
		
		Write
		\[
		w_{a, b} = \mathbf{1}_{a=1}\mathbf{1}_{(b,q)=1} \mathbf{1}_{ab \leq y_1} \mu(b) \left(1-\frac{\log b}{\log y_1}\right)
		\]
		for the coefficients of $M_1(\chi)$. Then by Lemma \ref{lemma: Conrey} on the M\"obius sums, we see that, for $uv \leq y_2 \leq y_1 q^{-\varepsilon_0}$, and a sufficiently small $d > 0$ (depending on $\varepsilon_0$),
		\begin{align} \label{eq:Wformulas}
			\begin{aligned}
				W_{u, v} &= \frac{\mathbf{1}_{u = 1}\mathbf{1}_{(v,q)=1} \mu(v)}{v} \left(\frac{vq}{\phi(vq)\log y_1}+O\left(\exp\left(-d \frac{(\log q)^{3/5}}{(\log \log q)^{1/5}}\right)\right)\right), \\
				W'_{u, v} &= - \frac{\mathbf{1}_{u = 1}\mathbf{1}_{(v,q)=1} \mu(v)}{v} \Biggl(\frac{vq}{\phi(vq) \log y_1} \left(\log\frac{y_1}{v} - 2 \gamma - 2\eta(vq)\right) \\
				& \qquad \qquad \qquad \qquad \qquad \qquad + O\left(\exp\left(-d \frac{(\log q)^{3/5}}{(\log \log q)^{1/5}}\right)\right)\Biggr).
			\end{aligned}
		\end{align}
		Substituting these into Lemma~\ref{le:PsiN1N2FirstSteps}, we see that
		\begin{align}
		\label{eq:PsiNM1formula}
		\begin{aligned}
&			\Psi_{N, M_1}(q) = \varphi^+(q) \sum_{v \leq y_2} \mu(v) \frac{(\log L^2(q) + \log \frac{y_1}{v}-2\gamma -2\eta(q)) Z_{1, v} -Z'_{1, v}}{\log y_1} \\
			& \quad + O\left(q^{\theta+1/2+\epsilon}+\varphi^+(q) \exp\left(-d \frac{(\log q)^{3/5}}{(\log \log q)^{1/5}}\right) \sum_{v \leq y_2} \left(|Z_{1, v}| + |Z'_{1, v}|\right)\right).
		\end{aligned}
		\end{align}
		
		In the main term we have
		\[
		\sum_{v \leq y_2} \mu(v) Z_{1, v} = \sum_{abv \leq y_2} \frac{\mu(v) z_{a, bv}}{abv} = \sum_{a \leq y_2} \frac{z_{a, 1}}{a},
		\]
		and
		\[
		\sum_{v \leq y_2} \mu(v) (Z_{1, v} \log v + Z'_{1, v}) = \sum_{abv \leq y_2}\frac{ \mu(v) z_{a, bv} \log(abv)}{abv} = \sum_{a \leq y_2} \frac{z_{a, 1} \log a}{a},
		\]
		so, recalling the definition of $L(q)$ from~\eqref{eq:definingL} we obtain the claimed main term.
		
	To prove Proposition~\ref{prop:PsiM1N1gcd1}(ii), note that, by the definitions of $Z_{u, v}$ and $Z'_{u, v}$, and the Cauchy-Schwarz inequality,
				\begin{align*}
			&\sum_{v \leq y_2} (|Z_{1, v}|+|Z'_{1, v}|) \leq 2\log y_2 \sum_{abv \leq y_2} \frac{|z_{a, bv}|}{abv} \leq 2 \log y_2 \sum_{ab \leq y_2} \tau(b)  \frac{|z_{a, b}|}{ab} \\
			&\leq 2 \log y_2 \left(\sum_{ab \leq y_2} \frac{\tau(b)^2}{ab}\right)^{1/2} \left(\sum_{ab \leq y_2} \frac{|z_{a, b}|^2}{ab}\right)^{1/2} \ll (\log y_2)^{7/2} \left(\sum_{ab \leq y_2} \frac{|z_{a, b}|^2}{ab}\right)^{1/2}.
		\end{align*}
Under the assumption of Proposition~\ref{prop:PsiM1N1gcd1}(ii), with $c < d/4$, the error term in \eqref{eq:PsiNM1formula} is acceptable. Proposition~\ref{prop:PsiM1N1gcd1}(iii) follows similarly using Lemma~\ref{lemma: Conrey}(iii).
		
		To prove Proposition~\ref{prop:PsiM1N1gcd1}(i), note that, by~\eqref{eq:X'dec},
		\begin{align*}
			\sum_{v \leq y_2} (|Z_{1, v}|+|Z'_{1, v}|) &\leq \sum_{v \leq y_2} |Z_{1, v}| + \sum_{mv \leq y_2} \Lambda(m) (|Z_{m, v}| + |Z_{1, mv}|) \leq 2 \log y_2 \sum_{uv \leq y_2} |Z_{u, v}|.
		\end{align*}
		Now, by the Cauchy-Schwarz inequality and Lemma~\ref{le:PsiN0N0}, for any $\epsilon>0$,
		$$
\begin{aligned}
			\sum_{uv \leq y_2} |Z_{u, v}| &\leq \left(\sum_{uv \leq y_2} \varphi(u) \varphi(v) |Z_{u, v}|^2\right)^{1/2} \left(\sum_{uv \leq y_2}\frac{1}{\varphi(u) \varphi(v)}\right)^{1/2} \\
			& \ll \lz\sqrt{\frac{q\Psi_{N, N}}{\phi(q)\varphi^+(q) \log q}}+q^{\theta/2+1/4+\epsilon}\pz \log q,
\end{aligned}
	$$
		and thus Proposition~\ref{prop:PsiM1N1gcd1}(i) also follows.
	\end{proof}
	\begin{remark} \label{rem:increase}
		The trick of taking a slightly longer Michel-Vanderkam mollifier was crucial for obtaining the good error terms in~\eqref{eq:Wformulas}. If we had $y_1 = y_2$, we would need to work with very small values of $y_1/v$ which would lead to worse error terms (as the error obtained from the proof of Lemma~\ref{lemma: Conrey} has a term of shape $(y/j)^{-\sigma}$) and consequently stricter requirements on the coefficients of $\Psi_N$.
		
		With a more careful treatment of the error term in the proof of Proposition~\ref{prop:PsiM1N1gcd1}(i), we could prove Proposition~\ref{prop:PsiM1N1gcd1}(i) with the condition $y_2 \leq y_1 q^{-\varepsilon}$ relaxed to $y_2 \leq y_1$. However, such a result is not directly applicable to prove Theorem~\ref{thm:main general result DirL} as the error term $o(\sqrt{\Psi_{N,N} \varphi^+(q)})$ might turn out to be troublesome: In the proof of Theorem~\ref{thm:main general result DirL}, Proposition~\ref{prop:PsiM1N1gcd1} is applied for $N(\chi) = N_0(\chi)$ defined in~\eqref{eq:N0def}, and one might in principle have that $\Psi_{N_0, N_0}$ is large when $\Psi_{N_\mathrm{B}, N_\mathrm{B}}$ is small. Writing
		\begin{equation*}
			N_{\mathrm{B}}(\chi)=\sum_{ab \leq Q^\theta} \frac{x_{a, b}\overline{\chi}(a)\chi(b)}{\sqrt{ab}} + \overline\varepsilon_\chi \sum_{ab \leq Q^\theta} \frac{y_{a, b}\chi(a)\overline{\chi}(b)}{\sqrt{ab}} =: N_1(\chi) + N_2(\chi),
		\end{equation*}
		say, and using the above-mentioned variant of Proposition~\ref{prop:PsiM1N1gcd1}(i) in place of Proposition~\ref{prop:PsiM1N1gcd1}(ii), the only possibility for Theorem~\ref{thm:main general result DirL} not to hold would be if there existed $N_{\mathrm{B}}(\chi)$ such that
		\[
		\Psi_{N_\mathrm{B}, N_\mathrm{B}} = \Psi_{N_1, N_1} + \Psi_{N_2, N_2} + 2\re \Psi_{N_1, N_2} \leq (1 + 1/\theta) \Psi_{N_{\mathrm{B}}}^2
		\]
		but $\Psi_{N_1, N_1}$ and $\Psi_{N_2, N_2}$ are unbounded. This scenario seems unlikely but it is not clear how to rule it out.  
	\end{remark}

We have now concluded the proofs of Proposition \ref{prop:PsiMN_1IS}(i) and Proposition \ref{prop:PsiM1,N1}, and it remains to prove Proposition \ref{prop:PsiMN_1IS}(ii).
\begin{proof}[Proof of Proposition \ref{prop:PsiMN_1IS}(ii)]
We could deduce Proposition \ref{prop:PsiMN_1IS}(ii) form Lemma~\ref{le:PsiN1N2FirstSteps} and Conrey's lemma~\cite[Lemma 10]{Conrey}, but for simplicity we refer to Michel and Vanderkam \cite{M-V}. Our $\Psi_{M,M}$ is equal to $\mathscr{Q}_1(P_k)/(\hat q \Gamma\bfrac{1}{4}^2)$ from \cite{M-V} with $k=0$, $P_k(x)=x$ and $\Delta=\frac{2\log y_1}{\log q}+O\bfrac{1}{\log q}$. The proof now follows from the evaluation of $\mathscr Q_1(P_k)$ in \cite[eq. (16)]{M-V}.
\end{proof}
	
	\section{Proof of Proposition \ref{prop:PsiM_2N_1}}\label{sec:PsiM2N1}
	We adapt the computations of Pratt (that concerned the special case $z_{a, b} = \frac{\Lambda(a)}{\log q} \mu(b) P_2(\frac{\log y_2/(ab)}{\log y_2})$) from \cite[Proof of Lemma 3.4]{Pratt} with some simplifications. Although this does not affect our computations, let us point out that in the published version of~\cite{Pratt}, the integrand on the right-hand side of the claim in Lemma 3.4 should be $P_2(1-x) P_3(1-x\frac{\theta_2}{\theta_3})$ instead of $P_2(x)P_3(x)$. 
	
	Similarly to the previous section, we notice that $\Psi_{M_2, N} = \overline{\Psi_{N, M_2}}$ and compute $\Psi_{N,M_2}$ (although now that $z_{a, b} \in \mathbb{R}$ we would avoid $\overline{z_{a, b}}$ also without doing this). Writing $w_b := \mu(b) (1-\frac{\log b}{\log y_1})$, we have
	$$
	\begin{aligned}
		\Psi_{N, M_2} &=\frac{1}{Q^2}\sum_{q\geq 1}\Phi\bfrac{q}{Q}\frac{q}{\phi(q)} \sum_{b_2 \leq y_1} \frac{w_{b_2}}{\sqrt{b_2}} \sum_{a_0 b_0 \leq y_2} \frac{z_{a_0, b_0}}{\sqrt{a_0 b_0}} \\
		&\cdot \sumplus_{\chi \mod{q}} |L(1/2, \chi)|^2 \varepsilon_\chi \chi(b_0 b_2) \overline{\chi}(a_0).
	\end{aligned}
	$$
	Using \cite[Lemma 4.1]{Pratt} (which is a slight variant of~\cite[(2.5)]{IS} with a similar proof),
	\begin{align}
		\label{eq:Psi23FE}
		\begin{aligned}
			\Psi_{N, M_2} &= \frac{2}{Q^2}\sum_{q\geq 1}\Phi\bfrac{q}{Q}\frac{q}{\phi(q)}\sum_{\substack{b_2 \leq y_1}} \frac{w_{b_2}}{\sqrt{b_2}} \sum_{a_0 b_0 \leq y_2} \frac{z_{a_0, b_0}}{\sqrt{a_0 b_0}} \cdot \sum_{\ell_1 \geq 1} \frac{1}{\sqrt{\ell_1}} \\
			&\quad\cdot \sumplus_{\chi \mod{q}} \varepsilon_\chi \chi(\ell_1 b_0 b_2) \overline{\chi}(a_0) \sum_{\ell_2 \geq 1} \frac{\overline{\chi}(\ell_2) V_2(\ell_1 \ell_2/q)}{\sqrt{\ell_2}},
		\end{aligned}
	\end{align}
	where
	\begin{align*}
		V_2(x) := \frac{1}{2\pi i} \int_{(1)} \frac{\Gamma^2\left(\frac{s}{2}+\frac{1}{4}\right)}{\Gamma^2\left(\frac{1}{4}\right)} \cdot \frac{G(s)}{s} \pi^{-s} x^{-s} ds  
	\end{align*}
	with $G(s)$ an even polynomial satisfying $G(0) = 1$ and vanishing to the second order at $s = 1/2$.
	
	As in \cite{Pratt}, the main contribution will come from the terms where \begin{align}
		\label{eq:divcond}
		\ell_1 b_0 b_2 \mid a_0.
	\end{align} We therefore write
	$$
	\Psi_{N, M_2}=\Psi_{N, M_2}^\M+\Psi_{N, M_2}^\E,
	$$ where $\Psi_{N, M_2}^\M$ is the part of \eqref{eq:Psi23FE} with the extra condition \eqref{eq:divcond}, while $\Psi_{N, M_2}^\E$ is the rest.
	
	We do not need the average over $q$ to evaluate the main terms, so we also write
	$$
	\Psi^{\mathcal M}_{N, M_2}=\frac{1}{Q^2}\sum_{q\geq 1}\Phi\bfrac {q}{Q}\frac{q}{\phi(q)}\Psi^{\mathcal M}_{N,M_2}(q),
	$$ with
	\begin{align} \label{eq:N1M2MT}
		\begin{aligned}
			\Psi^{\mathcal M}_{N, M_2}(q)&=2 \sum_{\substack{b_2 \leq y_1}} \frac{w_{b_2}}{\sqrt{b_2}} \sum_{\substack{a_0 b_0 \leq y_2 \\ b_0 b_2 \mid a_0}} \frac{z_{a_0, b_0}}{\sqrt{a_0 b_0}} \cdot \sum_{\substack{\ell_1 \geq 1\\ \ell_1b_0 b_2|a_0}} \frac{1}{\sqrt{\ell_1}} \\
			& \cdot \sumplus_{\chi \mod{q}} \varepsilon_\chi \chi(\ell_1 b_0 b_2) \overline{\chi}(a_0) \sum_{\ell_2 \geq 1} \frac{\overline{\chi}(\ell_2) V_2(\ell_1 \ell_2/q)}{\sqrt{\ell_2}}.
		\end{aligned}
	\end{align}
	
	\subsection{Main terms}
	Now we evaluate the main terms $\Psi_{N, M_2}^\M(q)$.
	
	From \cite[Proposition 4.3]{Pratt}, we obtain 
	\begin{align} \label{eq:revAFE}
		\sum_{\ell_2 \geq 1} \frac{\overline{\chi}(\ell_2) V_2(\ell_1 \ell_2/q)}{\sqrt{\ell_2}} = L(1/2, \overline{\chi})-\overline{\varepsilon_\chi} \sum_{n \geq 1} \frac{\chi(n)}{\sqrt{n}} F(n/\ell_1),
	\end{align}
	where
	\begin{align*}
		F(x) := \frac{1}{2\pi i} \int_{(1)} \frac{\Gamma\left(\frac{s}{2}+\frac{1}{4}\right) \Gamma\left(-\frac{s}{2}+\frac{1}{4}\right)}{\Gamma^2(\frac{1}{4})} \cdot \frac{G(s)}{s} x^{-s} ds.
	\end{align*}
	For $x\ll 1$, by shifting the integral to $\re(s)=-1/2$ (note that the poles at $s=\pm1/2$ are cancelled by the vanishing condition on $G(1/2)$), we have 
	$$
	F(x)=1+O(x^{1/2}).
	$$ 
	We choose the polynomial $G(s)$ so that it vanishes at all the poles of $\Gamma(\frac{\pm s}{2} + \frac{1}{4})$ in the region $|s| \leq A$ with $A>100$ large but fixed. Then, for $x\gg1$, we have $F(x)\ll x^{-100}$ say, by shifting the integral to the right. We also find that $F(x)=1-F(1/x)$ by shifting the integral to $\re(s)=-1$ and changing  variables $s\mapsto-s$.
	
	By \cite[Lemma 6.1]{Pratt} we obtain from~\eqref{eq:revAFE} that
	\begin{align}
		\label{eq:Vl1l2dec}
		\sum_{\ell_2 \geq 1} \frac{\overline{\chi}(\ell_2) V_2(\ell_1 \ell_2/q)}{\sqrt{\ell_2}} = \sum_{n \geq 1} \frac{\overline{\chi}(n)}{\sqrt{n}} V_1(n/\sqrt{q}) + \overline{\varepsilon_\chi}\sum_{n \geq 1} \frac{\chi(n)}{\sqrt{n}} \left(V_1(n/\sqrt{q}) - F(n/\ell_1)\right), 
	\end{align}
	where
	\begin{align*}
		V_1(x) := \frac{1}{2\pi i} \int_{(1)} \frac{\Gamma\left(\frac{s}{2}+\frac{1}{4}\right)}{\Gamma(\frac{1}{4})} \cdot \frac{G_1(s)}{s} \pi^{-s/2} x^{-s} ds,
	\end{align*}
	with $G_1(s)$ an even polynomial satisfying $G_1(0) = 1$. Here, by shifting the contour to $\re(s) = -1/3$ for $x \leq 1$ and to $\re(s) = 150$ for $x > 1$, we see that
	\begin{equation} \label{eq:V1est}
		V_1(x) =\begin{cases}
			1+O(x^{1/3}),&\hbox{for $x\leq 1$;}\\
			O(x^{-100}), &\hbox{for $x> 1$.}
		\end{cases}
	\end{equation}
	
	Substituting~\eqref{eq:Vl1l2dec} into~\eqref{eq:N1M2MT} gives
	$$
	\begin{aligned}
		\Psi_{N, M_2}^\M(q)&= 2\sum_{\substack{b_2 \leq y_1 \\ \ell_1 \geq 1}} \sum_{\substack{a_0 b_0\leq y_2 \\ \ell_1 b_0 b_2 \mid a_0}} \frac{w_{b_2} z_{a_0, b_0}}{\sqrt{b_2a_0b_0\ell_1}}
		\cdot  \sumplus_{\chi \mod{q}} \varepsilon_\chi \chi(\ell_1 b_0 b_2) \overline{\chi}(a_0) \\
		& \qquad \cdot\lz\sum_{n\geq 1} \frac{\overline{\chi}(n)}{\sqrt{n}} V_1(n/\sqrt{q}) + \overline{\varepsilon_\chi}\sum_{n\geq 1} \frac{\chi(n)}{\sqrt{n}} \left(V_1(n/\sqrt{q}) - F(n/\ell_1)\right)\pz\\
		&=: K_1+K_2,
	\end{aligned}
	$$
	say. Rearranging, 
	\begin{align*}
		K_1 &= 2\sum_{\substack{b_2 \leq y_1 \\ \ell_1 \geq 1}} \sum_{\substack{a_0 b_0\leq y_2 \\ \ell_1 b_0 b_2 \mid a_0}} \frac{w_{b_2} z_{a_0, b_0}}{\sqrt{b_2a_0b_0\ell_1}} \sum_{n\geq 1} \frac{V_1(n/\sqrt{q})}{\sqrt{n}}  \sumplus_{\chi \mod{q}} \varepsilon_\chi \chi(b_0 b_2 \ell_1) \overline{\chi}(a_0n).
	\end{align*}
	By \eqref{eq:orthogwitheps}, the sum over $\chi$ is $\ll q^{1/2+\epsilon/4}$. Using that $\ell_1 b_0 b_2 \mid a_0$ and $z_{a_0, b_0} \ll q^{\varepsilon/4}$, we obtain
	\begin{align*}
		K_1 \ll  q^{1/2+\epsilon/2} \sum_{a_0 \leq y_2} \frac{\tau_4(a_0)}{\sqrt{a_0}} \sum_{n \geq 1} \frac{V_1(n/\sqrt q)}{\sqrt{n}} \ll Q^{\theta/2+3/4+\epsilon},
	\end{align*}
	which is $O(Q^{1-\varepsilon})$ for $\theta<1/2-4\varepsilon$.
	
	Hence
	$$
	\Psi_{N, M_2}^\M(q)=K_2+O(Q^{ 1-\epsilon}).
	$$
	By orthogonality of even characters \eqref{eq:orthogonality}, we have
	\begin{align*}
		K_2 &= 2\sum_{\substack{b_2 \leq y_1 \\ \ell_1 \geq 1}} \sum_{\substack{a_0 b_0\leq y_2\\ \ell_1 b_0 b_2 \mid a_0}} \frac{w_{b_2} z_{a_0, b_0}}{\sqrt{b_2a_0b_0\ell_1}} \sum_{n\geq 1} \frac{V_1(n/\sqrt{q}) - F(n/\ell_1)}{\sqrt{n}} \\
		& \qquad \cdot \sumplus_{\chi \mod{q}} \chi(b_0 b_2 \ell_1 n) \overline{\chi}(a_0) \\
		&= \sum_{q = wv} \mu(v) \varphi(w) \sum_{\substack{b_2 \leq y_1 \\ \ell_1 \geq 1 \\ (b_2 \ell_1, q) = 1}} \sum_{\substack{a_0 b_0\leq y_2\\ \ell_1 b_0 b_2 \mid a_0 \\ (a_0 b_0, q) = 1}} \frac{w_{b_2} z_{a_0, b_0}}{\sqrt{b_2a_0b_0\ell_1}} \\
		& \quad \cdot \sum_{\substack{(n, q) = 1 \\ b_2 b_0 \ell_1 n \pm a_0 \equiv 0 \mod{w}}} \frac{V_1(n/\sqrt{q}) - F(n/\ell_1)}{\sqrt{n}}.
	\end{align*}
	
	Since $\ell_1 b_0 b_2 \mid a_0$, the congruence condition on the last line becomes
	\begin{equation}
		\label{eq:nmodw}
		n \equiv \pm \frac{a_0}{b_0 b_2 \ell_1} \mod{w}.
	\end{equation}
	Consider first the contribution of $w \leq q^{3/4}$. They contribute
	\begin{align*}
		\ll Q^\epsilon\sum_{\substack{w \mid q \\ w \leq q^{3/4}}} w \sum_{a_0 \leq y_2} \frac{\tau_4(a_0)}{\sqrt{a_0}} \left(\frac{\sqrt{Q^{1/2+\varepsilon}}}{w} + 1\right) \ll Q^{1/4+\theta/2+2\varepsilon}+ Q^{3/4+\theta/2+2\varepsilon} \ll Q^{ 1-\varepsilon}
	\end{align*}
	for $\theta < 1/2-6\varepsilon$.
	
	Hence we can concentrate on $w > q^{3/4}$. Note that when studying $K_2$ we can essentially restrict to the case where $n \leq q^{1/2+\varepsilon}$ and $\ell_1 b_0 b_2 \leq a_0 \leq Q^\theta \leq Q^{1/2}$. Hence~\eqref{eq:nmodw} actually implies $n = \frac{a_0}{b_0 b_2 \ell_1}$. This yields
	\begin{align*}
		K_2 &= \sum_{q = wv} \mu(v) \varphi(w) \sum_{\substack{b_2 \leq y_1 \\ (b_2, q) = 1}} \frac{w_{b_2}}{\sqrt{b_2}} \sum_{\substack{a_0 b_0 \leq y_2 \\ (a_0b_0,q)=1}} \frac{z_{a_0, b_0}}{\sqrt{a_0 b_0}} \cdot \sum_{\substack{\ell_1 \geq 1 \\ (\ell_1, q) = 1}} \frac{1}{\sqrt{\ell_1}} \\
		& \quad \sum_{\substack{n \geq 1 \\ (n, q) = 1 \\ b_0 b_2 \ell_1 n = a_0 }} \frac{V_1(n/\sqrt{q}) - F(n/\ell_1)}{\sqrt{n}} + O(Q^{ 1-\varepsilon}).
	\end{align*}
	Substituting $a_0 = b_0 b_2 \ell_1 n,$ using~\eqref{eq:evencharscount}, and recalling that $y_2 \leq y_1$, we see that
	\begin{align*}
		K_2 &=  2 \phi^+(q) \sum_{\substack{b_2 n \ell_1 b_0^2 \leq y_2 \\ (b_0 b_2\ell_1 n, q) = 1}} w_{b_2} \frac{z_{b_0 b_2 \ell_1 n, b_0}(V_1(n/\sqrt{q}) - F(n/\ell_1))}{b_0 b_2 \ell_1 n} + O(Q^{1-\varepsilon}).
	\end{align*}
	Using from~\eqref{eq:V1est} that $V_1(n/\sqrt{q}) = 1+O(q^{(\theta-1/2)/3})$ and the symmetry in $n, \ell_1$ together with  $F(x) = 1-F(x^{-1})$, we obtain, for $\theta < 1/2-7\varepsilon$,
	\begin{align*}
		K_2 = \phi^+(q) \sum_{\substack{b_2 n \ell_1 b_0^2 \leq y_2, \\ (b_0 b_2 \ell_1 n,q)=1}} w_{b_2} \frac{z_{b_0 b_2 \ell_1 n, b_0}}{b_0 b_2 \ell_1 n} + O(Q^{1-\varepsilon}).
	\end{align*}
	Hence, recalling that $w_b = \mu(b)(1-\frac{\log b}{\log y_1})$, we obtain
	\begin{align}
		\label{eq:PsiN0M2Mformula}
		\Psi_{N, M_2}^\M (q) &=\phi^+(q) \sum_{\substack{b_2 n \ell_1 b_0^2 \leq y_2 \\ (b_0 b_2\ell_1 n,q)=1}}  \mu(b_2) \left(1-\frac{\log b_2}{\log y_1}\right) \frac{z_{b_0 b_2 \ell_1 n, b_0}}{b_0 b_2 \ell_1 n} + O(Q^{1-\varepsilon}).
	\end{align}
	Writing $u = b_2 \ell_1$ and $r = b_2 n\ell_1$, the sum above is
	\begin{align*}
		\phi^+(q)\sum_{\substack{n u b_0^2 \leq y_2,\\(nub_0,q)=1}}  \frac{z_{b_0 n u, b_0}}{u b_0 n}\sum_{u = b_2 \ell_1} \mu(b_2) - \frac{\phi^+(q)}{\log y_1}\sum_{\substack{r b_0^2 \leq y_2,\\(rb_0,q)=1}}  \frac{z_{b_0 r, b_0}}{b_0 r} \sum_{r = b_2 n \ell_1} \mu(b_2) \log b_2 .
	\end{align*}
	Now the first term vanishes unless $u = 1$. The Dirichlet series corresponding to the convolution in the sum over $r=b_2 n \ell_1$ is
	\begin{align*}
		-\left(\frac{1}{\zeta}\right)' \cdot \zeta^2 = \zeta', \quad \text{and thus} \quad \sum_{r = b_2 n \ell_1} \mu(b_2) \log b_2 = -\log r,
	\end{align*}
	which yields the final result
	\begin{align*}
		\Psi_{N, M_2}^\M(q) &=\phi^+(q)\sum_{\substack{n b_0^2 \leq y_2,\\ (nb_0,q)=1}}  \frac{z_{b_0 n, b_0}}{b_0 n} + \phi^+(q)\sum_{\substack{r b_0^2 \leq y_2, \\ (rb_0,q)=1}} \frac{z_{b_0 r, b_0}}{b_0 r} \frac{\log r}{\log y_1} + O(Q^{ 1-\varepsilon})\\
		&= \phi^+(q)\sum_{\substack{n b_0^2 \leq y_2,\\ (nb_0,q)=1}}  \frac{z_{b_0 n, b_0}}{b_0 n}\left(1+\frac{\log n}{\log y_1}\right) + O(Q^{ 1-\varepsilon})
	\end{align*}
	corresponding to the main term in Proposition~\ref{prop:PsiM_2N_1}.
	
	\subsection{Error terms} \label{se:PropM2N1Errors}
	
	In this section we estimate the error term that arises from the terms with
	\begin{align*}
		\ell_1 b_0 b_2 \nmid a_0.
	\end{align*}
	Recall that 
	$$
	\begin{aligned}
		\Psi_{N, M_2}^\E&=\frac{2}{Q^2}\sum_{q\geq 1}\Phi\bfrac qQ \frac{q}{\phi(q)}\sum_{\substack{b_2 \leq y_1}} \frac{w_{b_2}}{\sqrt{b_2}} \sum_{a_0 b_0 \leq y_2} \frac{z_{a_0, b_0}}{\sqrt{a_0 b_0}}\\
		& \cdot \sum_{\substack{\ell_1,\ell_2\geq 1 \\ \ell_1 b_0 b_2 \nmid a_0}} \frac{V_2(\ell_1\ell_2/q)}{\sqrt{\ell_1\ell_2}}  \sumplus_{\chi \mod{q}} \varepsilon_\chi \chi(\ell_1 b_0 b_2) \overline{\chi}(a_0\ell_2).
	\end{aligned}
	$$
	
	We start by applying~(\ref{eq:orthogwitheps}). Recalling $z_{a, b} \in \mathbb{R}$, this yields
	\begin{align*}
		\begin{aligned}
			\Psi^{\E}_{N, M_2} &= \frac{2}{Q^2} \, \re \Biggl( \sum_{\substack{v, w \geq 1 \\ (v, w)=1}} \Phi\left(\frac{vw}{Q}\right) \frac{\mu^2(v) (vw)^{1/2} \varphi(w)}{\varphi(vw)}  \sum_{\substack{b_2 \leq y_1 \\ (b_2, vw) = 1}} \frac{w_{b_2}}{\sqrt{b_2}}  \\
			& \qquad \cdot \sum_{\substack{a_0 b_0 \leq y_2 \\ (a_0 b_0, vw) = 1}} \frac{z_{a_0, b_0}}{\sqrt{a_0 b_0}} \sum_{\substack{\ell_1, \ell_2 \geq 1 \\ (\ell_1 \ell_2, vw) = 1 \\ \ell_1 b_0 b_2 \nmid a_0}} \frac{V_2(\ell_1 \ell_2/(vw))}{\sqrt{\ell_1 \ell_2}} e\left(\frac{a_0 \ell_2 \overline{\ell_1 b_0 b_2 v}}{w}\right)\Biggr).
		\end{aligned}
	\end{align*}   
	
	Now we use dyadic splitting and a smooth partition of unity to localize our variables into dyadic (or $4$-adic) intervals and use inverse Mellin transform to separate variables in smooth functions. We denote
	$$\begin{aligned}
		\mathcal{E} &:=  \frac{W^{1/2}}{\left(A_0 B_0 B_2 L_1 L_2 V \right)^{1/2}} \Biggl| \sum_{\substack{v \sim V, w \in \mathbb{N} \\ (v, w)=1}} \alpha(v) G_3\left(\frac{w}{W}\right) \sum_{\substack{b_2 \sim B_2 \\ (b_2, vw) = 1}} \beta(b_2) \\
		&\cdot \sum_{\substack{a_0 \sim A_0, b_0 \sim B_0 \\ (a_0 b_0, vw) = 1}} z'_{a_0, b_0} \sum_{\substack{\ell_1, \ell_2 \geq 1 \\ (\ell_1 \ell_2, vw) = 1 \\ \ell_1 b_0 b_2 \nmid a_0}} G_1\left(\frac{\ell_1}{L_1}\right) G_2\left(\frac{\ell_2}{L_2}\right) e\left(\frac{a_0 \ell_2 \overline{\ell_1 b_0 b_2 v}}{w}\right)\Biggr|,
	\end{aligned}$$ where $z'_{a_0, b_0} := \sqrt{\frac{A_0 B_0}{a_0 b_0}} z_{a_0, b_0}$,
	$V,W,L_1,L_2,A_0,B_0, B_2 \geq 1/2$ satisfy 
	\begin{align} \label{eq:VariableRanges}
		VW \asymp Q, \quad L_1 L_2 \leq Q^{1+\varepsilon/10}, \quad A_0 B_0, B_2 \leq Q^\theta,
	\end{align} $G_1,G_2,G_3$ are smooth functions supported in $[1/2,2],$ $\alpha(v),\beta(b_2)$ are arbitrary bounded coefficients, and the notation $a \sim A$ means that $a \in (A, 2A]$. The following lemma shows that it is enough to prove that $\mathcal E=O(Q^{2-\epsilon/10}).$ 
	\begin{lemma}\label{le:dyadic decomposition}
		For some $V,W,L_1,L_2,A_0,B_0,B_2,\alpha(v),\beta(b_2),G_1,G_2,G_3$ as above, we have
		$$\Psi^{\mathcal E}_{N,M_2}\ll Q^{-2+\epsilon/100} \mathcal E.$$
	\end{lemma}
	\begin{proof}
		Let $G$ be a smooth function supported in $[1/2,2]$ that satisfies 
		$$
		\sum_{n\geq 0}G\bfrac{x}{2^n} = 1 \quad \text{for all $x \geq 1$}
		$$ 
		(see for example \cite[Lemme 2]{Fouvry} for a construction of such a function).
		Utilizing this for the variables $\ell_1, \ell_2$ and $w$ and making a dyadic splitting for the variables $a_0, b_0, b_2, v$, we obtain \begin{equation}\label{eq:partition}
			\Psi^{\mathcal E}_{N,M_2}\leq \frac{2}{Q^2}\sumd_{V,W,L_1,L_2,A_0,B_0, B_2}|\Psi(V,W,L_1,L_2,A_0,B_0, B_2)|,
		\end{equation} where $\sumd$ means that the variables run over powers of two, and 
		$$
		\begin{aligned}
			&\Psi(V,W,L_1,L_2,A_0,B_0, B_2)\\
			&=\sum_{\substack{v \sim V \\ w \geq 1 \\ (v, w)=1}} \Phi\left(\frac{vw}{Q}\right) \frac{\mu^2(v) (vw)^{1/2}}{\varphi(v)}\sum_{\substack{b_2 \leq y_1 \\ b_2 \sim B_2 \\ (b_2, vw) = 1}} \frac{w_{b_2}}{\sqrt{b_2}} \sum_{\substack{a_0 b_0 \leq y_2 \\ a_0 \sim A_0, b_0 \sim B_0 \\ (a_0 b_0, vw) = 1}} \frac{z_{a_0, b_0}}{\sqrt{a_0 b_0}} \\
			& \cdot \sum_{\substack{\ell_1, \ell_2 \geq 1 \\ (\ell_1 \ell_2, vw) = 1 \\ \ell_1 b_0 b_2 \nmid a_0}} \frac{V_2(\ell_1 \ell_2/(vw))}{\sqrt{\ell_1 \ell_2}} e\left(\frac{a_0 \ell_2 \overline{\ell_1 b_0 b_2 v}}{w}\right) G\bfrac{w}{W}G\bfrac{\ell_1}{L_1}G\bfrac{\ell_2}{L_2}.
		\end{aligned}
		$$
		Due to the support of $\Phi$, the conditions $a_0b_0 \leq y_2, \, b_2\leq y_1$ and the fast decay of $V_2(x)$, up to an acceptable error term we can restrict the sum in \eqref{eq:partition} to the ranges in \eqref{eq:VariableRanges}. Recalling that the sum in \eqref{eq:partition} runs over powers of 2, we thus have
		$$
		\Psi^{\mathcal E}_{N,M_2}\ll \frac{(\log Q)^{7}}{Q^2} \max_{V,W,L_1,L_2,A_0,B_0, B_2} \lab\Psi(V,W,L_1,L_2,A_0,B_0, B_2)\rab,
		$$ where the maximum runs over the variables satisfying \eqref{eq:VariableRanges}.

		Now we use Mellin inversion to separate variables in the smooth functions $\Phi$ and $V_2$. Writing $\Phi, V_2$ in terms of their Mellin transforms $\widetilde\Phi,\widetilde V_2$, we obtain
		$$
		\begin{aligned}
			&\Psi(V,W,L_1,L_2,A_0,B_0, B_2) \\
			&=\bfrac{1}{2\pi i}^2\int_{(c_1)} \int_{(c_2)} \Psi(V,W,L_1,L_2,A_0,B_0, B_2;s_1,s_2) ds_2ds_1,
		\end{aligned}
		$$ where $s_j=c_j+it_j$, and 
		$$
		\begin{aligned}
			&\Psi(V,W,L_1,L_2,A_0,B_0, B_2;s_1,s_2)\\	
			&:=\widetilde\Phi(s_1) \widetilde V_2(s_2) Q^{s_1} \sum_{\substack{v \sim V \\ w \geq 1 \\ (v, w)=1}} \frac{\mu^2(v) (vw)^{1/2+s_2-s_1}}{\varphi(v)}\sum_{\substack{b_2 \leq y_1 \\ b_2 \sim B_2 \\ (b_2, vw) = 1}} \frac{w_{b_2}}{\sqrt{b_2}} \sum_{\substack{a_0 b_0 \leq y_2 \\ a_0 \sim A_0, b_0 \sim B_0 \\ (a_0 b_0, vw) = 1}} \frac{z_{a_0, b_0}}{\sqrt{a_0 b_0}} \\
			&\cdot \sum_{\substack{\ell_1, \ell_2 \geq 1 \\ (\ell_1 \ell_2, vw) = 1 \\ \ell_1 b_0 b_2 \nmid a_0}} \frac{(\ell_1 \ell_2)^{-s_2}}{\sqrt{\ell_1 \ell_2}} e\left(\frac{a_0 \ell_2 \overline{\ell_1 b_0 b_2 v}}{w}\right) G\bfrac{w}{W}G\bfrac{\ell_1}{L_1}G\bfrac{\ell_2}{L_2}.
		\end{aligned}
		$$
		Thanks to the fast decay of the Mellin transforms, we can limit the integrals to $|\im(s_1)|,|\im(s_2)|\leq Q^{\epsilon/1000}$. Setting also $c_1=c_2=1/\log Q$, we obtain
		\begin{align*}
			&|\Psi(V,W,L_1,L_2,A_0,B_0, B_2)| \\
			&\ll Q^{\epsilon/500}\max_{|t_1|,|t_2|\leq Q^{\frac{\epsilon}{1000}}} |\Psi(V,W,L_1,L_2,A_0,B_0, B_2;\tfrac{1}{\log Q}+it_1,\tfrac{1}{\log Q}+it_2)|.
		\end{align*}
		This is $\ll Q^{\varepsilon/250} \mathcal{E}$ where in the definition of $\mathcal{E}$ we take 
		\begin{align*}
			\alpha(v) &= \frac{\mu^2(v) v^{1/2+it_2-it_1} V^{1/2-\varepsilon/500}}{\varphi(v)}, \quad \beta(b_2) := w_{b_2} \frac{B_2^{1/2}}{b_2^{1/2}}, \\ 
			G_1\left(\frac{\ell_1}{L_1}\right) &= G\left(\frac{\ell_1}{L_1}\right) \bfrac{L_1}{\ell_1}^{1/2+1/\log Q+it_2}|\widetilde{V}_2(1/\log Q+it_2)|^{1/3}, \\
			G_2\left(\frac{\ell_2}{L_2}\right) &= G\left(\frac{\ell_2}{L_2}\right) \bfrac{L_2}{\ell_2}^{1/2+1/\log Q+it_2}|\widetilde{V}_2(1/\log Q+it_2)|^{1/3}, \\
			G_3\left(\frac{w}{W}\right) &= G\left(\frac{w}{W}\right) \bfrac{w}{W}^{1/2+it_2-it_1} |\widetilde{V}_2(1/\log Q+it_2)|^{1/3} |\widetilde{\Phi}(1/\log Q+it_1)|;
		\end{align*}
		by the fast decay of the Mellin transforms $G_j^{(\nu)}(x) \ll_\nu 1$ for all $j \in \{1, 2, 3\}$ and all $\nu \geq 0$.
	\end{proof}

	Now it suffices to show that always $\mathcal{E} = O(Q^{2-\varepsilon/10})$. We split into three cases according to the sizes of the parameters.
	
	\textbf{Case 1: $L_2 > L_1 Q^{-\varepsilon}$ or $L_1 B_0 B_2 \leq A_1 Q^{\varepsilon}$.}
	
	At first we only consider the sum over $\ell_2$, which already provides enough cancellation in this case. We remove the condition $(\ell_2, v)=1$ using M\"obius inversion (introducing a variable $d$ and replacing $\ell_2$ by $d\ell_2$) and then split the variable $\ell_2$ into residue classes $\mod{w}$. This gives, for $(v, w) = 1$,
	\begin{align*}
		&\sum_{\substack{\ell_2 \\ (\ell_2, vw) = 1}} G_2\left(\frac{\ell_2}{L_2}\right) e\left(\frac{a_0 \ell_2 \overline{\ell_1 b_0 b_2 v}}{w}\right) = \sum_{\substack{d \mid v}} \mu(d) \sum_{\substack{\ell_2 \\ (\ell_2, w) = 1}} G_2\left(\frac{d \ell_2}{L_2}\right) e\left(\frac{a_0 d \ell_2 \overline{\ell_1 b_0 b_2 v}}{w}\right) \\
		&= \sum_{d \mid v} \mu(d) \sum_{\substack{a \mod{w} \\ (a, w) = 1}} e\left(\frac{a_0 d a \overline{\ell_1 b_0 b_2 v}}{w}\right)  \sum_{\substack{\ell_2 \equiv a \mod{w}}} G_2\left(\frac{d \ell_2}{L_2}\right). 
	\end{align*}
	Applying Poisson summation, this equals
	\begin{align} \label{eq:Case1AfterPoisson}
		\sum_{d \mid v} \mu(d) \sum_{\substack{a \mod{w} \\ (a, w) = 1}} e\left(\frac{a_0 d a \overline{\ell_1 b_0 b_2 v}}{w}\right)  \frac{L_2}{dw} \sum_{|h| \leq \frac{Wd}{L_2} Q^\varepsilon} \widehat{G_2}\left(\frac{hL_2}{dw}\right) e\left(\frac{ah}{w}\right) + O(Q^{-100}).
	\end{align}
	The contribution from $h=0$ equals
	\begin{align*}
		\frac{L_2}{w} \sum_{d \mid v} \frac{\mu(d)}{d} \sum_{\substack{a \mod{w} \\ (a, w) = 1}} e\left(\frac{a_0 d a \overline{\ell_1 b_0 b_2 v}}{w}\right) \widehat{G}(0).
	\end{align*}
	Since $(a_0d, w) = 1$, the Ramanujan sum over $a$ evaluates to $\mu(w)$ (see e.g.~\cite[formula (3.4)]{IK}) and thus the contribution from $h = 0$ equals
	\begin{align*}
		\frac{L_2}{w}\cdot \frac{\varphi(v)}{v} \mu(w) \widehat{G_2}(0),
	\end{align*}
	so, using~\eqref{eq:VariableRanges}, this term contributes to $\mathcal{E}$
	\begin{align*}
		&\ll \left(\frac{W}{A_0 B_0 B_2 L_1 L_2 V} \right)^{1/2} \sum_{\substack{v \sim V, w \sim W}} \sum_{\substack{b_2 \sim B_2}} \sum_{\substack{a_0 \sim A_0, b_0 \sim B_0}} |z'_{a_0, b_0}| \sum_{\substack{\ell_1\sim L_1}} \frac{L_2}{W} \\
		&\ll (VW)^{1/2} B_2^{1/2}(A_0B_0)^{1/2+\epsilon}L_1^{1/2}L_2^{1/2}  \ll Q^{1+\theta+\epsilon}\ll Q^{1+\theta+\epsilon},
	\end{align*}
	which is acceptable as $\theta < 1/2$.
	
	The contribution of the remaining $h$ to~\eqref{eq:Case1AfterPoisson} is
	\begin{align*}
		\frac{L_2}{w}  \sum_{d \mid v} \frac{\mu(d)}{d}  \sum_{0 < |h| \leq \frac{Wd}{L_2} Q^\varepsilon} \widehat{G_2}\left(\frac{hL_2}{dw}\right) \sum_{\substack{a \mod{w} \\ (a, w) = 1}} e\left(\frac{a(a_0 d \overline{\ell_1 b_0 b_2 v}+h)}{w}\right) .
	\end{align*}
	Here
	\begin{align*}
		\sum_{\substack{a \mod{w} \\ (a, w) = 1}} e\left(\frac{a(a_0 d \overline{\ell_1 b_0 b_2 v}+h)}{w}\right) &= \sum_{\substack{a \mod{w} \\ (a, w) = 1}} e\left(\frac{a(a_0 d +h \ell_1 b_0 b_2 v)}{w}\right).
	\end{align*}
	By a bound for Ramanujan sums (see for example~\cite[(3.3)]{IK}) the sum over $a$ has absolute value at most $(w, a_0 d+h\ell_1 b_0 b_2 v)$. Now $\ell_1 b_0 b_2 \nmid a_0$ and $d \mid v,$ so $\ell_1 b_0 b_2 v \nmid a_0 d$ and thus always $a_0 d +h \ell_1 b_0 b_2 v \neq 0$.
	
	Given this, we can bound the contribution of $h \neq 0$ to $\mathcal{E}$ by
	\begin{align*}
		&\ll \left(\frac{W}{A_0 B_0 B_2 L_1 L_2 V} \right)^{1/2} \sum_{\substack{v \sim V, w \sim W}} \sum_{d \mid v} \sum_{\substack{b_2 \sim B_2}} \sum_{\substack{a_0 \sim A_0, b_0 \sim B_0}}  |z'_{a_0, b_0}| \\
		&\qquad \cdot \frac{L_2}{d W}  \sum_{\substack{\ell_1\sim L_1}} \sum_{\substack{0 < |h| \leq \frac{Wd}{L_2}Q^\varepsilon \\ a_0d + h\ell_1 b_0 b_2 v \neq 0}} (w, a_0d + h\ell_1 b_0 b_2 v) \\
		&\ll \frac{Q^{2\epsilon} (VW)^{1/2} B_2^{1/2}(A_0B_0)^{1/2+\epsilon}L_1^{1/2} }{L_2^{1/2}} \max_{1 \leq k \leq q^{ 4+2\epsilon}} \sum_{w \sim W} (w, k) \\
		&\ll \frac{Q^{3\epsilon} V^{1/2} W^{3/2} B_2^{1/2}(A_0B_0)^{1/2+\epsilon}L_1^{1/2} }{L_2^{1/2}} \ll Q^{3/2+4\epsilon}\bfrac{B_2A_0B_0 L_1}{L_2}^{1/2}.
	\end{align*}
	If either $L_2 > L_1 Q^{-\varepsilon}$ or $L_1 B_0 B_2 \leq A_0 Q^\varepsilon$, this is $\ll Q^{ 3/2+\theta+5\varepsilon}$ which is $O(Q^{2-\epsilon})$ for $\theta<1/2-6\varepsilon$ and thus the claim holds in case 1.

	\textbf{Case 2: $V \geq Q^\varepsilon$, $L_2 \leq L_1 Q^{-\varepsilon},$ and $L_1 B_0 B_2 > A_0 Q^{\varepsilon}$.} 
	
	Notice first that since $L_1 B_0 B_2 \geq A_0 Q^\varepsilon$, we can drop the condition $\ell_1 b_0 b_2 \nmid a_0$. In this case we move back to multiplicative characters. Using \eqref{eq:orthogonality of all characters with root number} we have
	\begin{align*}
		e\left(\frac{a_0 \ell_2 \overline{\ell_1 b_0 b_2 v}}{w}\right) = \frac{\sqrt{w}}{\varphi(w)} \sum_{\chi \mod{w}} \overline{\varepsilon_{\chi}} \chi(a_0 \ell_2) \overline{\chi}(\ell_1 b_0 b_2 v),
	\end{align*}
	so that in this case
	\begin{align*}
		\mathcal{E} &\ll \frac{\log Q}{(A_0 B_0 B_2 L_1 L_2 V)^{1/2}} \sum_{\substack{v \sim V \\ w \sim W}} \sum_{\chi \mod{w}} \left|\sum_{\ell_1, \ell_2 \geq 1} \chi(\ell_2) \overline{\chi}(\ell_1) G_1\left(\frac{\ell_1}{L_1}\right) G_2\left(\frac{\ell_2}{L_2}\right) \right| \\
		& \qquad \cdot \left|\sum_{\substack{a_0 \sim A_0, b_0 \sim B_0 \\ b_2 \sim B_2}} \beta(b_2)z'_{a_0,b_0} \chi(a_0) \overline{\chi}(b_0 b_2)  \right|.
	\end{align*}
	Applying the Cauchy-Schwarz inequality and Lemma \ref{le:general orthogonality}, we obtain
	\begin{align*}
		\mathcal{E} &\ll \frac{Q^{\epsilon/10}}{(A_0 B_0 B_2 L_1 L_2 V)^{1/2}} \sum_{\substack{v \sim V \\ w \sim W}} (L_1 L_2 + W)^{1/2} (L_1 L_2)^{1/2} (A_0 B_0 B_2 + W)^{1/2} (A_0 B_0 B_2)^{1/2} \\
		&\ll  Q^{\epsilon/10} V^{1/2} W \left(L_1L_2+W\right)^{1/2}\left(A_0B_0B_2+W\right)^{1/2} \\
		& \ll \frac{ Q^{1+\epsilon/10}}{V^{1/2}} Q^{1/2+\epsilon/20} \lz Q^{\theta}+\frac{Q^{1/2}}{V^{1/2}}\pz
	\end{align*}
	which is $O(Q^{2-\epsilon/2})$ when $\theta<1/2-\epsilon$ and $V \geq Q^\varepsilon$.
	
	\textbf{Case 3: $V < Q^\varepsilon$, $L_2 \leq L_1 Q^{-\varepsilon}$ and $L_1 B_0 B_2 \geq Q^{\varepsilon} A_0$.}
	
	Notice that in this case in particular $L_2 \leq Q^{1/2}$ and we automatically have $\ell_1 b_0 b_2 \nmid a_0$, so we can forget about that condition and we need to estimate
	\begin{align*}
		\mathcal{E} &=  \left(\frac{W}{A_0 B_0 B_2 L_1 L_2 V} \right)^{1/2} \Biggl| \sum_{\substack{v \sim V, w \in \mathbb{N} \\ (v, w)=1}} \alpha(v) G_3\left(\frac{w}{W}\right) \sum_{\substack{b_2 \sim B_2 \\ (b_2, vw) = 1}} \beta(b_2) \\
		&\cdot \sum_{\substack{a_0 \sim A_0, b_0 \sim B_0 \\ (a_0 b_0, vw) = 1}} z'_{a_0, b_0} \sum_{\substack{\ell_1, \ell_2 \geq 1 \\ (\ell_1 \ell_2, vw) = 1}} G_1\left(\frac{\ell_1}{L_1}\right) G_2\left(\frac{\ell_2}{L_2}\right) e\left(\frac{a_0 \ell_2 \overline{\ell_1 b_0 b_2 v}}{w}\right)\Biggr|.
	\end{align*}
	In order to apply the Kloosterman sum bound (Lemma \ref{le:Klo}), we need to remove some cross-conditions between the variables, in particular co-primality conditions concerning the variables $\ell_1$ and $w$ that we will use as smooth variables. We start by removing the condition $(v, \ell_1) = 1$ introducing $\mu(t)$. Additionally we combine $r = a_0 \ell_2$. This yields
	\begin{align*}
		\mathcal{E} &=  \left(\frac{W}{A_0 B_0 B_2 L_1 L_2 V} \right)^{1/2} \Biggl|\sum_{t \geq 1} \mu(t) \sum_{\substack{v \sim V/t, w\geq 1 \\ (tv, w)=1}} \alpha(vt) G_3\left(\frac{w}{W}\right) \sum_{\substack{b_2 \sim B_2 \\ (b_2, vtw) = 1}} \beta(b_2) \\
		&\cdot \qquad \sum_{\substack{r \asymp A_0 L_2, \, b_0 \sim B_0 \\ (r b_0, vtw) = 1}} y_{r, b_0} \sum_{\substack{\ell_1\geq 1 \\ (\ell_1, w) = 1}} G_1\left(\frac{\ell_1t}{L_1}\right) e\left(\frac{r \overline{\ell_1 t^2 b_0 b_2 v}}{w}\right)\Biggr|,
	\end{align*}
	where
	$$
	y_{r,b_0}=\sum_{\substack{r=a_0\ell_2\\ a_0\sim A_0}}z'_{a_0,b_0} G_2\bfrac{\ell_2}{L_2}.
	$$
	Next we remove the condition $(r, w) = 1$ introducing $\mu(f)$, so that
	\begin{align*}
		\mathcal{E} &=  \left(\frac{W}{A_0 B_0 B_2 L_1 L_2 V} \right)^{1/2} \Biggl|\sum_{f, t \geq 1} \mu(f) \mu(t) \sum_{\substack{v \sim V/t, w \geq 1 \\ (tv, fw)=1}} \alpha(vt) G_3\left(\frac{wf}{W}\right) \\
		&\quad \cdot\sum_{\substack{b_2 \sim B_2 \\ (b_2, vtfw) = 1}} \beta(b_2)  \sum_{\substack{r \asymp A_0 L_2/f,\ b_0 \sim B_0 \\ (b_0, vtfw) = 1 \\ (rf, vt) = 1}} y_{rf, b_0} \sum_{\substack{\ell_1 \geq 1 \\ (\ell_1, wf) = 1}} G_1\left(\frac{\ell_1t}{L_1}\right) e\left(\frac{r \overline{\ell_1 t^2 b_0 b_2 v}}{w}\right)\Biggr|.
	\end{align*}
	Finally we need to remove the condition $(\ell_1, f) = 1$ introducing $\mu(h)$, so that
	\begin{align*}
		\mathcal{E} &=  \left(\frac{W}{A_0 B_0 B_2 L_1 L_2 V} \right)^{1/2} \Biggl|\sum_{f, h, t \geq 1} \mu(f) \mu(h) \mu(t) \sum_{\substack{v \sim V/t,\ w\geq 1 \\ (tv, fhw)=1}} \alpha(vt) G_3\left(\frac{wfh}{W}\right) \\
		&\quad \cdot \sum_{\substack{b_2 \sim B_2 \\ (b_2, vtfhw) = 1}} \beta(b_2) \sum_{\substack{r \asymp A_0 L_2/fh,\ b_0\sim B_0 \\ (b_0, vtfhw) = 1 \\ (r, vt) = 1}} y_{rfh, b_0} \sum_{\substack{\ell_1 \geq 1 \\ (\ell_1, w) = 1}} G_1\left(\frac{\ell_1ht}{L_1}\right) e\left(\frac{r \overline{\ell_1 h t^2 b_0 b_2 v}}{w}\right)\Biggr|.
	\end{align*}
	
	We split the variables $f, h, t$ dyadically, so that $f \sim F$, $h \sim H$, and $t \sim T$. Note that we can assume that $T \ll V \leq Q^\varepsilon$. In order to easily apply Lemma~\ref{le:Klo}, we note that $\mathcal{E} \ll (\log Q)^{O(1)} \mathcal{E}'$, where
	\[
	\mathcal{E}' := \left(\frac{W}{A_0 B_0 B_2 L_1 L_2 V} \right)^{1/2} F \Biggl| \sum_{\widetilde{r} \asymp \widetilde{R}} \sum_{\widetilde{n} \asymp \widetilde{N}} \widetilde{b}_{\widetilde{n}, \widetilde{r}} \sum_{\substack{\widetilde{c}, \widetilde{d} \\ (\widetilde{c}, \widetilde{r}\widetilde{d}) = 1}} \widetilde{g}\left(\frac{\widetilde{c}}{\widetilde{C}}, \frac{\widetilde{d}}{\widetilde{D}}\right) e\left(\widetilde{n} \frac{\overline{\widetilde{r} \widetilde{d}}}{\widetilde{c}}\right)\Biggr|,
	\]
	where $\widetilde{g}$ is a smooth function with compact support, we substituted
	\begin{align*}
		\widetilde{n} = r, \quad \widetilde{r} = ht^2b_0 b_2 v, \quad \widetilde{d}= \ell_1, \quad \text{and} \quad \widetilde{c} = w,
	\end{align*}
	and wrote
	\begin{align*}
		\widetilde{N} = \frac{A_0 L_2}{FH}, \quad \widetilde{R} = HT^2B_0 B_2\frac{V}{T} = HTVB_0 B_2, \quad \widetilde{D}= \frac{L_1}{HT}, \quad \text{and} \quad \widetilde{C} = \frac{W}{FH},
	\end{align*}
	and
	\[
	\widetilde{b}_{\widetilde{n}, \widetilde{r}} = \sum_{\substack{\widetilde n=r \\ \widetilde{r} = ht^2b_0 b_2v \\ h \sim H, t \sim T \\ b_0 \sim B_0, b_2 \sim B_2, v \sim V/T }} z_{r, h, t, b_0, b_2, v, f}
	\]
	for certain coefficients $z_{r, h, t, b_0, b_2, v, f} \ll Q^{\varepsilon/30}$. We have ignored cross-conditions in smooth functions but they can be removed by the standard trick of applying Mellin transforms as in Lemma \ref{le:dyadic decomposition}.
	Now
	\begin{align*}
		\sum_{\widetilde{n}} \sum_{\widetilde{r}} \left| \widetilde{b}_{\widetilde{n}, \widetilde{r}}\right|^2 \ll Q^{\varepsilon/15} \sum_{\substack{\widetilde{n} \asymp \widetilde{N} \\ t^2 u \asymp \widetilde{R} \\ t \sim T}} \tau_4(u) \ll  Q^{\varepsilon/10} \widetilde{N} \frac{\widetilde{R}}{T} \ll \frac{VA_0B_0B_2L_2}{F}
	\end{align*}
	and thus by Lemma~\ref{le:Klo} we obtain that
	\begin{align*}
		\mathcal{E}^2 &\ll Q^{\varepsilon/5} \frac{F W}{ L_1} \Biggl(\frac{W}{FH}\left(HTVB_0 B_2+\frac{A_0 L_2}{FH}\right)\left(\frac{W}{FH}+VB_0 B_2L_1\right) \\
		& \qquad + \frac{W^2 L_1}{F^2 H^3 T} \sqrt{\left(HTVB_0 B_2 + \frac{A_0 L_2}{FH}\right) HTVB_0 B_2} + \frac{L_1^2A_0 B_0 B_2 VL_2}{TFH^2}\Biggr).
	\end{align*}
	Using $T \ll V \leq Q^\varepsilon$, ignoring some denominators, and multiplying out, we obtain
	\begin{align*}
		\mathcal{E}^2 &\ll Q^{\varepsilon/5} \Biggl(W^3 V^2 B_0 B_2+ W^2V^3 (B_0 B_2)^2 +  W^3 A_0 \frac{L_2}{L_1} +  W^2 V L_2 A_0B_0B_2 \\
		& \qquad +  W^3 V B_0 B_2 + W^3 \sqrt{L_2A_0B_0 B_2 V} +  W V L_1L_2A_0B_0 B_2\Biggr).
	\end{align*}
	Recalling our assumptions in Case 3, it is straightforward to see that $\mathcal{E}^2 \ll Q^{4-\varepsilon}$ for $\theta < 1/2-10\varepsilon$ as requested.
	
	\section{The unbalanced case} \label{se:unbalanced}
	We end with a brief discussion of the optimality of the unbalanced Michel-Vanderkam mollifier.
	Consider the unbalanced mollifier
	\begin{equation}\label{eq:M unbalanced} \begin{aligned}
			M(\chi) &= \sum_{b\leq q^{\theta_1}}\frac{\mu(b)\chi(b)}{\sqrt b}\left(1-\frac{\log b}{\log q^{\theta_1}}\right)+\alpha \overline\epsilon_\chi\sum_{b \leq q^{\theta_2}}\frac{\mu(b)\overline\chi(b)}{\sqrt b}\left(1-\frac{\log b}{\log q^{\theta_2}}\right) \\
			&=: M_1(\chi) + \alpha M_2(\chi),
		\end{aligned}
	\end{equation}
	say, with $\theta_2 \leq \theta_1<1/2$, which gives the best known non-vanishing proportion for fixed prime moduli (see~\cite{Khan-Mil-Ngo}). To study this case, we need the following moment estimates.
	\begin{proposition}\label{prop:PsiM,PsiMM unbalanced}
		Let $q\geq 3$, and let $M_1(\chi)$ and $M_2(\chi)$ be as in \eqref{eq:M unbalanced}. Then the following hold.
		\begin{enumerate}[(i)]
			\item For $j=1,2$, we have $\Psi_{M_j}(q)=\phi^+(q)(1+o(1))$.
			\item For $j=1,2$, we have $\Psi_{M_j,M_j}(q)=\phi^+(q)\lz1+\frac1{\theta_j}+o(1)\pz$.
			\item Assume also that $q$ is a prime, and $2\theta_1+\theta_2<1$ and $4\theta_1+6\theta_2<3$. Then
			$\Psi_{M_1,M_2}(q)=\phi^+(q)(1+o(1)).$
		\end{enumerate}
	\end{proposition}
	\begin{proof}
		Part (i) follows from Proposition \ref{prop:Psi_N_1} with $\theta=\theta_1$ and $z_{a,b}=\mathbf{1}_{a=1} \mathbf{1}_{b\leq q^{\theta_j}}\mu(b)(1-\frac{\log b}{\log q^{\theta_j}})$, whereas part (ii) follows from Proposition~\ref{prop:PsiMN_1IS}(ii), noting that $\Psi_{M_2, M_2}(q) = \Psi_{\varepsilon_\chi M_2, \varepsilon_\chi M_2}(q)$.
		
		For Part (iii), we refer to \cite{Khan-Mil-Ngo}. In particular, our $\frac{\Psi_{M_1,M_2}}{2\phi^+(q)}$ equals to the term \cite[(3.4)]{Khan-Mil-Ngo} with $c_1=c_2=1$, where $M$ from \cite[(3.4)]{Khan-Mil-Ngo} corresponds to $q^{\theta_2}$, and $MR$ from \cite[(3.4)]{Khan-Mil-Ngo} corresponds to our $q^{\theta_1}$. The main term of \cite[(3.4)]{Khan-Mil-Ngo} is evaluated in \cite[below (3.8)]{Khan-Mil-Ngo}. Now with $\theta,\alpha$ as in \cite[below (3.12)]{Khan-Mil-Ngo}, we have $\theta_1=\theta+\alpha$, $\theta_2=\theta$, so our assumptions imply that the lengths of our mollifiers satisfy \cite[(3.13)]{Khan-Mil-Ngo}, which are the conditions under which the error terms can be proved to be negligible.
	\end{proof}
	We remark that the main term in $\Psi_{M_1,M_2}(q)$ could be obtained from Proposition \ref{prop:PsiM_2N_1}, but we only estimated the error terms there on average (in a larger range of $\theta_j$).
	
	Let us first use Theorem \ref{thm:criterion2} to find the optimal value of $\alpha$. Inserting the moment calculations from Proposition~\ref{prop:PsiM,PsiMM unbalanced} there, we obtain that the optimal value of $\alpha$ is
	\[
	\alpha = \frac{1-(1+\frac{1}{\theta_1})}{1-(1+\frac{1}{\theta_2})} + o(1) = \frac{\theta_2}{\theta_1} + o(1),
	\]
	which is in agreement with the result of Khan-Mili\'cevi\'c-Ngo \cite{Khan-Mil-Ngo}.
	
	\smallskip
	
	We now assume that $\theta_2<\theta_1$, set $\alpha=\frac{\theta_2}{\theta_1}$ and proceed to investigate whether adding a Bui-type piece would be useful in this context. We define
	\[
	N(\chi) := \sum_{ab\leq q^{\theta_1-\varepsilon_1}}\frac{x_{a, b} \overline{\chi}(a) \chi(b)}{\sqrt{ab}},
	\]
	with a small $\varepsilon_1 \in (0, \theta_1-\theta_2)$. It will show up that with adequately selected $x_{a, b} \in \mathbb{R}$ with $x_{a, b} \ll \exp\lz c\frac{(\log q)^{3/5}}{(\log \log q)^{1/5}}\pz$ for a small constant $c > 0$, the non-vanishing proportion can be improved (without increasing the lengths of the mollifiers) if one is able to compute the moments and the main terms are as expected.
	
	To show this, we use Corollary~\ref{cor:EfficientN}(i), which implies that if
	\begin{equation} \label{eq:unbalancedclaim}
		|\overline{\Psi_M(q)} \Psi_{M, N}(q) - \overline{\Psi_N(q)} \Psi_{M, M}(q)| \gg |\Psi_N(q)| \Psi_{M, M}(q),
	\end{equation} and $N(\chi)$ is efficient, then adding a multiple of $N(\chi)$ is useful.
	
	By Proposition \ref{prop:PsiM,PsiMM unbalanced}, we find that
	\[
	\begin{aligned}
		\Psi_M(q) &= \Psi_{M_1}(q)+\alpha\Psi_{M_2}(q)=\varphi^+(q)(1+\alpha+o(1)), \\
		\Psi_{M, M}(q)& = \Psi_{M_1,M_1}(q)+2\alpha\Psi_{M_1,M_2}(q)+\alpha^2\Psi_{M_2,M_2}(q)\\
		&=\phi^+(q)\lz(1+\alpha)^2+\frac{1}{\theta_1}+\frac{\alpha^2}{\theta_2}\pz.
	\end{aligned}
	\]
	Furthermore, by Propositions~\ref{prop:Psi_N_1} and~\ref{prop:PsiM1,N1}, for any $\varepsilon > 0$,
	\[
	\Psi_N(q) = \phi^+(q) \sum_{\substack{kb^2\leq q^{\theta_1-\varepsilon_1} \\ (kb, q) = 1}}\frac{x_{kb,b}}{kb} + O(q^{\theta_1/2+1/2+\varepsilon}),
	\]
	and
	\[
	\Psi_{M_1, N}(q) = \phi^+(q) \sum_{\substack{kb^2\leq q^{\theta_1-\varepsilon_1} \\ (kb, q) = 1}}\frac{x_{kb,b}}{kb}\lz1+\frac{1}{\theta_1}-\frac{\log k}{\log q^{\theta_1}}\pz+o\lz \phi^+(q)\pz.
	\]
	
	Now it remains to calculate $\Psi_{M_2,N}(q)$. We cannot directly use Proposition~\ref{prop:PsiM_2N_1}, as the error terms there were only estimated on average over $q$. However, we can follow the proof and see what one would expect.
	From the proof (see in particular~\eqref{eq:PsiN0M2Mformula}) one can see that the main term of $\Psi_{M_2, N}(q)$ is supposed to be
	\begin{align*}
		&\varphi^+(q) \sum_{\substack{b_2 \leq q^{\theta_2} \\ b_2 n \ell_1 b_0^2 \leq q^{\theta_1-\varepsilon_1} \\ (b_0 b_2 \ell_1 n, q) = 1}} \mu(b_2) \left(1-\frac{\log b_2}{\log q^{\theta_2}}\right) \frac{x_{b_0 b_2 \ell_1 n, b_0}}{b_0 b_2 \ell_1 n} \\
		&= \varphi^+(q) \sum_{\substack{b_2 n \ell_1 b_0^2 \leq q^{\theta_1-\varepsilon_1} \\ (b_0 b_2 \ell_1 n, q) = 1}} \mu(b_2) \left(1-\frac{\log b_2}{\log q^{\theta_2}}\right) \frac{x_{b_0 b_2 \ell_1 n, b_0}}{b_0 b_2 \ell_1 n} \\
		& - \varphi^+(q) \sum_{\substack{q^{\theta_2} < b_2 \leq q^{\theta_1-\varepsilon_1} \\ b_2 n \ell_1 b_0^2 \leq q^{\theta_1-\varepsilon_1} \\ (b_0 b_2 \ell_1 n, q) = 1}} \mu(b_2) \left(1-\frac{\log b_2}{\log q^{\theta_2}}\right) \frac{x_{b_0 b_2 \ell_1 n, b_0}}{b_0 b_2 \ell_1 n} =: \Psi_{M_2, N}^{\mathcal{M}_1}(q) + \Psi_{M_2, N}^{\mathcal{M}_2}(q),
	\end{align*}
	say. Now, as in the argument following~\eqref{eq:PsiN0M2Mformula}, we have
	\[
	\Psi_{M_2, N}^{\mathcal{M}_1}(q) = \varphi^+(q) \sum_{\substack{kb^2 \leq q^{\theta_1-\varepsilon_1} \\ (kb, q) = 1}} \frac{x_{kb, b}}{kb} \left(1+\frac{\log k}{\log q^{\theta_2}}\right).
	\]
	Noting that $\alpha\cdot\frac{\log k}{\log q^{\theta_2}}=\frac{\log k}{\log q^{\theta_1}}$, we obtain
	\[
	\Psi_{M_1, N}(q) + \alpha\Psi_{M_2, N}^{\mathcal{M}_1}(q) =  \varphi^+(q) \sum_{\substack{k^2 b \leq q^{\theta_1-\varepsilon_1} \\ (bk, q) = 1}} \frac{x_{kb, k}}{kb}\lz1+\alpha+\frac{1}{\theta_1}\pz +o\lz \phi^+(q)\pz.
	\]
	The above computations now imply that the main terms in $\Psi_M(q)(\Psi_{M_1,N}(q)+\alpha\Psi^{\mathcal M_1}_{M_2,N}(q))-\Psi_{N}(q)\Psi_{M,M}(q)$ vanish, so assuming that the error terms work out, the main term of the left hand side of~\eqref{eq:unbalancedclaim} is
	\begin{equation}\label{eq:unbalanced difference}
		\Psi_M(q) \cdot\alpha\Psi_{M_2, N}^{\mathcal{M}_2}(q)= -\alpha(1+\alpha)  \varphi^+(q)^2 \sum_{\substack{q^{\theta_2} < b_2 \leq q^{\theta_1-\varepsilon_1} \\ b_2 n \ell_1 b_0^2 \leq q^{\theta_1-\varepsilon_1} \\ (b_0 b_2 \ell_1 n, q) = 1}} \mu(b_2) \left(1-\frac{\log b_2}{\log q^{\theta_2}}\right) \frac{x_{b_0 b_2 \ell_1 n, b_0}}{b_0 b_2 \ell_1 n}.
	\end{equation}
	Taking for instance $x_{a, b} =\frac{ \Lambda(a)}{\log q} \mu(b) P_2\lz1-\frac{\log ab}{\log q^{\theta_1-\epsilon_1}}\pz$ for some polynomial $P_2(x)$ like Bui (see~\eqref{eq:BuiMollifier}), one can check that for most choices of the polynomial $P_2(x)$ with $P_2(0)=0$, this is not $o(|\Psi_N(q)| \Psi_{M, M}(q))$ and $N(\chi)$ is efficient. Thus by Corollary~\ref{cor:EfficientN}(i) we can see that adding a Bui piece should be useful in the unbalanced case. 
	
	We also remark that the main term in \eqref{eq:unbalanced difference} vanishes for instance if $x_{a,b}$ is supported on $a< q^{\theta_2}$, which for example implies that it is not helpful to add any mollifier of the form
	\begin{equation*}
		\sum_{b\leq q^{\theta_1-\varepsilon_1}}\frac{x_b\chi(b)}{\sqrt b}
	\end{equation*} to \eqref{eq:M unbalanced}.
	
	\bibliography{../refs}

\begin{thebibliography}{10}

\bibitem{Bombieri}
Enrico Bombieri.
\newblock A lower bound for the zeros of {R}iemann's zeta function on the
  critical line (following {N}. {L}evinson).
\newblock In {\em S\'eminaire {B}ourbaki (1974/1975: {E}xpos\'es {N}os.
  453--470)}, volume Vol. 514 of {\em Lecture Notes in Math.}, pages Exp. No.
  465, pp. 176--182. Springer, Berlin-New York, 1976.

\bibitem{Bui}
H.~M. Bui.
\newblock Non-vanishing of {D}irichlet {$L$}-functions at the central point.
\newblock {\em Int. J. Number Theory}, 8(8):1855--1881, 2012.

\bibitem{BCY}
H.~M. Bui, Brian Conrey, and Matthew~P. Young.
\newblock More than 41\% of the zeros of the zeta function are on the critical
  line.
\newblock {\em Acta Arith.}, 150(1):35--64, 2011.

\bibitem{Ch}
S.~Chowla.
\newblock {\em The {R}iemann hypothesis and {H}ilbert's tenth problem}.
\newblock Mathematics and its Applications, Vol. 4. Gordon and Breach Science
  Publishers, New York-London-Paris, 1965.

\bibitem{Conrey}
Brian Conrey.
\newblock Zeros of derivatives of {R}iemann's {$\xi $}-function on the critical
  line.
\newblock {\em J. Number Theory}, 16(1):49--74, 1983.

\bibitem{CIS}
J.~Brian Conrey, Henryk Iwaniec, and Kannan Soundararajan.
\newblock Critical zeros of {D}irichlet {$L$}-functions.
\newblock {\em J. Reine Angew. Math.}, 681:175--198, 2013.

\bibitem{DI}
J.-M. Deshouillers and H.~Iwaniec.
\newblock Kloosterman sums and {F}ourier coefficients of cusp forms.
\newblock {\em Invent. Math.}, 70(2):219--288, 1982/83.

\bibitem{DPR}
Sary Drappeau, Kyle Pratt, and Maksym Radziwi\l\l.
\newblock One-level density estimates for {D}irichlet {$L$}-functions with
  extended support.
\newblock {\em Algebra Number Theory}, 17(4):805--830, 2023.

\bibitem{Fouvry}
\'Etienne Fouvry.
\newblock Sur le probl\`eme des diviseurs de {T}itchmarsh.
\newblock {\em J. Reine Angew. Math.}, 357:51--76, 1985.

\bibitem{IvicBook}
Aleksandar Ivi\'c.
\newblock {\em The {R}iemann zeta-function}.
\newblock Dover Publications, Inc., Mineola, NY, 2003.
\newblock Theory and applications, Reprint of the 1985 original [Wiley, New
  York; MR0792089 (87d:11062)].

\bibitem{IS}
H.~Iwaniec and P.~Sarnak.
\newblock Dirichlet {$L$}-functions at the central point.
\newblock In {\em Number theory in progress, {V}ol. 2
  ({Z}akopane-{K}o\'{s}cielisko, 1997)}, pages 941--952. de Gruyter, Berlin,
  1999.

\bibitem{IK}
Henryk Iwaniec and Emmanuel Kowalski.
\newblock {\em Analytic number theory}, volume~53 of {\em American Mathematical
  Society Colloquium Publications}.
\newblock American Mathematical Society, Providence, RI, 2004.

\bibitem{Khan-Mil-Ngo}
Rizwanur Khan, Djordje Mili\'cevi\'c, and Hieu~T. Ngo.
\newblock Nonvanishing of {D}irichlet {$L$}-functions, {II}.
\newblock {\em Math. Z.}, 300(2):1603--1613, 2022.

\bibitem{Khan-Ngo}
Rizwanur Khan and Hieu~T. Ngo.
\newblock Nonvanishing of {D}irichlet {$L$}-functions.
\newblock {\em Algebra Number Theory}, 10(10):2081--2091, 2016.

\bibitem{Kou}
Dimitris Koukoulopoulos.
\newblock {\em The distribution of prime numbers}, volume 203 of {\em Graduate
  Studies in Mathematics}.
\newblock American Mathematical Society, Providence, RI, [2019] \copyright
  2019.

\bibitem{KM2000}
E.~Kowalski and P.~Michel.
\newblock A lower bound for the rank of {$J_0(q)$}.
\newblock {\em Acta Arith.}, 94(4):303--343, 2000.

\bibitem{Levinson}
Norman Levinson.
\newblock More than one third of zeros of {R}iemann's zeta-function are on
  {$\sigma =1/2$}.
\newblock {\em Advances in Math.}, 13:383--436, 1974.

\bibitem{M-V}
Philippe Michel and Jeffrey VanderKam.
\newblock Non-vanishing of high derivatives of {D}irichlet {$L$}-functions at
  the central point.
\newblock {\em J. Number Theory}, 81(1):130--148, 2000.

\bibitem{Pascadi}
Alexandru Pascadi.
\newblock Large sieve inequalities for exceptional maass forms and
  applications.
\newblock {\em Pre-print, arXiv:2404.04239}, 2024.

\bibitem{Pratt}
Kyle Pratt.
\newblock Average nonvanishing of {D}irichlet {$L$}-functions at the central
  point.
\newblock {\em Algebra Number Theory}, 13(1):227--249, 2019, see also
  arXiv:1804.01445.

\bibitem{PRZZ}
Kyle Pratt, Nicolas Robles, Alexandru Zaharescu, and Dirk Zeindler.
\newblock More than five-twelfths of the zeros of {$\zeta$} are on the critical
  line.
\newblock {\em Res. Math. Sci.}, 7(2):Paper No. 2, 74, 2020.

\bibitem{Radziwill}
Maksym Radziwill.
\newblock Limitations to mollifying $\zeta(s)$.
\newblock {\em Pre-print, arXiv:1207.6583}, 2012.

\bibitem{Selberg1}
Atle Selberg.
\newblock Contributions to the theory of the {R}iemann zeta-function.
\newblock {\em Arch. Math. Naturvid.}, 48(5):89--155, 1946; Collected papers,
  vol. I, 214--280, Springer, Berlin 1989.

\bibitem{Selberg2}
Atle Selberg.
\newblock Contributions to the theory of {D}irichlet's {$L$}-functions.
\newblock {\em Skr. Norske Vid.-Akad. Oslo I}, 1946(3):62, 1946; Collected
  papers, vol. I, 281--340, Springer, Berlin 1989.

\bibitem{Selberg3}
Atle Selberg.
\newblock The zeta-function and the {R}iemann hypothesis.
\newblock In {\em C. {R}. {D}ixi\`eme {C}ongr\`es {M}ath. {S}candinaves 1946},
  pages 187--200. Gjellerup, Copenhagen, 1947; Collected papers, vol. I,
  341--355, Springer, Berlin 1989.

\bibitem{Sound95}
K.~Soundararajan.
\newblock Mean-values of the {R}iemann zeta-function.
\newblock {\em Mathematika}, 42(1):158--174, 1995.

\bibitem{sound}
K.~Soundararajan.
\newblock Nonvanishing of quadratic {D}irichlet {$L$}-functions at
  {$s=\frac12$}.
\newblock {\em Ann. of Math. (2)}, 152(2):447--488, 2000.

\bibitem{Titchmarsh}
E.~C. Titchmarsh.
\newblock {\em The theory of the {R}iemann zeta-function}.
\newblock The Clarendon Press, Oxford University Press, New York, second
  edition, 1986.
\newblock Edited and with a preface by D. R. Heath-Brown.

\end{thebibliography}
	\bibliographystyle{plain}
	
\end{document}